\documentclass[preprint, 11pt]{elsarticle}
\usepackage{lineno}
\pdfoutput=1

\usepackage{amsmath}
\usepackage{amssymb}
\usepackage{amsfonts}
\usepackage{mathtools}
\usepackage{mathrsfs}
\usepackage{color}
\usepackage[dvipsnames]{xcolor}
\usepackage{graphicx}
\usepackage[percent]{overpic}
\usepackage{url}
\usepackage{caption}
\usepackage{subcaption}
\usepackage{enumerate}
\usepackage{overpic}
\usepackage{amsthm}
\usepackage{moreverb}
\usepackage{multirow}
\usepackage{multicol}
\usepackage[pdftex,colorlinks,bookmarksopen,bookmarksnumbered,citecolor=red,urlcolor=red]{hyperref}
\usepackage{mwe}
\usepackage{graphbox}
\usepackage{textcomp}
\usepackage{physics}
\usepackage{tikz}

\def\0{\boldsymbol{0}}

\def\xbf{\boldsymbol{x}}

\def\calcf{\mathbf{\mathcal{C}}}
\def\calsf{\mathbf{\mathcal{S}}}
\def\caldf{\mathbf{\mathcal{D}}}
\def\caluf{\mathbf{\mathcal{U}}}
\def\calsigf{\mathbf{\bm{\Sigma}}}
\def\calvf{\mathbf{\mathcal{V}}}
\def\calif{\mathbf{\mathcal{I}}}
\def\calof{\mathbf{\mathcal{O}}}
\def\calff{\mathbf{\mathcal{F}}}
\def\calmf{\mathbf{\mathcal{M}}}

\def\tildeq{\widetilde{q}}
\def\tildep{\widetilde{\psi}}

\def\cl {\nonumber \\}
\def\el {\nonumber }

\newtheorem{rem}{Remark}[section]

\newcommand{\bm}[1]{\mbox{\boldmath{$#1$}}}
\definecolor{ForestGreen}{RGB}{34,139,34}
\DeclarePairedDelimiterX{\inp}[2]{\langle}{\rangle}{#1, #2}
\def\div{\nabla\cdot}

\usepackage{verbatim}
\usepackage{graphicx}
\usepackage{epstopdf}

\begin{document}

\begin{frontmatter}


\title{Data-driven reduced order modeling of
a two-layer quasi-geostrophic ocean model}





\author[Houston]{Lander Besabe}
\ead{lybesabe@central.uh.edu}

\author[SISSA]{Michele Girfoglio\corref{mycorrespondingauthor}}
\ead{mgirfogl@sissa.it}
\cortext[mycorrespondingauthor]{Corresponding author}

\author[Houston]{Annalisa Quaini}
\ead{aquaini@central.uh.edu}

\author[SISSA]{Gianluigi Rozza}
\ead{grozza@sissa.it}

\address[SISSA]{SISSA, International School for Advanced Studies, Mathematics Area, mathLab, via Bonomea, Trieste 265 34136, Italy}
\address[Houston]{Department of Mathematics, University of Houston, 3551 Cullen Blvd, Houston TX 77204, USA}

\begin{abstract}
The two-layer quasi-geostrophic equations (2QGE) is a simplified model that describes the dynamics of a stratified, wind-driven ocean in terms of potential vorticity and stream function. Its numerical simulation  is plagued by a high computational cost due to the size
of the typical computational domain and the need for high resolution to capture the full spectrum of turbulent scales.
In this paper, we present a data-driven reduced order model (ROM) for the 2QGE that drastically reduces the computational time to predict ocean dynamics, especially when there are variable physical parameters.
The main building blocks of our ROM are: i) proper orthogonal decomposition (POD) and ii) long short-term memory (LSTM) recurrent neural networks. Snapshots data are collected from a high-resolution simulation for part of the time interval of interest and for given parameter values in the case of variable parameters.
POD is applied to each field variable to extract the dominant modes and a LSTM model is trained on the modal coefficients associated with the snapshots for each variable. Then, the trained LSTM models predict the modal coefficients for the remaining part of the time interval of interest and for a new parameter value. 
To illustrate the predictive performance of our POD-LSTM ROM and the corresponding time savings, we consider an extension of the so-called double-gyre wind forcing test. We show that the POD-LSTM ROM is accurate in predicting 
both time-averaged fields
and time-dependent quantities (modal coefficients, enstrophy, and kinetic energy), even when retaining only 10-20\% of the singular value energy
of the system. The computational speed up for the prediction is about up
to 1E+07 compared to a finite volume based full order method.
\end{abstract}


\end{frontmatter}

\section{Introduction}


Ocean flows are notorious for being computationally costly to simulate for several reasons, which include
the size of the domain (surface area of the order of 1E+06~Km$^2$ and depth of up to 11 Km), the wide spectra of spatial scales, and the need for long time integration to observe long term flow behavior. The spectra of 
spatial scales depend on the flow 
regime that is typically characterized by two non-dimensional numbers: the Reynolds number $Re$, which is the ratio between the inertial and the viscous forces; and the Rossby number $Ro$, which is the ratio between the inertial force and the Coriolis force.
When $Re$ is large and $Ro$ is small, 
as is typical for ocean flows, 
the smallest spatial scales are
particularly small and thus can only be resolved by extremely fine computational mesh, which are costly.

In order to reduce the computational cost, assumptions are introduced at
the modeling level.
Since the ocean depth is much smaller
than its horizontal length, one can replace the 
incompressible Navier-Stokes equations
with the shallow water equations, which are 
obtained from the Navier-Stokes equations
by averaging on the vertical coordinate. 
By assuming hydrostatic and geostrophic balance, i.e., very small $Ro$, one can further simplify
the shallow water equations to get the
single-layer quasi-geostrophic equations (QGE), also known as the barotropic vorticity equations. 
While the QGE are a widely used model for ocean 
dynamics, they do not take into account 
the fact that the ocean is a stratified fluid due to differences in temperature and salinity which result into varying densities. 
Thus, a slightly more realistic model comes from 
adding the effect of stratification, which leads to the two-layer quasi-geostrophic equations (2QGE). We note that, although the 2QGE is a simplified model of the ocean, it is still rather complex and computationally expensive 
when $Re$ grows large and $Ro$ grows small.

When the flow depends on parameters (e.g., $Re$, $Ro$, layers depths) known to vary in a given range,  as is the case for ocean flows, 
reduced order modeling (ROM) can be applied
to radically reduce the computational time.
In broad terms, ROMs are inexpensive surrogates for expensive models.
The surrogate is constructed by extracting the dominant system dynamics from selected high-resolution simulations using, e.g., Proper Orthogonal Decomposition (POD) \cite{Weller2010, Giere2016}. 
This computationally intensive initial phase, called offline, is performed once and for all. 
The expense incurred in this construction process is then amortized over many surrogate solutions. 
In this second phase, called online, the cheap
surrogate solutions can be computed in different ways,
e.g.,  projection methods (e.g., \cite{Girfoglio_JCP}), 
interpolation (e.g., \cite{Xiao2015}), or machine learning 
techniques (e.g., \cite{Maulik2020}).
If projection methods are applied, the ROM is called 
intrusive, otherwise it is called non-intrusive.
ROMs are widely used in many query contexts in fields ranging
ranging from structural mechanics \cite{Casenave2019,Guo2018,Lin1995,Oliver2017} to
biomedical \cite{Fresca2021,Fresca2020,Pfaller2020,Lucas2019,Tezzele2018b,BALLARIN2016609} and naval engineering \cite{Tezzele2022,GADALLA2021104819,jmse9020185,Tezzele2018}. 
See, e.g., \cite{peter2021modelvol1,benner2020modelvol2,benner2020modelvol3,hesthaven2016certified, malik2017reduced,rozza2008reduced,Rozza2022}
for reviews on ROM.

Much interesting work has been done on the construction of intrusive and non-intrusive ROMs for the single-layer QGE. One of the earliest works on this topic \cite{San2015} constructs a ROM for the QGE using POD to generate the reduced space for vorticity 
and Galerkin projection for the online phase. 
Large eddy simulation inspired
closure models are utilized to stabilize 
this projection-based ROM. In \cite{Mou-2020}, two methodologies to capture the unresolved modes 
(i.e., the modes discarded during the offline phase)
are introduced, both based on
data-driven corrections. It is shown that these 
corrections significantly
improve the accuracy of the standard POD-Galerkin projection ROM. 
A hybrid ROM approach for the QGE that 
uses a convex combination of Galerkin projection and extreme learning machine models is introduced in \cite{Rahman-2018}. 
In \cite{Rahman-2019} a non-intrusive ROM is constructed by using short-term memory-based (LSTM) architecture trained on the modal coefficients of the high-resolution simulations performed offline. 
A similar approach to \cite{Rahman-2019}, but with
LSTM  replaced with dynamic mode decomposition \cite{kutz2016dynamic}, is explored in \cite{Golzar2023}. 
The only work that applied ROM to a 2QGE model, 
specifically the stochastic 2QGE, is \cite{Qi2016}, 
which levereges a truncated spectral expansion of the baroclinic-barotropic components of the 2QGE 
for the analysis of the statistical response of the model with respect to important physical quantities. 

In this paper, we develop a non-intrusive (also called data-driven) ROM for the 2QGE that uses POD for the generation of a basis for the reduced order space 
in the offline phase
and a simple LSTM architecture to find the coefficients
for the reduced basis functions to obtain
the ROM approximation of each variable in the online phase.
Our approach could be seen as an extension of the methodology in 
\cite{Rahman-2019} to the 2QGE. However, unlike \cite{Rahman-2019}, 
where POD is applied only to the snapshots
matrix for the vorticity, we apply POD to generate different
bases for vorticity and stream function. We show that this leads
to improved accuracy.
For a classical benchmark that displays a very large range of fluctuating 
spatial and temporal scales, we show that our ROM is not only accurate but
extremely efficient, with up to 1E+07 speed up factor in the predictive phase.

The rest of the paper is organized as follows. Sec.~\ref{sec:2qge} describes the 2QGE model. Sec.~\ref{sec:fom} presents the time and space discretization of the model and the segregated algorithm used to obtain the high-resolution solutions in the 
offline phase. 
Sec.~\ref{sec:ROM} presents the details of our POD-LSTM ROM. Sec.~\ref{sec:num_res} discusses the numerical results obtained for an extension of the double-wind gyre forcing benchmark. 
Conclusions are outlined in Sec. \ref{sec:conclusion}.


\section{Governing equations} \label{sec:2qge}

We consider an ocean composed of two layers, 
each with uniform depth, density, and temperature. 
These layers, which are called isopycnals, occupy
a two-dimensional fixed spatial domain 
$\Omega=[x_0,x_f]\times[-L/2,L/2]$. 
We will call layer 1 the top layer.
We assume that the depths of the layers, $H_1$ and $H_2$, are much smaller than the meridional length $L$. The 2QGE describe the dynamics of these two layers, under some 
simplifying assumptions. We refer the reader to, e.g., \cite{Marshall1997,Chassignet1998,Berloff1999,DiBattista2001,BERLOFF_KAMENKOVICH_PEDLOSKY_2009}, 
for a thorough description of such assumptions.

Below, we state the 2QGE in non-dimensional form. 
For this purpose, we need to introduce some notation. 
Let $y$ be the non-dimensional vertical coordinate in the 2D domain
and $\delta = \frac{H_1}{H_1+H_2}$ the aspect ratio of the depths of the isopycnal layers.
We denote with $\sigma$ the friction coefficient at the bottom of the ocean and with $\nu$ the eddy viscosity coefficient. 
The 2QGE model linearizes the Coriolis frequency $f$ as follows: 
$f = f_0 + \beta y$, where $f_0$ is the local Earth rotation rate at $y = 0$, 
which is the center of the basin, and $\beta$ is the gradient of the Coriolis frequency. 
In the non-dimensionalization process, the Rossby number $Ro$, the Reynolds number $Re$, and the Froude number $Fr$ appear in the model. They are 
defined as 
\begin{equation*}
    Re=\frac{UL}{\nu},\quad Ro = \frac{U}{\beta L^2}, \quad Fr = \frac{f_0^2U}{g'\beta H},
\end{equation*}
where $U$ is the characteristic velocity scale, 
$H=H_1 + H_2$ is the total ocean depth,
$g' = g\Delta\rho/\rho_1$  is the reduced gravity, $g$ being the gravitational constant, $\Delta\rho$ the density difference between the two layers, and $\rho_1$ is the density of the top layer. 
Finally, let $(0,T]$ be a time interval of interest. 
The non-dimensional 2QGE read: find 
potential vorticities $q_l$ and stream functions $\psi_l$, for $l = 1, 2$, such that in $\Omega\times(0,T]$:
\begin{align}
    &\frac{\partial q_1}{\partial t} + \nabla \cdot \left(\left(\nabla \times \bm{\Psi}_1 \right)q_1\right) + \frac{Fr}{Re~\delta}\Delta(\psi_2-\psi_1) - \frac{1}{Re}\Delta q_1 = F, \label{eq:qge2_1}\\
    &\frac{\partial q_2}{\partial t} + \nabla \cdot \left(\left(\nabla \times \bm{\Psi}_2 \right)q_2\right) + \frac{Fr}{Re~(1-\delta)}\Delta(\psi_1-\psi_2) + \sigma\Delta\psi_2- \frac{1}{Re}\Delta q_2 = 0, \label{eq:qge2_2}\\
    &q_1 = Ro\Delta\psi_1 + y + \frac{Fr}{\delta}\left(\psi_2 - \psi_1 \right), \label{eq:kin_eq1} \\
    &q_2 = Ro\Delta\psi_2 + y + \frac{Fr}{1-\delta}\left(\psi_1 - \psi_2 \right), \label{eq:kin_eq2}
\end{align}
where  $\bm{\Psi}_l=(0, 0, \psi_l)$ and $F$ is forcing coming from the wind.
More details on this model can be found in, 
e.g., \cite{Besabe2024, San2012, Salmon1978, Medjo2000, Fandry1984, Mu1994}.

Following \cite{San2012}, we prescribe free-slip and impenetrable boundary conditions and start
the system from a rest state, i.e., we impose the following boundary and initial conditions:
\begin{align} \label{qge-bdry}
    \psi_l = 0, &\quad \mbox{on }\partial\Omega\times(0,T),\\
    q_l = y, &\quad \mbox{on }\partial\Omega\times(0,T),\\
    q_l(x,y) = q_l^0 = y, &\quad \mbox{in }\Omega\times\{0\},
\end{align}
for $l=1,2$.

\section{Full Order Method}\label{sec:fom}

Let us start with the time discretization of problem \eqref{eq:qge2_1}-\eqref{eq:kin_eq2}. 
Let $\Delta t = {T}/{N}$, where $N\in\mathbb{N}$, $t_n = t_0 + n\Delta t$, for $n = 1, \dots, N$, and denote by $f^n$ the approximation of the quantity $f$ at time $t_n$. Using Backward Difference Formula of order 1, the time-discrete problem reads: given the initial condition $q_l^0 = q_l(t_0)$, solve for $(q_l^{n+1}, \psi_l^{n+1})$ for $l = 1, 2$ and $n\geq0$:
\begin{align}
    &\frac{1}{\Delta t}q_1^{n+1} + \div\left(\left(\nabla \times \bm{\Psi}_1^{n}\right)q_1^{n+1}\right) + \frac{Fr}{Re~\delta} \Delta \left( \psi_2^{n+1} - \psi_1^{n+1} \right) - \frac{1}{Re}\Delta q_1^{n+1} = b_1^{n+1}, \label{eq:2qge-time-disc1}
     \\
    & Ro\Delta \psi_1^{n+1} + y + \frac{Fr}{\delta}\left(\psi_2^{n+1} - \psi_1^{n+1}\right) = q_1^{n+1}, \label{eq:2-qgetime-disc2}\\
    & \frac{1}{\Delta t}q_2^{n+1} + \div\left(\left(\nabla \times \bm{\Psi}_2^{n}\right)q_2^{n+1}\right) + \frac{Fr}{Re~(1-\delta)}\Delta \left( \psi_1^{n+1} - \psi_2^{n+1} \right) - \frac{1}{Re}\Delta q_2^{n+1}, \cl 
    &\quad + \sigma \Delta \psi_2^{n+1}= b_2^{n+1} \label{eq:2qge-time-disc3} \\
    & Ro\Delta \psi_2^{n+1} + y + \frac{Fr}{1-\delta}\left(\psi_1^{n+1} - \psi_2^{n+1}\right) = q_2^{n+1} ,\label{eq:2qge-time-disc4}
\end{align}
where $b_1^{n+1} = F + q_1^n/\Delta t$ and $b_2^{n+1} = q_2^n/\Delta t$. Note that we have linearized the convective terms in \eqref{eq:2qge-time-disc1} and \eqref{eq:2qge-time-disc3} with a first-order extrapolation.
The choice of the BDF1 scheme for time discretization is out of convenience. Of course, other discretization schemes are possible. See, e.g., \cite{Girfoglio_JCAM2023}, where we used the BDF2 scheme for the single layer quasi-geostrophic equations.

To contain the computational cost, we adopt a 
segregated algorithm to solve problem 
\eqref{eq:2qge-time-disc1}-\eqref{eq:2qge-time-disc4}. The algorithm is as follows:
given $(q_l^n,\psi_l^n)$, at time $t_{n+1}$ perform the following steps:
\begin{itemize}
    \item[-] Step 1: find the potential vorticity of the top layer $q_1^{n+1}$ such that
    \begin{align}\label{eq:seg-alg1}
        \frac{1}{\Delta t}q_1^{n+1} + \div\left(\left(\nabla \times \bm{\Psi}_1^{n}\right)q_1^{n+1}\right) - \frac{1}{Re}\Delta q_1^{n+1} = b_1^{n+1} 
        - \frac{Fr}{Re~\delta} \Delta \left( \psi_2^{n} - \psi_1^{n} \right).
    \end{align}
    \item[-] Step 2: find the stream function of the top layer $\psi_1^{n+1}$ such that:
    \begin{equation} \label{eq:seg-alg2}
         Ro\Delta \psi_1^{n+1} + y - \frac{Fr}{\delta} \psi_1^{n+1} = q_1^{n+1} - \frac{Fr}{\delta}\psi_2^{n}.
    \end{equation}
    \item[-] Step 3: find the potential vorticity  of the bottom layer $q_2^{n+1}$ such that:
    \begin{align} \label{eq:seg-alg3}
        &\frac{1}{\Delta t}q_2^{n+1} + \div\left(\left(\nabla \times \bm{\Psi}_2^{n}\right)q_2^{n+1}\right) - \frac{1}{Re}\Delta q_2^{n+1}= b_2^{n+1} -\sigma \Delta \psi_2^{n} \cl
        & \quad - \frac{Fr}{Re~(1 - \delta)} \Delta \left( \psi_1^{n+1} - \psi_2^{n} \right).
    \end{align}   
    \item[-] Step 4: find the stream function of the bottom layer $\psi_2^{n+1}$ such that:
    \begin{equation}\label{eq:seg-alg4}
         Ro\Delta \psi_2^{n+1} + y - \frac{Fr}{1-\delta} \psi_2^{n+1} = q_2^{n+1} - \frac{Fr}{1-\delta}\psi_1^{n+1} .
    \end{equation}
\end{itemize}
The advantage of this algorithm is that it decouples the computation of each variable. 

For the space discretization of eq.~\eqref{eq:seg-alg1}--\eqref{eq:seg-alg4}, we adopt a Finite Volume (FV) discretization. To begin, we partition the computational domain $\Omega$ into $N_C$ control volumes with $\Omega_i, i=1,\dots,N_C$,  such that $\bigcup_i \Omega_i=\Omega$ and $\Omega_i\cap\Omega_j = \emptyset$ for all $i\neq j$. Let
$\textbf{A}_j$ denote the surface vector of each face of the control volume, i.e., the surface
area of the face multiplied by the outward normal.
The FV approximation is obtained by integrating the 2QGE
in space and applying the Gauss divergence theorem:
\begin{align}
    &\frac{1}{\Delta t}q_{1,k}^{n+1} + \sum_j\phi_{1,j}^{n}q_{1,k}^{n+1,j} - \frac{1}{Re}\sum_j \left(\nabla q_{1,k}^{n+1}\right)_j\cdot \textbf{A}_j = b_{1,k}^{n+1} - \frac{Fr}{Re~\delta}\sum_j \left(\nabla(\psi_{2,k}^n - \psi_{1,k}^n)\right)_j\cdot \textbf{A}_j 
\label{eq:sp-disc-q1}\\
    & Ro\sum_j\left(\nabla \psi_{1,k}^{n+1}\right)_j\cdot\textbf{A}_j + y_k + \frac{Fr}{\delta}\left(\psi_{2,k}^n-\psi_{1,k}^{n+1}\right) = q_{1,k}^{n+1}, \label{eq:sp-disc-psi1} \\
    & \frac{1}{\Delta t}q_{2,k}^{n+1} + \sum_j\phi_{2,j}^{n}q_{2,k}^{n+1,j} - \frac{1}{Re}\sum_j \left(\nabla q_{2,k}^{n+1}\right)_j\cdot \textbf{A}_j = b_{2,k}^{n+1} \cl 
    & \quad + \left(\frac{Fr}{Re~(1-\delta)} - \sigma\right) \sum_j (\nabla\psi_{2,k}^n)_j\cdot \textbf{A}_j - \frac{Fr}{Re~(1-\delta)}\sum_j \left(\nabla\psi_{1,k}^{n+1}\right)_j\cdot \textbf{A}_j, \label{eq:sp-disc-q2} \\
     & Ro\sum_j\left(\nabla \psi_{2,k}^{n+1}\right)_j\cdot\textbf{A}_j + y_k - \frac{Fr}{1-\delta}\psi_{2,k}^{n+1}= q_{2,k}^{n+1} - \frac{Fr}{1-\delta}\psi_{1,k}^{n+1}, \label{eq:disc-psi2}
\end{align}
where $q_{l,k}$ and $\psi_{l,k}$ are the average potential vorticity and average stream function 
of layer $l$ over the control volume $\Omega_k$, 
$q_{l,k}^{j}$ is the potential vorticity associated to the centroid of the $j$-th face and normalized by the volume of $\Omega_k$,
$b_{l,k}$ is the average discrete forcing, and $y_k$ is the vertical coordinate of the centroid of the control volume $\Omega_k$.
For the discretization of the convective terms in \eqref{eq:sp-disc-q1} and \eqref{eq:sp-disc-q2}, we have introduced:
\begin{equation}\label{eq:flux}
    \phi_{l,j}=\left(\nabla \times \bm{\Psi}_{l,j}\right)\cdot\textbf{A}_j \text{ with } \bm{\Psi}_{l,j} = (0,0,\psi_{l,j})^T,
\end{equation}
which is calculated through a linear interpolation over neighboring cells using a central difference scheme which is of second-order. Also the diffusive terms are approximated using a second order central differencing scheme.  For further information on the numerical scheme, we refer the reader to \cite{Girfoglio2019, GirfoglioPSIZETA}.

It was shown in \cite{Besabe2024} that the numerical scheme presented in this section yields stable and accurate numerical solutions if the mesh is fine enough, i.e., mesh size $h<\delta_M$ where
$\delta_M$ is the Munk scale defined as:
\begin{align}
\delta_M = L \, \sqrt[3]{\dfrac{Ro}{Re}}. \label{eq:munk}
\end{align}
If one wishes to attain a stable solution with a coarse mesh, i.e., with size $h>\delta_M$, additional artificial dissipation needs to be introduced in the the system as shown in, e.g., \cite{San2012, Besabe2024}.


Lastly, we note that all numerical simulations are executed using GEA \cite{GEA, GirfoglioFVCA10}, an open-source software package which is built upon C++ finite volume library OpenFOAM\textsuperscript{\textregistered} \cite{Weller1998}.

\section{Reduced Order Model}\label{sec:ROM}

This section discusses the main building blocks
of our data-driven reduced order model for the 2QGE: proper orthogonal decomposition (POD) and the long short-term memory (LSTM) architecture.

Let us consider first the case where time is the only parameter. 
We assume that the solution $(q_l,\psi_l), l = 1,2$, of the system \eqref{eq:qge2_1}--\eqref{eq:kin_eq2} can be written as follows
\begin{align}
    q_l(t,\xbf)\approx q_{l}^r(t,\xbf) &= \tildeq_l(\xbf) + \sum_{i=1}^{N_{q_l}^r} \alpha_{l,i}(t)\varphi_{l,i}(\xbf),\label{eq:q_approx} \\
    \psi_l(t,\xbf)\approx \psi_{l}^r(t, \xbf) &= \tildep_l(\xbf) + \sum_{i=1}^{N_{\psi_l}^r} \beta_{l,i}(t)\xi_{l,i}(\xbf), \label{eq:psi_approx}
\end{align}
where $q_{l}^r$ and $\psi_{l}^r$ are the reduced order approximations and
$\tildeq_l$ and $\tildep_l$ are the time-averaged potential vorticity and stream function for each layer computed by
\begin{equation}\label{eq:time_av}
    \tildeq_l(\xbf) = \frac{1}{N^t}\sum_{p=1}^{N^t} q_l(t_p,\xbf), \quad \tildep_l (\xbf) = \frac{1}{N^t}\sum_{p=1}^{N^t} \psi_l(t_p,\xbf).
\end{equation}

Next, we consider the case where there are physical parameters of interest 
(e.g., the aspect ratio $\delta$, friction coefficient $\sigma$, etc.)
that vary in a given range.
Let us collect all such parameters in a vector $\bm{\mu}$ belonging
to parameter space $\caldf$. 
We assume that, for a given $\bm{\mu}$, the solution $(q_l,\psi_l), l = 1,2$, of the system \eqref{eq:qge2_1}--\eqref{eq:kin_eq2} can be written as follows
\begin{align}
    q_l(t,\xbf,  \bm{\mu})\approx q_{l}^r(t,\xbf, \bm{\mu}) &= \tildeq_l^0(\xbf) + \sum_{i=1}^{N_{q_l}^r} \alpha_{l,i}(t,\bm{\mu})\varphi_{l,i}(\xbf),\label{eq:q_approx_mu} \\
    \psi_l(t,\xbf, \bm{\mu})\approx \psi_{l}^r(t, \xbf, \bm{\mu}) &= \tildep_l^0(\xbf) + \sum_{i=1}^{N_{\psi_l}^r} \beta_{l,i}(t,\bm{\mu})\xi_{l,i}(\xbf), \label{eq:psi_approx_mu}
\end{align}
where, again, $q_{l}^r$ and $\psi_{l}^r$ are the reduced order approximations, while
$\tildeq_l^0$ and $\tildep_l^0$ are zero-th
order approximations of $\tildeq_l(\xbf,\bm{\mu})$ and 
$\tildep_l(\xbf, \bm{\mu})$. See Sec.~\ref{sec:pod} for more details on this. 

Sec.~\ref{sec:pod} describes how to find the reduced basis functions $\varphi_{l,i}$ in \eqref{eq:q_approx},\eqref{eq:q_approx_mu} and $\xi_{l,i}$ in
\eqref{eq:psi_approx}, \eqref{eq:psi_approx_mu}, and the cardinalities 
$N_{q_l}^r$ and $N_{\psi_l}^r$
of the reduced bases. 
Sec.~\ref{sec:pod-lstm} explains
how we find coefficients $\alpha_{l,i}$ in \eqref{eq:q_approx},\eqref{eq:q_approx_mu} and $\beta_{l,i}$ in
\eqref{eq:psi_approx},\eqref{eq:psi_approx_mu} for any new value of $t$ and $\bm{\mu}$.

\subsection{Construction of the POD bases for the reduced spaces} \label{sec:pod}

Let us start again from the case where time is the only parameter. The solutions associated to the sampling of the time interval $\{ t_p \}_{p = 1}^{N^t}$, also called snapshots, will be used to generate a basis for the reduced spaces via POD.
The reader interested in learning about
the mathematics, applications, and variations of POD is referred to, e.g., \cite{Christensen1999,Chatterjee2000,Kunisch2003,Sieber2016,Rozza2022}. Below, we describe how we use
POD to find the reduced basis functions.

From the snapshots, we compute
the fluctuations from the time-averaged fields defined by
\begin{align} 
    q_l' (t_p, \xbf) &= q_l(t_p,\xbf) - \tildeq_l(\xbf), \label{eq:fluct_q} \\
    \psi_l' (t_p, \xbf) &= \psi_l(t_p, \xbf) - \tildep_l(\xbf), \label{eq:fluct_psi}
\end{align}
for $l = 1,2$ and $p = 1, \dots, N^t$. 
The time-averaged fields are computed as in 
\eqref{eq:time_av}. Let $N_h$ be the number of degrees of freedom for the FOM solution.
Fluctuations \eqref{eq:fluct_q} and \eqref{eq:fluct_psi} are collected in 
snapshot matrices. For $\Phi\in\{q_l',\psi_l'\}$, $l = 1,2$, the snapshot matrix $\calsf_\Phi\in\mathbb{R}^{N_h\times N^t}$ is given by
\begin{equation}\label{eq:snapshots}
    \calsf_\Phi = \left[\Phi(t_1),\dots, \Phi(t_{N^t})\right]. 
\end{equation}

Singular value decomposition (SVD) is applied
to matrix $\calsf_\Phi$ for $\Phi\in\{q_l',\psi_l'\}$, $l = 1,2$:
\begin{equation} \label{eq:svd}
    \calsf_\Phi = \caluf_\Phi\calsigf_\Phi\calvf_\Phi^T
\end{equation}
where $\caluf_\Phi \in \mathbb{R}^{N_h\times N^t}$ and $\calvf_\Phi \in \mathbb{R}^{N^t \times N_h}$ are matrices whose columns are the left and right singular vectors of $\calsf_\Phi$, respectively, and $\calsigf_\Phi \in \mathbb{R}^{N^t \times N^t}$ is the diagonal matrix whose diagonal entries
are the singular values $\sigma_\Phi^i$ of $\calsf_\Phi$ arranged in descending order. 
The POD basis functions are the first $N_\Phi^r$ left singular vectors in $\caluf_\Phi$, where $N_\Phi^r$ is chosen such that for a
user-defined threshold $0<\delta_\Phi < 1$:
\begin{equation} \label{eq:energy}
    \frac{\sum_{i=1}^{N_\Phi^r} \sigma_\Phi^i}{\sum_{i=1}^{N^t} \sigma_\Phi^i} \geq \delta_\Phi.
\end{equation}
Eq.~\eqref{eq:energy} means that by selecting the first $N_\Phi^r$ POD modes as basis functions
we retain a user-defined fraction ($\delta_\Phi$) 
of the singular value energy of the system.
We will denote with $\caluf_{\Phi, N_{\Phi}^r}$ the matrix consisting of the first $N_{\Phi}^r$ columns of $\caluf_\Phi$.
The reduced basis space is the span of the 
columns of $\caluf_{\Phi, N_{\Phi}^r}$, which represent the most dominant spatial structures from each field of interest. Obviously, one succeeds in generating a
\emph{reduced} basis if $N_\Phi^r \ll N^t$, i.e., if with a small number of POD modes one can capture
reasonably well the system dynamics. 

Now that we have the basis functions, we can use eq.~\eqref{eq:q_approx}--\eqref{eq:psi_approx} to approximate the snapshots:
\begin{align}
    q_l(t_p,\xbf)\approx q_{l}^r(t_p,\xbf) &= \tildeq_l(\xbf) + \sum_{i=1}^{N_{q_l}^r} \alpha_{l,i}(t_p)\varphi_{l,i}(\xbf),\label{eq:q_approx_snap} \\
    \psi_l(t_p,\xbf)\approx \psi_{l}^r(t_p, \xbf) &= \tildep_l(\xbf) + \sum_{i=1}^{N_{\psi_l}^r} \beta_{l,i}(t_p)\xi_{l,i}(\xbf), \label{eq:psi_approx_snap}
\end{align}
for $l = 1,2$ and $p = 1, \dots, N^t$, where
the coefficients $\alpha_{l,i}(t_p)$ and $\beta_{l,i}(t_p)$ for each time instant $t_p$ are the entries of $p$-th column of matrix $\calcf_\Phi = \caluf_{\Phi, N_{\Phi}^r}^T\calsf_\Phi\in\mathbb{R}^{N_\Phi^r\times N_t}$.

The next subsection explains how we find coefficients
$\alpha_{l,i}$ and $\beta_{l,i}$ for times that do not belong to partition $\{ t_p \}_{p = 1}^{N^t}$.

\begin{rem}\label{rem1}
An alternative to the generation of separate POD bases for
vorticities and stream functions is to apply 
POD only to the snapshot matrices of the vorticities  to get basis
functions $\varphi_{l,i}$. Then, 
the reduced basis functions for the
stream functions are obtained from solving Poisson problems: 
\begin{equation} \label{eq:poisson}
    \Delta \xi_{l,i} = -\varphi_{l,i},
\end{equation}
for $i = 1, ..., N_{q_l}^r$ and $l = 1,2$.
Then, using linearity argument, one only needs to compute the modal coefficients once per time step for the vorticities and use them for the reconstruction and prediction of the stream functions too. This means that
eq.~\eqref{eq:q_approx}-\eqref{eq:psi_approx} become
\begin{align}
    q_l(t,\xbf)\approx q_{l}^r(t,\xbf) &= \tildeq_l(\xbf) + \sum_{i=1}^{N_{q_l}^r} \alpha_{l,i}(t)\varphi_{l,i}(\xbf),\label{eq:q_approx2} \\
    \psi_l(t,\xbf)\approx \psi_{l}^r(t, \xbf) &= \tildep_l(\xbf) + \sum_{i=1}^{N_{q_l}^r} \alpha_{l,i}(t)\xi_{l,i}(\xbf). \label{eq:psi_approx2}
\end{align}
This strategy has been used for the single layer quasi-geostrophic equations in \cite{Rahman-2019, Mou-2020, Rahman-2018, Golzar2023, QGE_review}.
\end{rem}

If, in addition to time, there are variable physical parameters ${\bm \mu}$, we need to 
sample $\caldf\times (0,T]$. For this, we collect the solutions associated to partition
$\{\bm{\mu}_k\}_{k=1}^{N^d} \times \{ t_p \}_{p = 1}^{N^t}$ of $\caldf\times (0,T]$.
From the snapshots, we compute
the fluctuations from the time-averaged fields defined by
\begin{align} 
    q_l' (t_p, \xbf, \bm{\mu}_k) &= q_l(t_p,\xbf, \bm{\mu}_k) - \tildeq(\xbf, \bm{\mu}_k), \label{eq:fluct_q_mu} \\
    \psi_l' (t_p, \xbf, \bm{\mu}_k) &= \psi_l(t_p, \xbf, \bm{\mu}_k) - \tildep(\xbf, \bm{\mu}_k), \label{eq:fluct_psi_mu}
\end{align}
for $l = 1,2$, $k = 1, \dots, N^d$, and $p = 1, \dots, N^t$. Let $N_s = N^d\times N^t$.
Fluctuations \eqref{eq:fluct_q_mu} and 
\eqref{eq:fluct_psi_mu}
are stored in snapshot matrices
$\calsf_\Phi\in\mathbb{R}^{N_h\times N_s}$:
\begin{equation}\label{eq:snapshots_mu}
    \calsf_\Phi = \left[\Phi(\bm{\mu}_1,t_1),\dots, \Phi(\bm{\mu}_1,t_{N^t}),\dots,\Phi(\bm{\mu}_{N^d},t_{N^t})\right],
\end{equation}
for $\Phi\in\{q_l',\psi_l'\}$, $l = 1,2$.
At this point, the procedure is the same as in the case of time as the only parameter, i.e., one applies SVD and retains the first 
$N_\Phi^r$ left singular vectors according to criterion \eqref{eq:energy} to generate the POD basis. With that, one uses \eqref{eq:q_approx_mu}-\eqref{eq:psi_approx_mu} to approximate the snapshots:
\begin{align}
    q_l(t_p,\xbf, \bm{\mu}_k)\approx q_{l}^r(t_p,\xbf, \bm{\mu}_k) &= \tildeq_l(\xbf, \bm{\mu}_k) + \sum_{i=1}^{N_{q_l}^r} \alpha_{l,i}(t_p, \bm{\mu}_k)\varphi_{l,i}(\xbf),\label{eq:q_approx_snap_mu} \\
    \psi_l(t_p,\xbf, \bm{\mu}_k)\approx \psi_{l}^r(t_p, \xbf, \bm{\mu}_k) &= \tildep_l(\xbf, \bm{\mu}_k) + \sum_{i=1}^{N_{\psi_l}^r} \beta_{l,i}(t_p, \bm{\mu}_k)\xi_{l,i}(\xbf), \label{eq:psi_approx_snap_mu}
\end{align}
for $l = 1,2$, $k = 1, \dots, N^d$, and $p = 1, \dots, N^t$, where
the coefficients $\alpha_{l,i}(t_p, \bm{\mu}_k)$ and $\beta_{l,i}(t_p, \bm{\mu}_k)$ for each time instant $t_p$ 
and parameter $\bm{\mu}_k$
are 
computed in the same way as when time is the only parameter.

For a time $t$ and parameter $\bm{\mu}$
not in the sampling set $\{\bm{\mu}_k\}_{k=1}^{N^d} \times \{ t_p \}_{p = 1}^{N^t}$, 
we need to find $\tildeq_l^0$ and $\tildep_l^0$
in \eqref{eq:q_approx_mu}-\eqref{eq:psi_approx_mu}, together with 
coefficients
$\alpha_{l,i}$ and $\beta_{l,i}$.
The zero-th order approximations $\tildeq_l^0$ and $\tildep_l^0$
are simply given by 
\begin{equation}\label{eq:time_av_mu}
    \tildeq_l^0(\xbf) = \frac{1}{N^t}\sum_{p=1}^{N^t} q_l(t_p,\xbf, \bm{\mu}_c), \quad \tildep_l^0(\xbf) = \frac{1}{N^t}\sum_{p=1}^{N^t} \psi_l(t_p,\xbf, \bm{\mu}_c),
\end{equation}
where $\bm{\mu}_c$ is the sample in 
$\{\bm{\mu}_k\}_{k=1}^{N^d}$ closest to 
$\bm{\mu}$ in the Euclidean norm. How to 
find coefficients $\alpha_{l,i}$ and $\beta_{l,i}$ is explained in the next subsection.

\subsection{The data-driven ROM based on LSTM} \label{sec:pod-lstm}


We develop a data-driven ROM using a deep learning technique \cite{Wang2018} that performs well in the prediction of a time-series behavior: a type of recurrent neural network (RNN) called long short-term memory. 
The LSTM architecture was introduced in \cite{LSTM1997} to resolve
the problem of vanishing or exploding gradients that RNNs often encounter \cite{Bengio1994}. 
The main mechanisms of the LSTM architecture are gating activation functions and cell states that serve as the memory of the network. 
It is thanks to the gating functions that 
vanishing gradients are avoided, while the fact that the gradient of the cell states is bounded by $1$ avoids exploding gradients \cite{Bayer2015}. 
Without incurring into gradient problems, 
the LSTM learns long term behavior better than regular RNNs.

\begin{figure}[h!]
    \centering
    \includegraphics[width=0.8\textwidth]{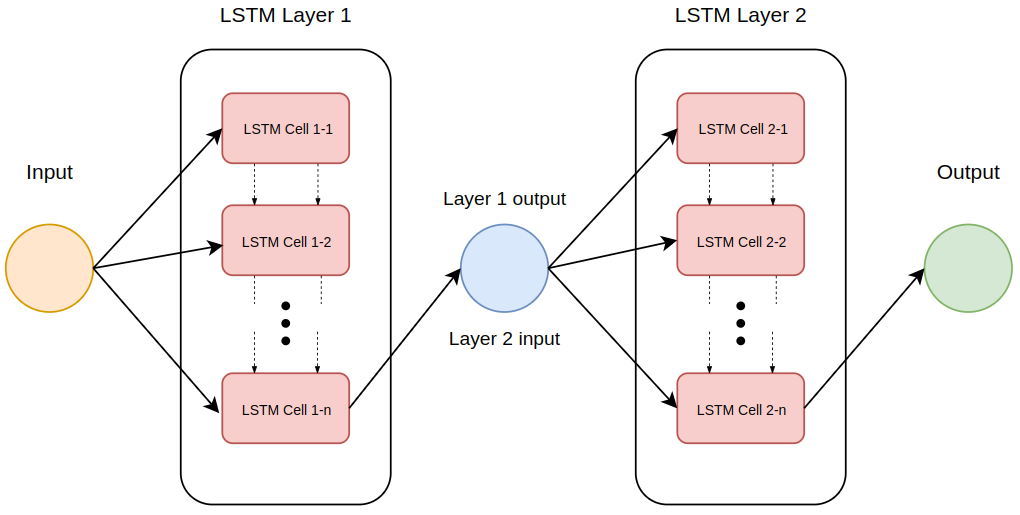}
    \caption{Illustration of an LSTM architecture with 2 layers and $n$ cells (or memory blocks) in each layer.}
    \label{fig:lstm-network}
\end{figure}

Like any RNN, LSTM contains recurrent layers and neurons. Within these layers, there are memory blocks, also called LSTM cells. An illustration of a LSTM network with two layers and $n$ memory blocks in each layer is shown as an example in Fig.~\ref{fig:lstm-network}. Every layer is its own network, which means that if we have two layers the input of the second layer is the output of the first layer, with the option to pass them through an activation function before passing it to the next layer.

Before explaining how each cell in Fig.~\ref{fig:lstm-network} works, we specify the input and output of our LSTM network. In order to do that, we need to introduce 
one important hyperparameter: the lookback window $\sigma_L$, which indicates how many steps in the history of the data series the LSTM considers, both in training and prediction.
Let us consider the reduced space for the potential vorticity in 
layer $l$ as representative.
For simplicity of notation, we denote with $N$
the number of retained modes $N_{q_l}^r$. 
For a given time $t_p$ and parameter vector $\bm{\mu}_k$, 
we store the time information $t_p,\dots, t_{p-\sigma_L+1}$, parameter vector $\bm{\mu}_k$, and coefficients 
$\alpha_{l,i}(t_p)$ in \eqref{eq:q_approx_snap}
in a $\sigma_L\times (N + N^d + 1)$ matrix
\begin{equation}\label{eq:input}
    \left[\begin{array}{ccccc}
        \bm{\mu}_k & t_p &\alpha_{l,1}(t_p) & \cdots & \alpha_{l,N}(t_p) \\
        \vdots & \vdots & \vdots & \ddots & \vdots \\
        \bm{\mu}_k & t_{p-\sigma_L+1} & \alpha_{l,1}(t_{p-\sigma_L + 1}) & \cdots & \alpha_{l,N}(t_{p-\sigma_L + 1})
    \end{array}\right].
\end{equation}
Through the LSTM architecture, we aim to match the
training matrix \eqref{eq:input}, i.e., the
input at time $t_p$ and for parameter vector $\bm{\mu}_k$, to the corresponding output training vector $\left[\alpha_{l,1} (t_{p+1}),\dots,\alpha_{l,K}(t_{p+1})\right]$, 
i.e., the coefficients for the reduced
order approximation at time $t_{p+1}$. Below, we
explain how each cell in Fig.~\ref{fig:lstm-network}
works.

\begin{figure}[h!]
    \centering
    \includegraphics[width=0.8\textwidth]{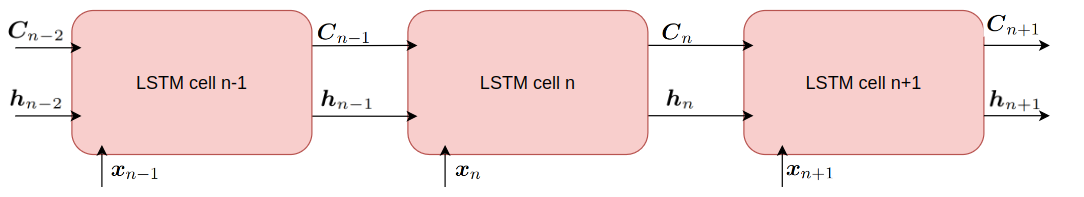}
    \caption{Connections among adjacent LSTM cells.}
    \label{fig:lstm-cells}
\end{figure}

\begin{figure}[h!]
    \centering
    \includegraphics[width=0.6\textwidth]{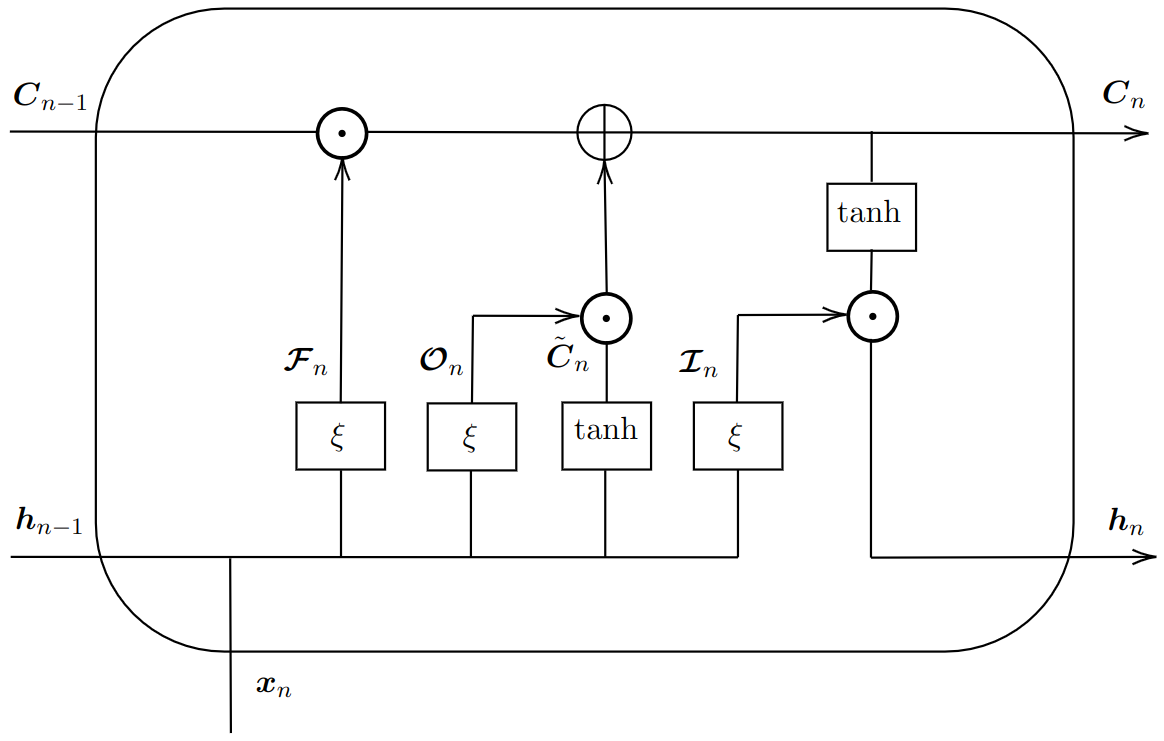}
    \caption{Schematic diagram of an LSTM cell.}
    \label{fig:lstm-diagram}
\end{figure}

The generic LSTM cell, denoted with $n$, manages the flow of training/testing data by selectively adding information through the input gate $\bm \calif_n$, forgetting some of such data in the forget gate $\bm \calff_n$, and then passing them through the next cell ($n+1$) through the output gate $\bm \calof_n$. 
For simplicity, let us denote the cell input with $\bm x_n$, which comes from either the input of the network or an intermediate input as shown in Fig.~\ref{fig:lstm-network}. 
The generic cell $n$ takes as input also the output of the preceding cell $\bm h_{n-1}$, which can be interpreted as the short term memory of the network. See Fig.~\ref{fig:lstm-network}. 
Let us denote with $\left[\bm h_{n-1}, \bm x_n\right]$ the concatenation of $\bm h_{n-1}$ and $\bm x_n$.
We have
\begin{align} 
   \bm \calif_n &= \xi\left(W_I \left[\bm h_{n-1}, \bm x_n\right] + \bm b_I\right), \cl
   \bm \calff_n &= \xi\left(W_F \left[\bm h_{n-1}, \bm x_n\right] + \bm b_F\right), \cl
   \bm \calof_n &= \xi \left( W_O \left[\bm h_{n-1}, \bm x_n\right] + \bm b_O \right), \el
\end{align}
where $\xi$ is the sigmoid function with range in $(0,1)$, $W_I$, $W_F$, and $W_O$
are matrices of weights, while $\bm b_I$, $\bm b_F$, and $\bm b_O$
are vectors of biases. The cell output $\bm h_n$ is computed as 
\begin{align} 
    \bm h_n &= \bm \calof_n \odot \tanh\left( \bm C_n\right), \el
\end{align}
where $\odot$ is the Hadamard product and $\bm C_n$ is the cell state given by 
\begin{align} 
    \bm C_n &= \bm \calff_n \odot \bm C_{n-1} + \bm \calif_n \odot \widetilde{\bm 
 C}_{n},\el
\end{align}
with $\widetilde{\bm C}_n$ representing the updated cell state
\begin{align} 
    \widetilde{\bm C}_n &= \tanh\left(
 W_C  \left[\bm h_{n-1}, \bm x_n\right] + \bm  b_C \right). \el
\end{align}
Here, $W_C$ is a matrix of weights and $\bm b_C$
is a vector of biases. The cell state $\bm C_n$ flows through the network passing some linear operations and can be seen as the long term memory of the network. See Fig.~\ref{fig:lstm-cells}.
A schematic diagram of an LSTM cell as described above is shown in 
Fig.~\ref{fig:lstm-diagram}. For a detailed discussion on the architecture and mathematics behind LSTM networks, see, e.g., \cite{Sherstinsky2020}.

Above, we have explained how the LSTM network
for the potential vorticity in 
layer $l$ works. We will call that network
$\calmf_{q_l}$, $l = 1, 2$.
The LSTM network for the stream function in 
layer $l$, denoted with $\calmf_{\psi_l}$, will be built analogously, using the same lookback window $\sigma_L$.

\section{Numerical Results} \label{sec:num_res}

To assess the performance of POD-LSTM ROM presented in Sec.~\ref{sec:ROM} in terms of both reconstruction and prediction, we consider an extension \cite{Besabe2024} of the so-called double-wind gyre forcing experiment, which is a classical benchmark for new numerical models of geophysical flows \cite{Nadiga2001, Holm2003, Greathbatch2000, Monteiro2014, Monteiro2015, Girfoglio_JCAM2023, QGE-review, San2015, Girfoglio2023}.

The computational domain for this benchmark is $\Omega = [0,1]\times[-1,1]$ and the forcing term in \eqref{eq:qge2_1} is $F = \sin(\pi y)$.
For the FOM, we use a mesh size $h=1/256$ which is about 3 times smaller than the Munk scale for the chosen value of eddy viscosity $\nu=155$ m\textsuperscript{2}s\textsuperscript{-1}. 
The time step is set to $\Delta t = 2.5\mathrm{E}-05$. 
We collect snapshots over time interval $[10,50]$ 
every $0.1$ unit of time, i.e., 401 snapshots for each variable ($q_1', q_2', \psi_1',$ and $\psi_2'$), and 
build the corresponding POD reduced spaces as explained in Sec.~\ref{sec:pod}. 

In Sec.~\ref{sec:asses-pod}, we fix all the physical parameters, keeping time as the only parameter. Specifically, we set
\begin{equation} \label{eq:params1}
    Ro = 0.001,\text{ }Re = 450,\text{ }Fr = 0.1,\text{ }\sigma = 0.005,\text{ and } \delta = 0.5.
\end{equation}
With the POD-LSTM ROM, we will reconstruct 
the time-averaged quantities $\widetilde{q}_l$ and $\widetilde{\psi}_l$, $l = 1,2$,
during the $[10,50]$ time interval to evaluate its performance in system identification. Additionally, 
we will use the POD-LSTM ROM to predict the same
quantities over $(50,100]$.
We will also consider two metrics to evaluate the system identification and the prediction: 
\begin{itemize}
    \item[-]  a root mean square error for the field $\Phi\in\{q_1,q_2,\psi_1,\psi_2\}$ \cite{Rahman-2019}:
\begin{equation} \label{eq:rmse}
    \varepsilon_{\Phi}^{(1)} = \sqrt{\frac{1}{N_x N_y}\sum_{i=1}^{N_x}\sum_{j=1}^{Ny}\abs{\widetilde{\Phi}_{i,j}^{\text{FOM}}-\widetilde{\Phi}_{i,j}^{\text{ROM}}}^2},
\end{equation}
where $N_x$ and $N_y$ are the number of cells in the $x$ and $y$ directions, respectively, $\widetilde{\Phi}^{\text{FOM}}$ and $\widetilde{\Phi}^{\text{ROM}}$ are the time-averaged
fields $\Phi$ computed by the FOM and POD-LSTM, respectively;
\item[-] the relative $L^2$ error:
\begin{equation} \label{eq:l2-error}
    \varepsilon_{\Phi}^{(2)} = \frac{\norm{\widetilde{\Phi}^{\text{FOM}}-\widetilde{\Phi}^{\text{ROM}}}_{L^2(\Omega)}}{\norm{\widetilde{\Phi}^{\text{FOM}}}_{L^2(\Omega)}}.
\end{equation}
\end{itemize}
In these metrics, the FOM solution is considered to be the \textit{true} solution. 
 
Finally, we will compare the evolution of enstrophy:
\begin{equation} \label{eq:enstrophy}
    \mathcal{E}_l = \frac{1}{2}\int_\Omega q_l^2\,d\Omega,
\end{equation}
and kinetic energy
\begin{equation} \label{eq:kin-energy}
    E_{l} = \frac{1}{2}\int_{\Omega} \left(\frac{\partial \psi_l}{\partial x}\right)^2 + \left(\frac{\partial \psi_l}{\partial y}\right)^2 \,d\Omega,
\end{equation}
given by FOM and ROM for $l = 1,2$.

In Sec.~\ref{sec:param-study}, we perform a parametric study, with the aspect ratio $\delta$
as the variable parameter, in addition to time. 
We will reconstruct the same quantities and use the same error metrics discussed above.

\subsection{Reconstruction and prediction with time as the only parameter}
\label{sec:asses-pod}

\begin{figure}[htb!]
\centering
    \begin{subfigure}[b]{0.47\textwidth}
         \centering
         \includegraphics[width=\textwidth]{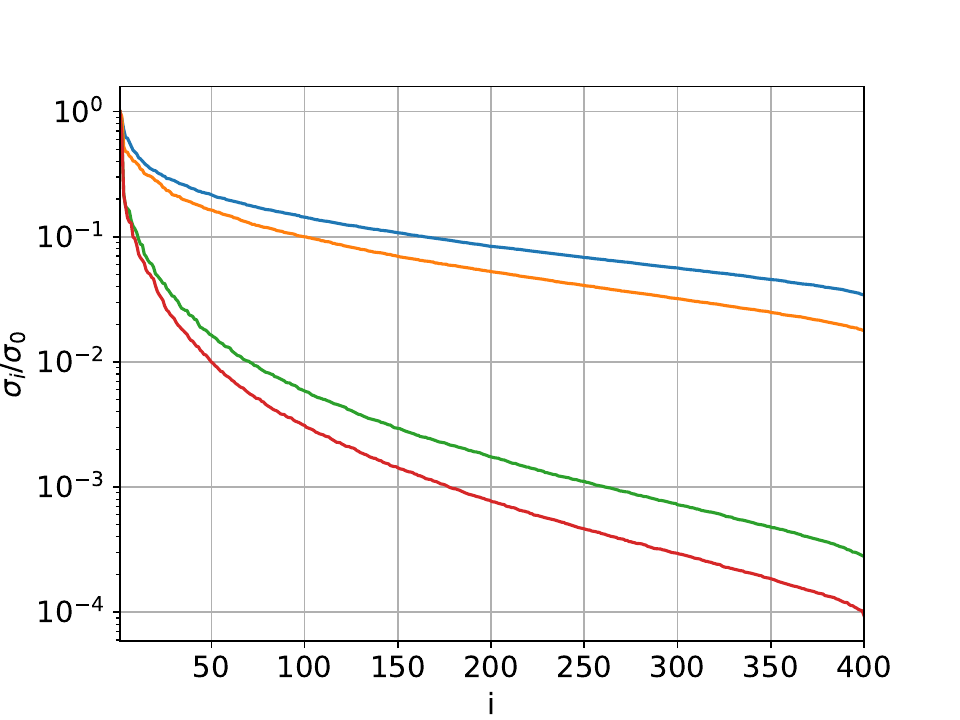}
         \caption{\scriptsize{Eigenvalue decay}}
     \end{subfigure}
    \begin{subfigure}[b]{0.47\textwidth}
         \centering
         \includegraphics[width=\textwidth]{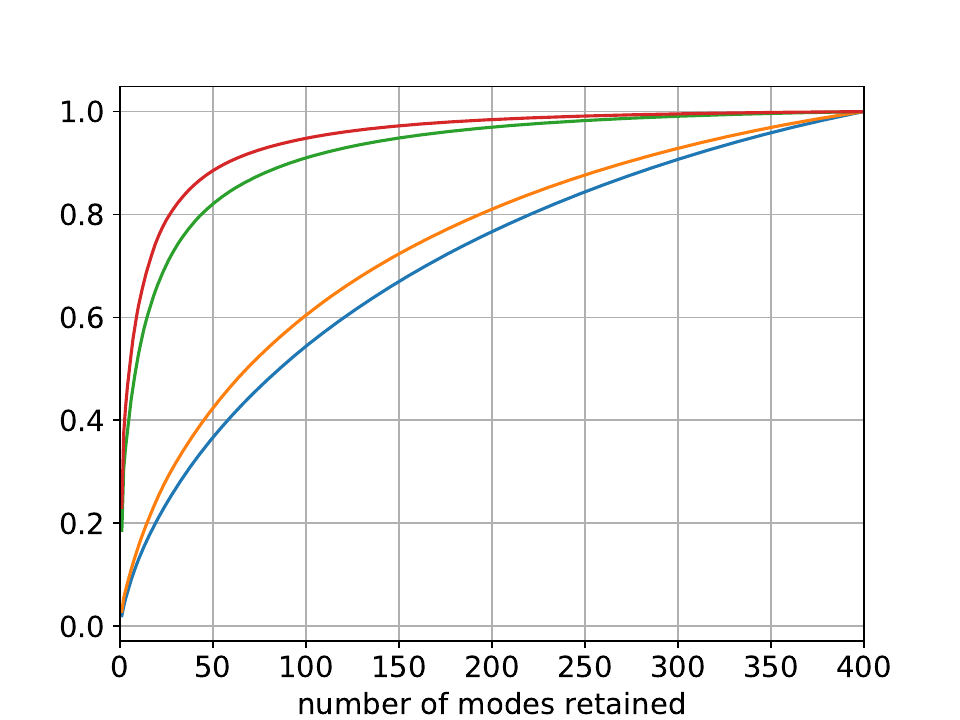}
         \caption{\scriptsize{Energy accumulation}}
     \end{subfigure}

     \begin{subfigure}[h]{0.4\textwidth}
         \centering
         \vspace{0.3cm}
         \includegraphics[width=\textwidth]{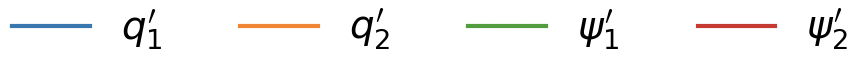}
     \end{subfigure}
\caption{(a) Singular value decay and (b) energy accumulation of the potential vorticity fluctuations and stream function fluctuations in both layers.}
\label{fig:energy-eigen}
\end{figure}

We applied the procedure described in Sec.~\ref{sec:pod} to generate to POD basis, starting from the 401 snapshots saved from the FOM simulation using the parameters specified in \eqref{eq:params1}.
The decay of the singular values of the potential vorticity fluctuations and stream function fluctuations in both layers is shown in Fig.~\ref{fig:energy-eigen} (a). We see that in both layers the decay for the potential vorticity fluctuation is much slower than the decay of the stream function fluctuation.
The energy accumulation for all the fluctuations
is shown in Fig.~\ref{fig:energy-eigen} (b).
Typically, one would want to retain $99.99\%$ of the singular value energy of the system.
For this test case, it would mean that we need 
$N_{q_1}^r = N_{q_2}^r = 399$, $N_{\psi_1}^r = 398$, $N_{\psi_2}^r = 395$.
Obviously, these values are not viable since there would likely be no reduction in the computational time
when compared to the FOM. 
If we were to retain only
$50\%$ of the energy, the number of basis functions would be reasonable for the stream function fluctuations ($N_{\psi_1}^r = 8$ and $N_{\psi_2}^r = 5$) but still too large 
for the potential vorticity fluctuations ($N_{q_1}^r = 86$ and $N_{q_2}^r = 68$).

Concerning the LSTM networks, for $q_1'$ and $q_2'$
we choose a single layer with $100$ cells, while
we use $3$ hidden layers with $50$ cells each
for $\psi_1'$ and $\psi_2'$. 
This difference in the networks is due to the 
fact that, as we shall see, the time-averaged
stream functions
$\widetilde{\psi}_1$ and $\widetilde{\psi}_2$ have more complex patterns than the time-averaged
vorticities $\widetilde{q}_1$ and $\widetilde{q}_2$. The simpler structure of the patterns in $\widetilde{q}_1$ and $\widetilde{q}_2$
is such that using more than one layer would lead
to overfitting. Concerning the optimization method, 
following \cite{Rahman-2019} we choose 
a type of stochastic gradient descent widely used in machine learning \cite{Kingma2014} called
Adam, which optimizes the weights by minimizing the
mean square error (MSE). We summarize the hyperparameters for our LSTM networks in Table \ref{tab:hyperparam}. 
We have intentionally left out of Table \ref{tab:hyperparam} the lookback window $\sigma_L$
because we will show how the networks are sensitive to it.

\begin{table}[h!]
    \centering
    \begin{tabular}{|l|c|c|}
        \hline 
        Hyperparameters & $\calmf_{q_l}$ & $\calmf_{\psi_l}$ \\
        \hline
        Number of layers & 1 & 3 \\
        \hline
        Number of cells per layer & 100 & 50 \\
        \hline
        Batch size & 8 & 16 \\
        \hline
        Epochs & $500$ & $500$ \\
        \hline
        Activation function & $\tanh$ & $\tanh$ \\
        \hline
        Validation & $20\%$ & $20\%$ \\
        \hline
        Testing : testing ratio & $1:4$ & $1:4$ \\
        \hline
        Loss function & MSE & MSE \\
        \hline
        Optimizer & Adam & Adam \\
        \hline
        Learning rate & 1E-02 & 1E-02 \\
        \hline
        Drop out probability & - & 0.1\\
        \hline
        Weight decay & 1E-05 & 1E-05 \\
        \hline
    \end{tabular}
    \caption{Hyperparameters for the LSTM network component of the POD-LSTM ROMs.}
    \label{tab:hyperparam}
\end{table}

\begin{figure}[htb!]
\centering
\begin{tabular}{cccccc}
       & FOM & \hspace{-0.4cm}$\sigma_L = 1$ & \hspace{-0.4cm} $\sigma_L = 2$ & \hspace{-0.4cm} $\sigma_L = 5$ & \hspace{-0.4cm} $\sigma_L=10$ \\
    $\tildeq_1$ & \includegraphics[align=c,scale = 0.3]{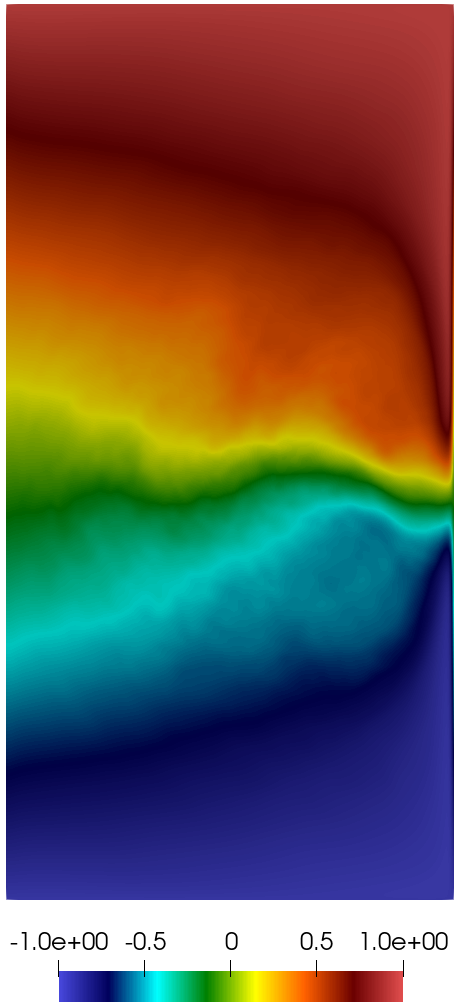} & \hspace{-0.4cm}\includegraphics[align=c,scale = 0.3]{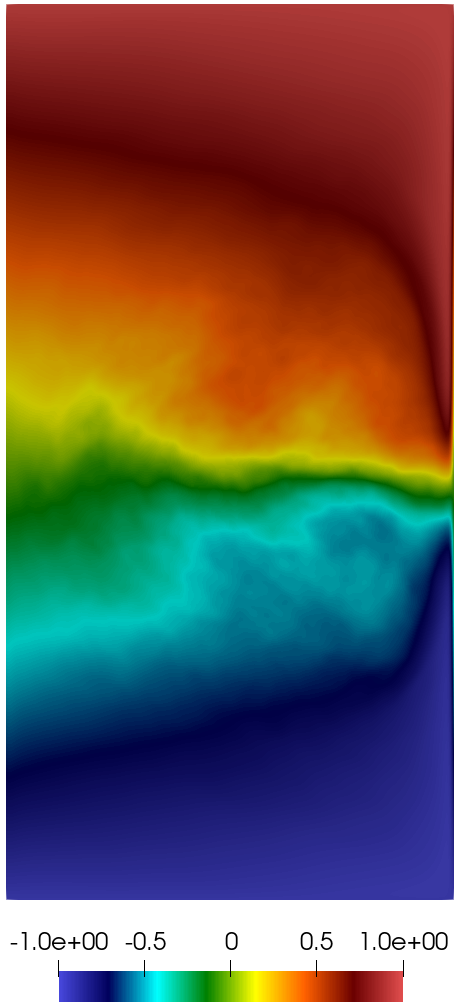} & \hspace{-0.4cm}\includegraphics[align=c,scale = 0.3]{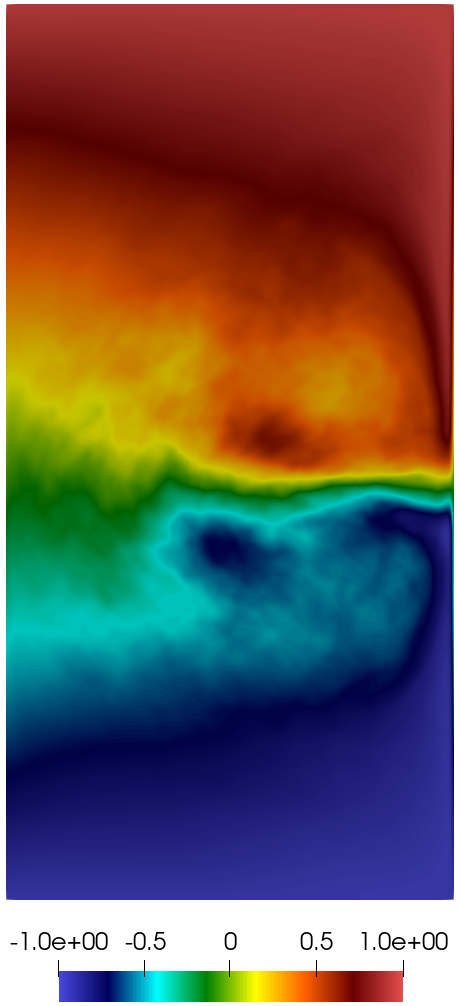} & \hspace{-0.4cm}\includegraphics[align=c,scale = 0.3]{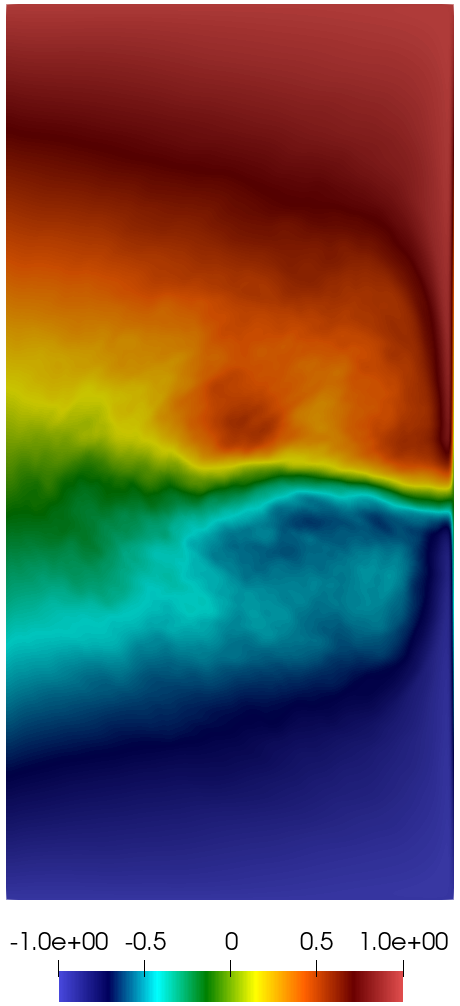} & \hspace{-0.4cm}\includegraphics[align=c,scale = 0.3]{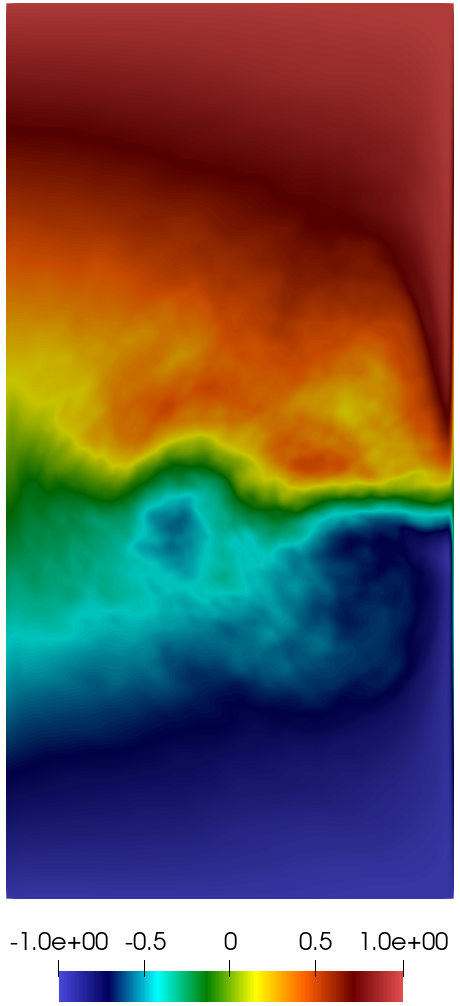} \\ 
    $\tildeq_2$ & \includegraphics[align=c,scale = 0.3]{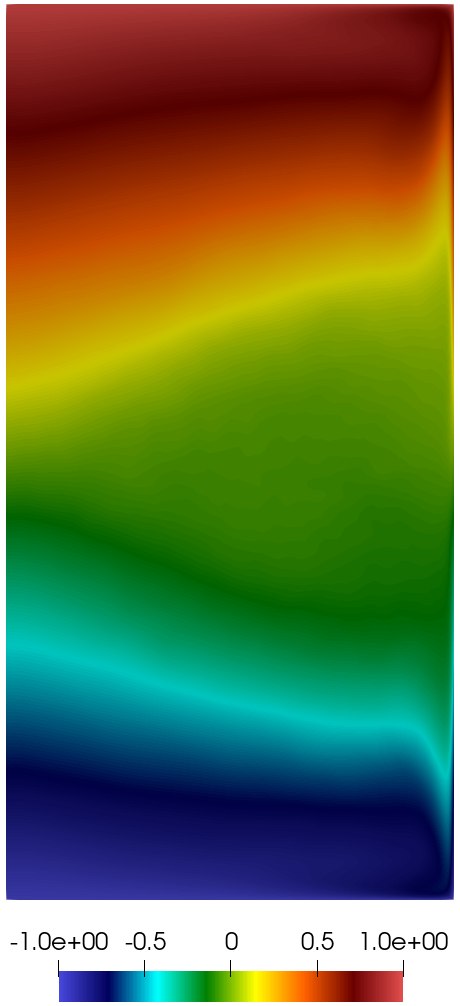} & \hspace{-0.4cm}\includegraphics[align=c,scale = 0.3]{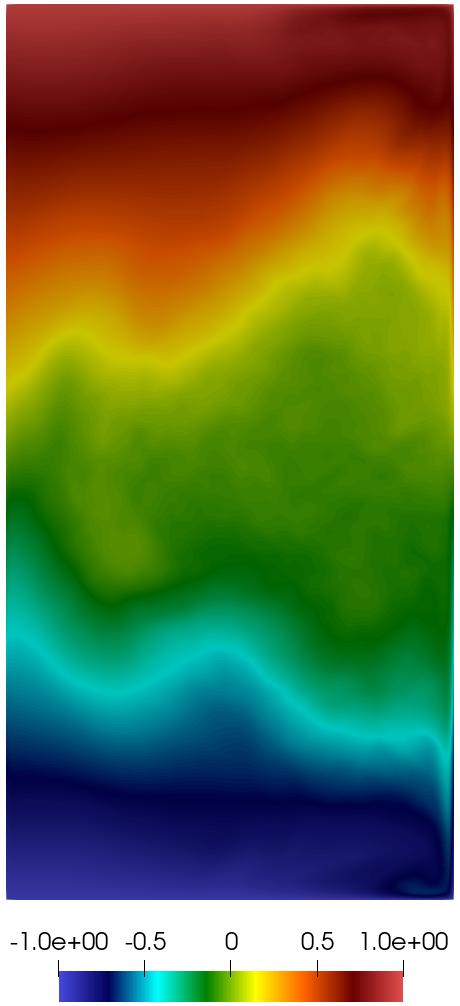} & \hspace{-0.4cm}\includegraphics[align=c,scale = 0.3]{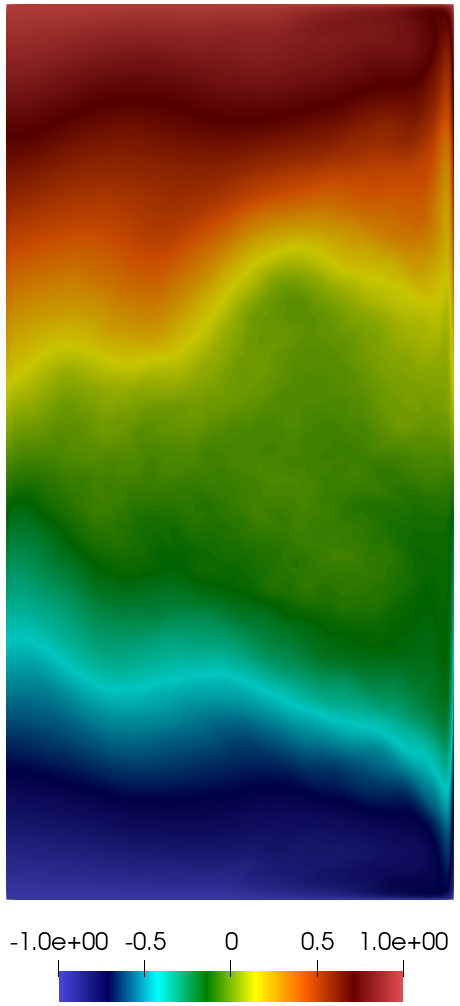} & \hspace{-0.4cm}\includegraphics[align=c,scale = 0.3]{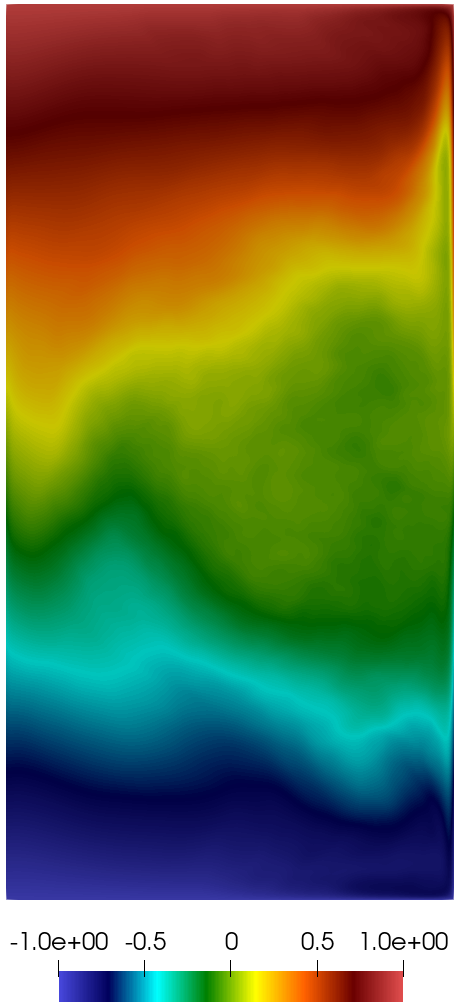} & \hspace{-0.4cm}\includegraphics[align=c,scale = 0.3]{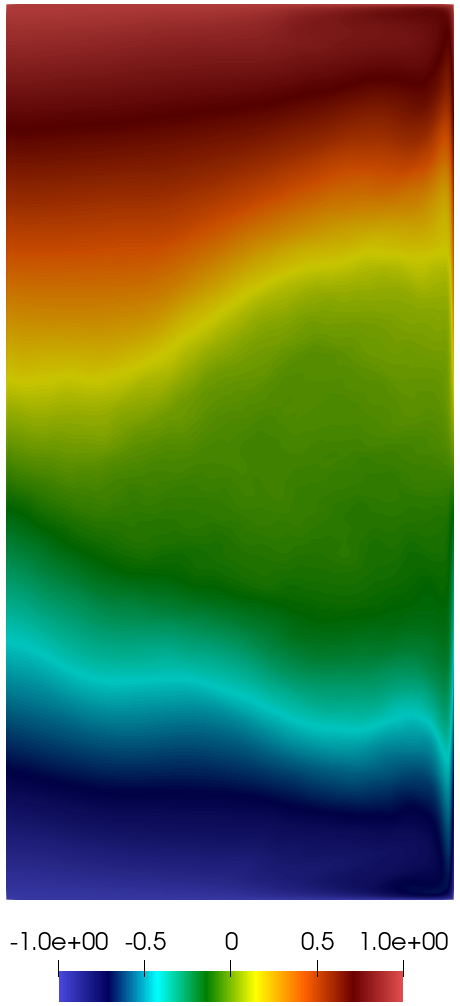} \\
\end{tabular}
\caption{$\widetilde{q}_1$ and $\widetilde{q}_2$ computed over the predictive time interval $[50,100]$ by the FOM (first column) and the POD-LSTM models $\calmf_{q_1}$ (first row) and $\calmf_{q_2}$ (second row) for different values of the lookback window $\sigma_L$ (second to fifth columns).}
\label{fig:q_mu}
\end{figure}

Let us start to assess the predictive performance of the POD-LSTM ROM for the vorticities. 
In order to prioritize computational time savings, we fix the number of modes to $N_{q_1}^r = N_{q_2}^r = 10$, which translates to $14\%$ of cumulative energy for $q_1'$ and $17\%$ for $q_2'$.
Fig. \ref{fig:q_mu} compares the time-averaged potential vorticity $\tildeq_l$, $l = 1, 2$, computed 
over the predictive time interval $[50,100]$
by the FOM (first column) and 
the POD-LSTM ROM with different lookback windows
(second to fifth columns).
The ROM is able to capture the main features of the FOM solution for $\tildeq_1$
for all values of the lookback window. 
However, we see differences
in how the ROM predicts the dark green region (small negative values) around the basin center (i.e., $y=0$), with value $\sigma_L=1$ providing the
most accurate reconstruction of the FOM. 
In the case of $\tildeq_2$, the POD-LSTM ROM struggles
to reconstruct the FOM solution with
$\sigma_L=1,2,5$ which give an incorrect oscillating behavior around the northern and southern regions of the ocean basin. The most accurate prediction of true $\tildeq_2$ is given by $\sigma_L=10$.

To quantify the comparison in Fig.~\ref{fig:q_mu}, Table \ref{tab:error_q_mu} reports the errors $\varepsilon_{q_l}^{(1)}$ \eqref{eq:rmse} and $\varepsilon_{q_l}^{(2)}$ \eqref{eq:l2-error}, for $l = 1,2 $. Both metrics of accuracy confirm that
$\sigma_L=1$ (resp., $\sigma_L=10$) provides the most accurate reconstruction of $\tildeq_1$ (resp., $\tildeq_2$).

\begin{table}[h!]
    \centering
    \begin{tabular}{|c|c|c|c|c|}
        \hline
        $\sigma_L$ & $\varepsilon_{q_1}^{(1)}$ & $\varepsilon_{q_1}^{(2)}$ & $\varepsilon_{q_2}^{(1)}$ & $\varepsilon_{q_2}^{(2)}$ \\
        \hline
        1 & 2.716E-02 & 4.024E-02 & 8.880E-02 & 8.866E-02 \\
        \hline
        2 & 7.837E-02 & 1.162E-01 & 6.902E-02 & 6.894E-02 \\
        \hline
        5 & 4.953E-02 & 7.340E-02 & 8.639E-02 & 8.620E-02 \\
        \hline
        10 & 7.465E-02 & 1.107E-01 & 4.424E-02 & 4.417E-02 \\
        \hline
    \end{tabular}
    \caption{Error metrics $\varepsilon_{q_l}^{(1)}$ \eqref{eq:rmse} and $\varepsilon_{q_l}^{(2)}$ \eqref{eq:l2-error}, $l = 1,2 $, for different
    values of the lookback window $\sigma_L$.}
    \label{tab:error_q_mu}
\end{table}

\begin{figure}[htb!]
    \centering
    \begin{subfigure}[h]{0.45\textwidth}
         \centering
         \includegraphics[width=\textwidth]{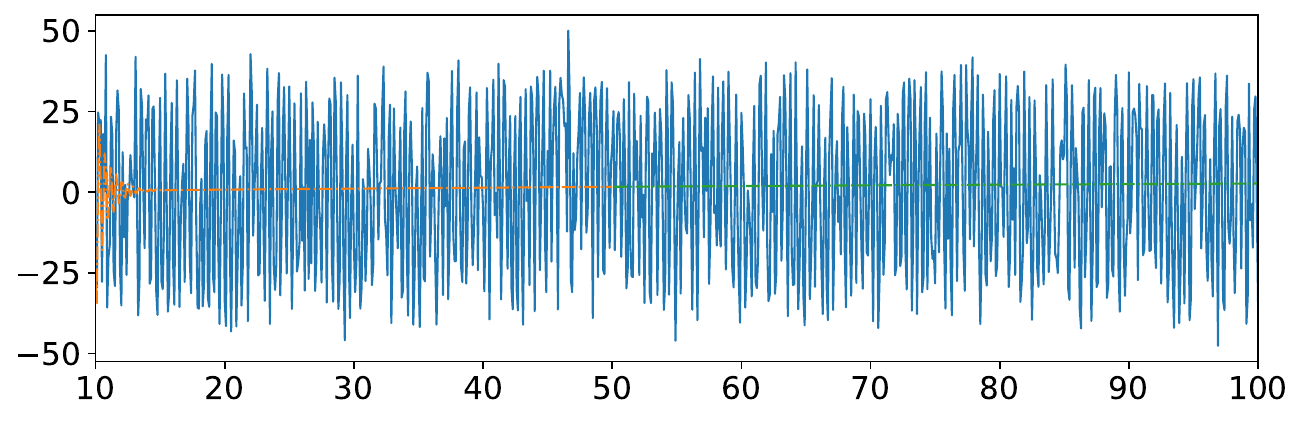}
         \caption{\scriptsize{$\alpha_{1,1}(t),\sigma_L=1$}}
     \end{subfigure}
     \begin{subfigure}[h]{0.45\textwidth}
         \centering
         \includegraphics[width=\textwidth]{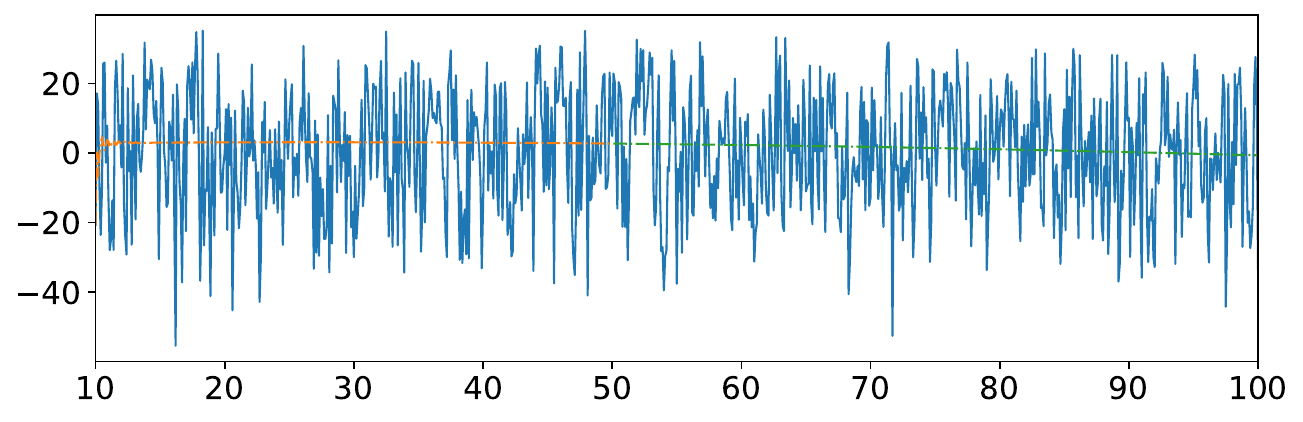}
         \caption{\scriptsize{$\alpha_{1,3}(t),\sigma_L=1$}}
     \end{subfigure}

     \begin{subfigure}[h]{0.45\textwidth}
         \centering
         \includegraphics[width=\textwidth]{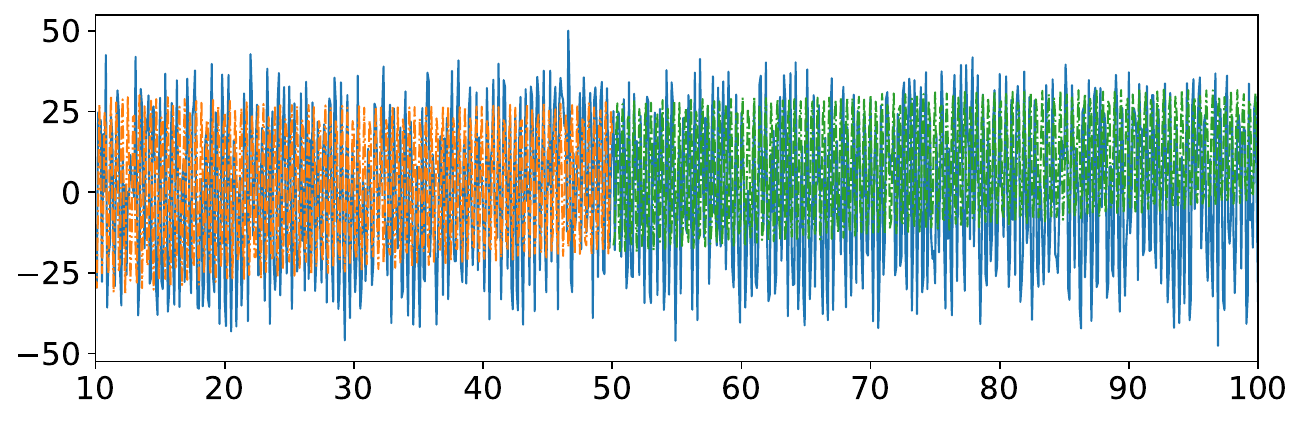}
         \caption{\scriptsize{$\alpha_{1,1}(t),\sigma_L=2$}}
     \end{subfigure}
     \begin{subfigure}[h]{0.45\textwidth}
         \centering
         \includegraphics[width=\textwidth]{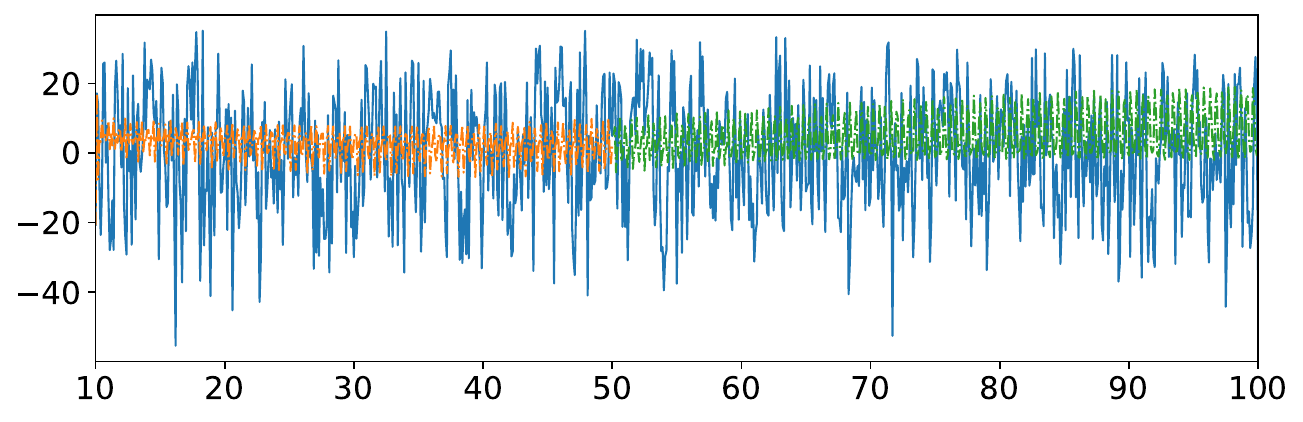}
         \caption{\scriptsize{$\alpha_{1,3}(t),\sigma_L=2$}}
     \end{subfigure}

     \begin{subfigure}[h]{0.45\textwidth}
         \centering
         \includegraphics[width=\textwidth]{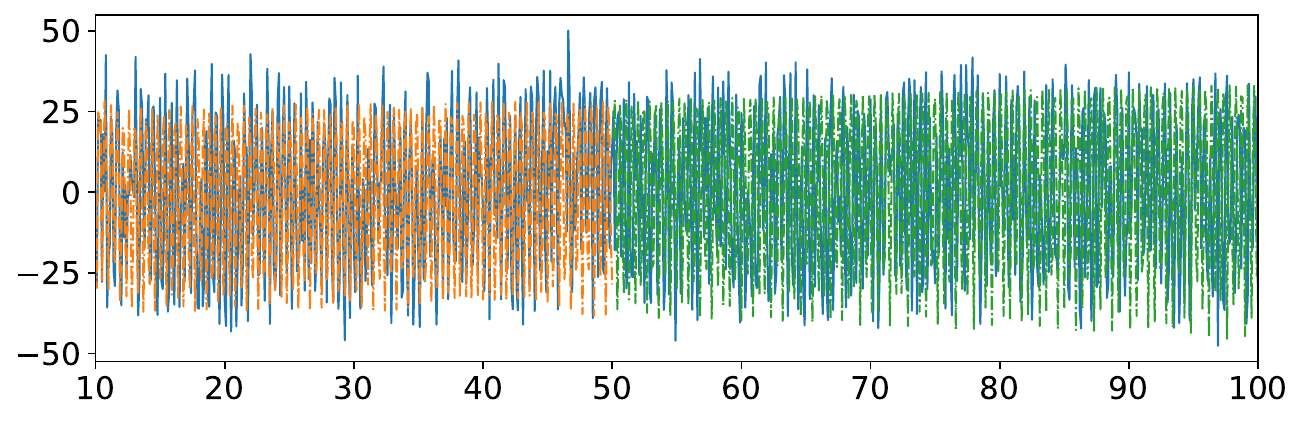}
         \caption{\scriptsize{$\alpha_{1,1}(t),\sigma_L=5$}}
     \end{subfigure}
     \begin{subfigure}[h]{0.45\textwidth}
         \centering
         \includegraphics[width=\textwidth]{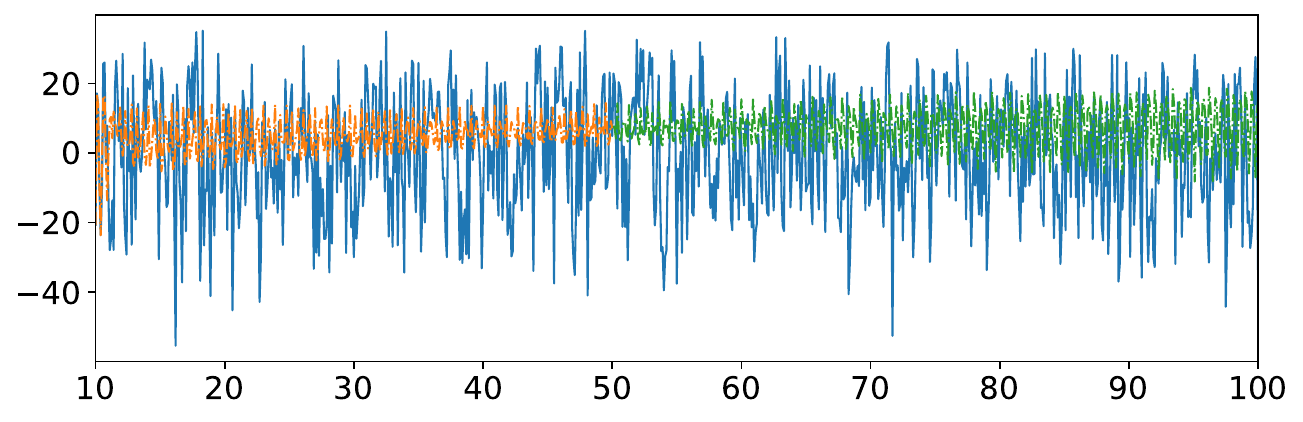}
         \caption{\scriptsize{$\alpha_{1,3}(t),\sigma_L=5$}}
     \end{subfigure}

     \begin{subfigure}[h]{0.45\textwidth}
         \centering
         \includegraphics[width=\textwidth]{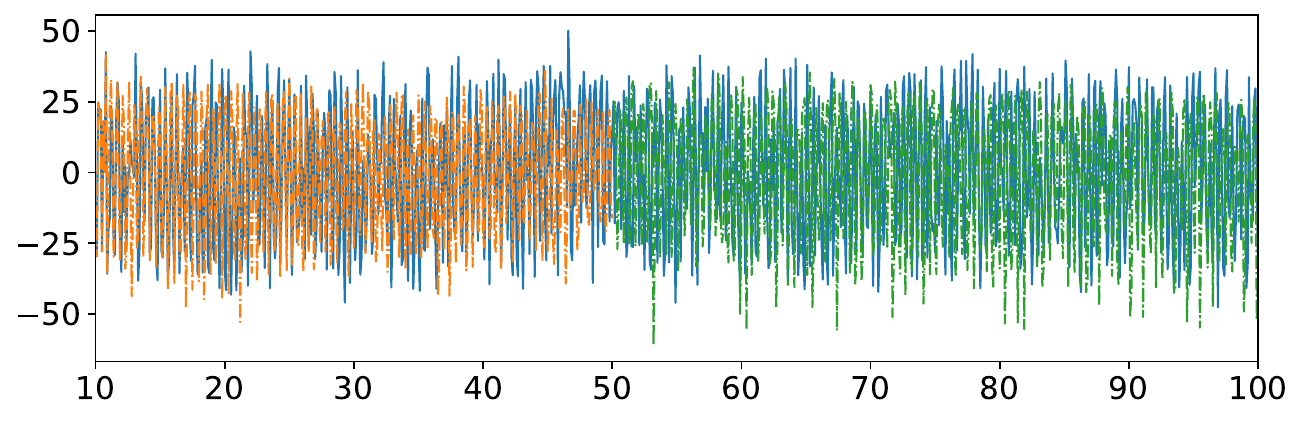}
         \caption{\scriptsize{$\alpha_{1,1}(t),\sigma_L=10$}}
     \end{subfigure}
     \begin{subfigure}[h]{0.45\textwidth}
         \centering
         \includegraphics[width=\textwidth]{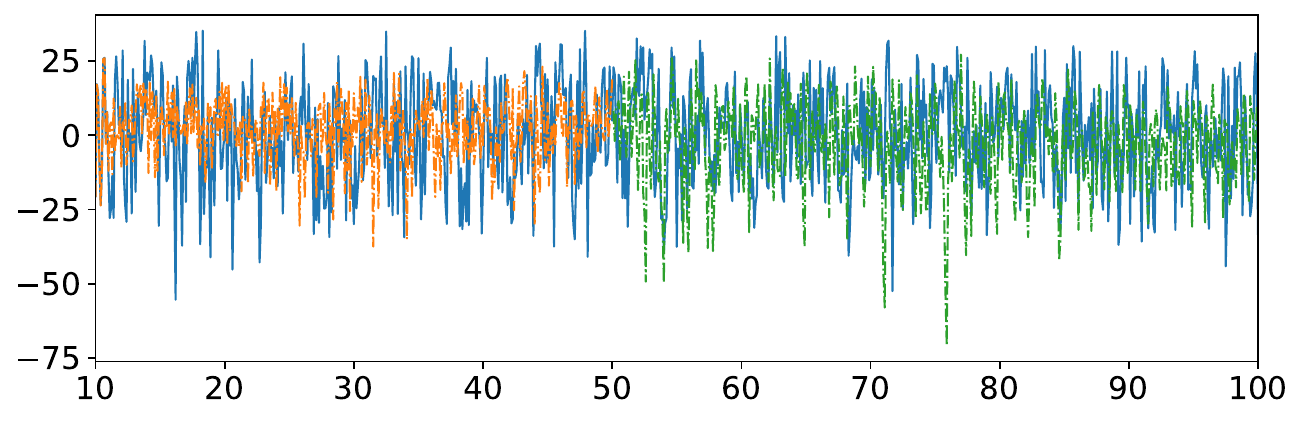}
         \caption{\scriptsize{$\alpha_{1,3}(t),\sigma_L=10$}}
     \end{subfigure}
     \begin{subfigure}[h]{0.45\textwidth}
         \centering
         \vspace{0.3cm}
         \includegraphics[width=\textwidth]{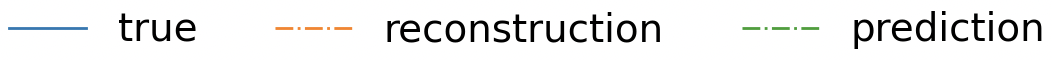}
     \end{subfigure}
\caption{Time evolution of the first $\alpha_{1,1}$ (left) and third $\alpha_{1,3}$ (right) modal coefficients for $q_1$ over the time interval $[10,100]$ for different values of the lookback window $\sigma_L$. 
}
\label{fig:q1-coeff-mu}
\end{figure}

\begin{figure}[htb!]
    \centering
    \begin{subfigure}[h]{0.45\textwidth}
         \centering
         \includegraphics[width=\textwidth]{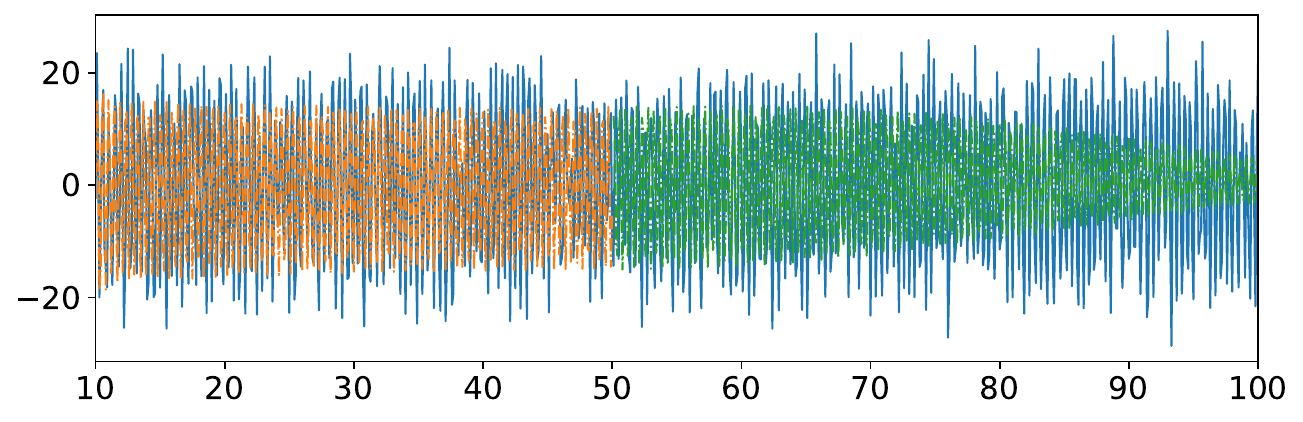}
         \caption{\scriptsize{$\alpha_{2,1}(t),\sigma_L=1$}}
     \end{subfigure}
     \begin{subfigure}[h]{0.45\textwidth}
         \centering
         \includegraphics[width=\textwidth]{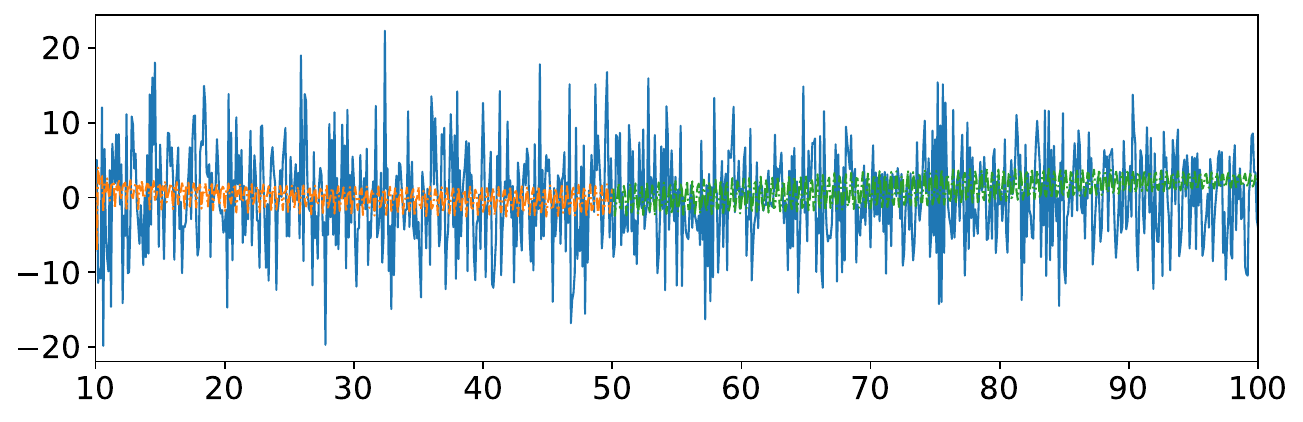}
         \caption{\scriptsize{$\alpha_{2,3}(t),\sigma_L=1$}}
     \end{subfigure}

     \begin{subfigure}[h]{0.45\textwidth}
         \centering
         \includegraphics[width=\textwidth]{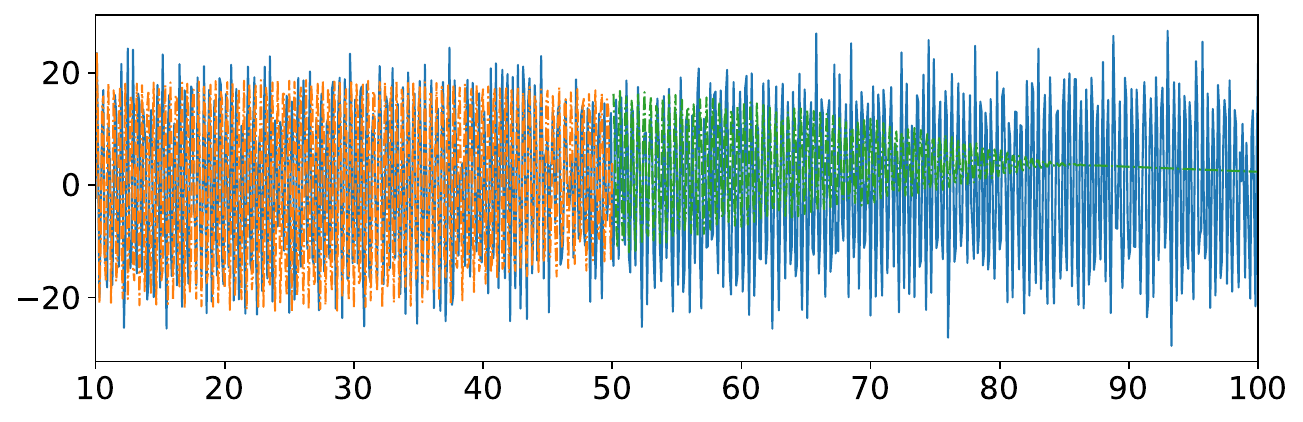}
         \caption{\scriptsize{$\alpha_{2,1}(t),\sigma_L=2$}}
     \end{subfigure}
     \begin{subfigure}[h]{0.45\textwidth}
         \centering
         \includegraphics[width=\textwidth]{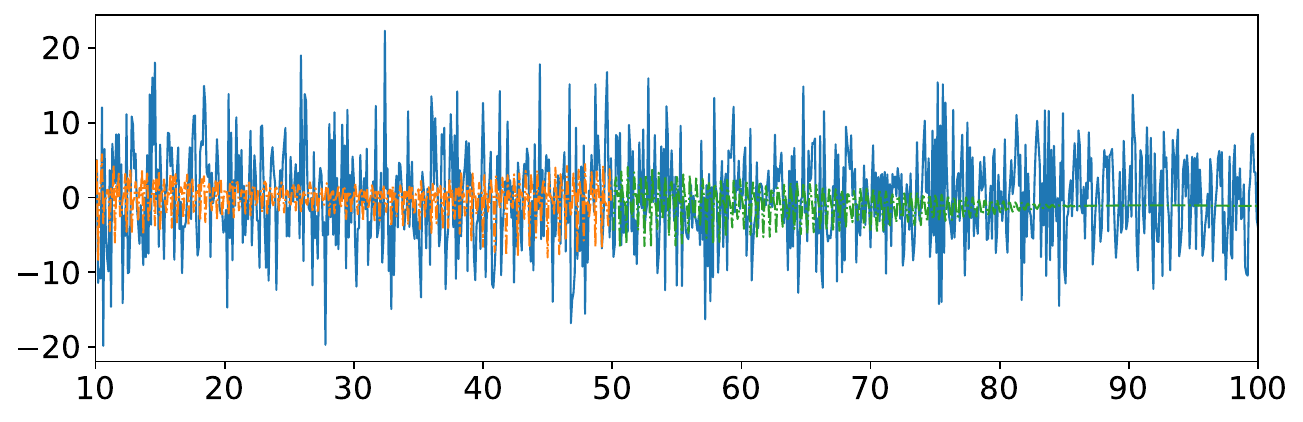}
         \caption{\scriptsize{$\alpha_{2,3}(t),\sigma_L=2$}}
     \end{subfigure}

     \begin{subfigure}[h]{0.45\textwidth}
         \centering
         \includegraphics[width=\textwidth]{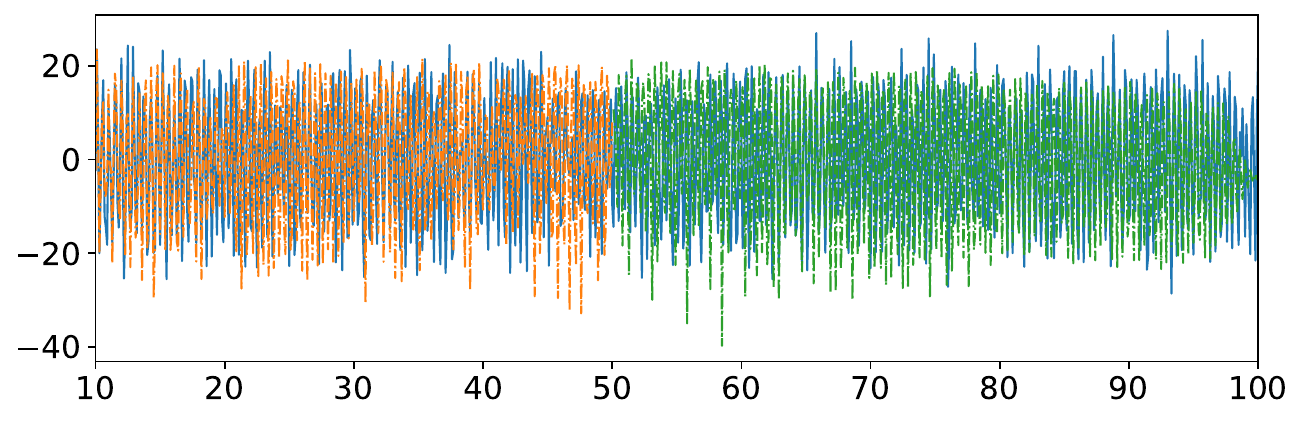}
         \caption{\scriptsize{$\alpha_{2,1}(t),\sigma_L=5$}}
     \end{subfigure}
     \begin{subfigure}[h]{0.45\textwidth}
         \centering
         \includegraphics[width=\textwidth]{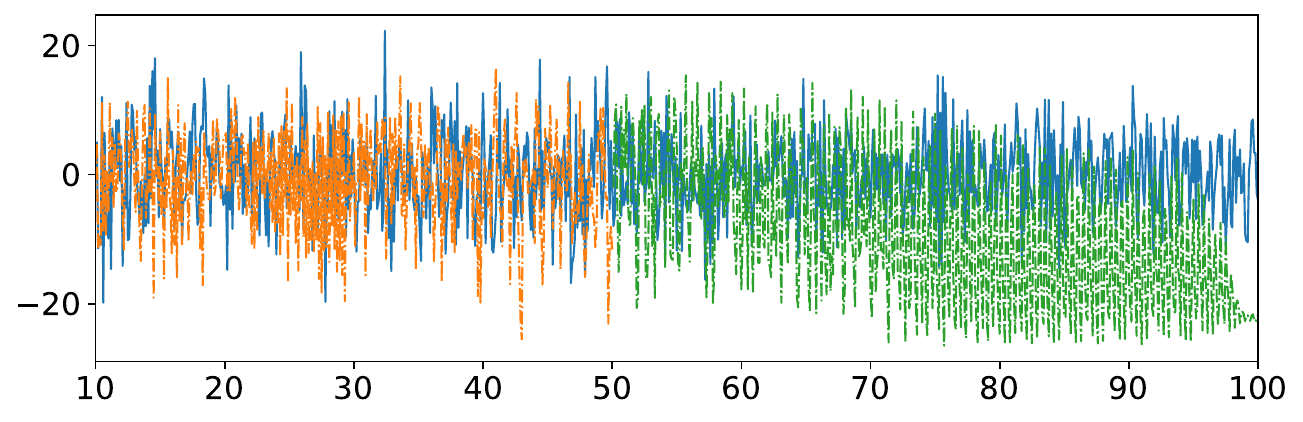}
         \caption{\scriptsize{$\alpha_{2,3}(t),\sigma_L=5$}}
     \end{subfigure}

     \begin{subfigure}[h]{0.45\textwidth}
         \centering
         \includegraphics[width=\textwidth]{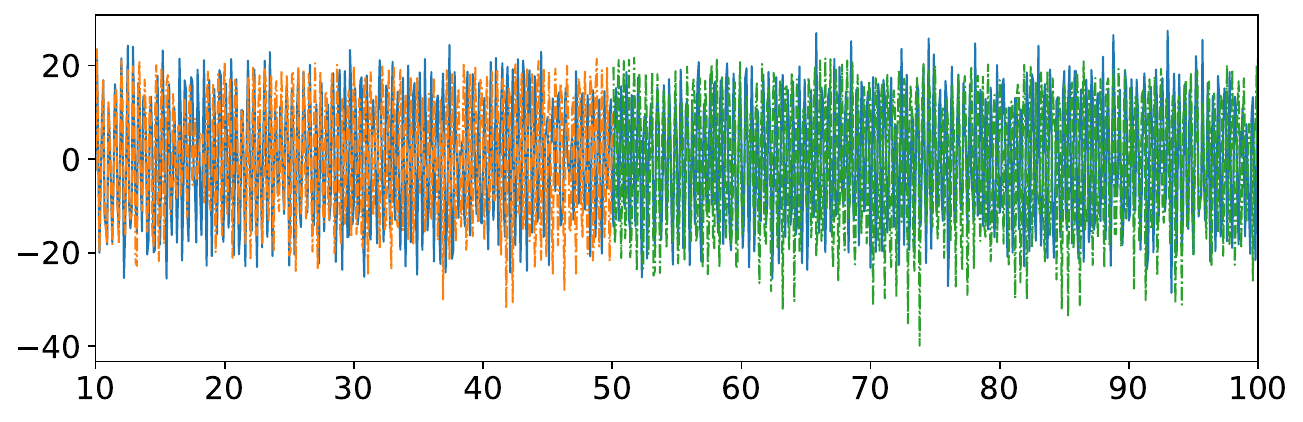}
         \caption{\scriptsize{$\alpha_{2,1}(t),\sigma_L=10$}}
     \end{subfigure}
     \begin{subfigure}[h]{0.45\textwidth}
         \centering
         \includegraphics[width=\textwidth]{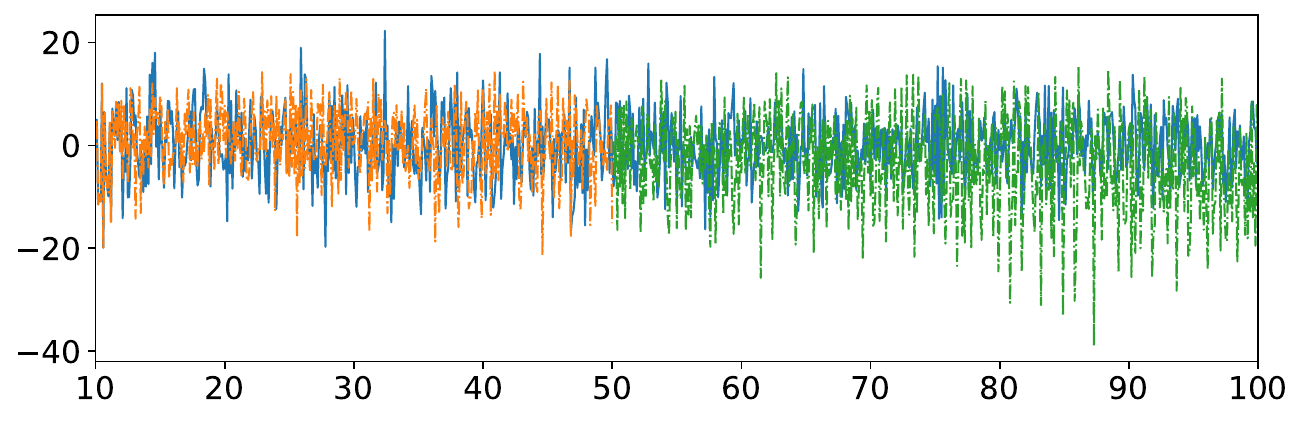}
         \caption{\scriptsize{$\alpha_{2,3}(t),\sigma_L=10$}}
     \end{subfigure}

     \begin{subfigure}[h]{0.45\textwidth}
         \centering
         \vspace{0.3cm}
         \includegraphics[width=\textwidth]{figs/coeff-legend.png}
     \end{subfigure}
\caption{Time evolution of the first $\alpha_{2,1}$ (left) and third $\alpha_{2,3}$ (right) modal coefficients for $q_2$ over the time interval $[10,100]$ for different values of the lookback window $\sigma_L$. 
}
\label{fig:q2-coeff-mu}
\end{figure}

So far, we have shown that $\calmf_{q_1}$ with $\sigma_L=1$ and $\calmf_{q_2}$ with $\sigma_L=10$
give rather accurate predictions of
the time-averaged fields $\tildeq_1$ and $\tildeq_2$, respectively. Next, we investigate how well they  predict the dynamics of the system over the unseen time interval $[50,100]$. Let us start with examining  
the first $\alpha_{l,1}$ and third $\alpha_{l,3}$ modal coefficients, $l = 1,2 $, computed using $\calmf_{q_1}$ and $\calmf_{q_2}$ with different lookback windows $\sigma_L$. 
We have chosen these coefficients, whose evolution in time is shown in in Fig. \ref{fig:q1-coeff-mu}--\ref{fig:q2-coeff-mu}, because they are representative
of all other modal coefficients. 
For $q_1$, we see that, although $\sigma_L=1$ does well in predicting the average, its time-dependent prediction is poor. This is in line with the findings in \cite{Rahman-2019}. 
As $\sigma_L$ increases, we see that both
the reconstruction ($t \in [10,50]$) and the prediction ($t \in (50,100]$) of 
$\alpha_{1,1}$ and $\alpha_{1,3}$
improves, with $\sigma_L=10$ providing the most accurate results for both coefficients. Similar conclusions can be drawn for $q_2$ from Fig.~\ref{fig:q2-coeff-mu}, i.e., the reconstructive and predictive ability of $\sigma_L=1$ is also poor and it tends to improve with the widening of the lookback window.
Nonetheless, some differences can be noted between 
Fig.~\ref{fig:q1-coeff-mu} and \ref{fig:q2-coeff-mu}. The oscillations in $\alpha_{2,1}$ and $\alpha_{2,3}$ computed with $\sigma_L=1$
are more sustained than the oscillations in
$\alpha_{1,1}$ and $\alpha_{1,3}$, while the oscillations in 
$\alpha_{2,1}$ and  $\alpha_{2,3}$ computed with $\sigma_L=2$ die down for $t > 80$. 
Moreover, we see that the average of the oscillations in $\alpha_{2,3}$ computed with $\sigma_L=5$ is off for $t > 70$ with respect
to the true average. We also note that 
$\sigma_L=10$ yields the best prediction of both $q_2$ and its time-averaged quantity $\tildeq_2$, while the values of $\sigma_L$ are different for the first layer: $\sigma_L=1$ for $q_1$ and $\sigma_L=10$ for $\tildeq_1$.

\begin{figure}[htb!]
    \centering
    \begin{subfigure}[h]{0.4\textwidth}
         \centering
         \includegraphics[width=\textwidth]{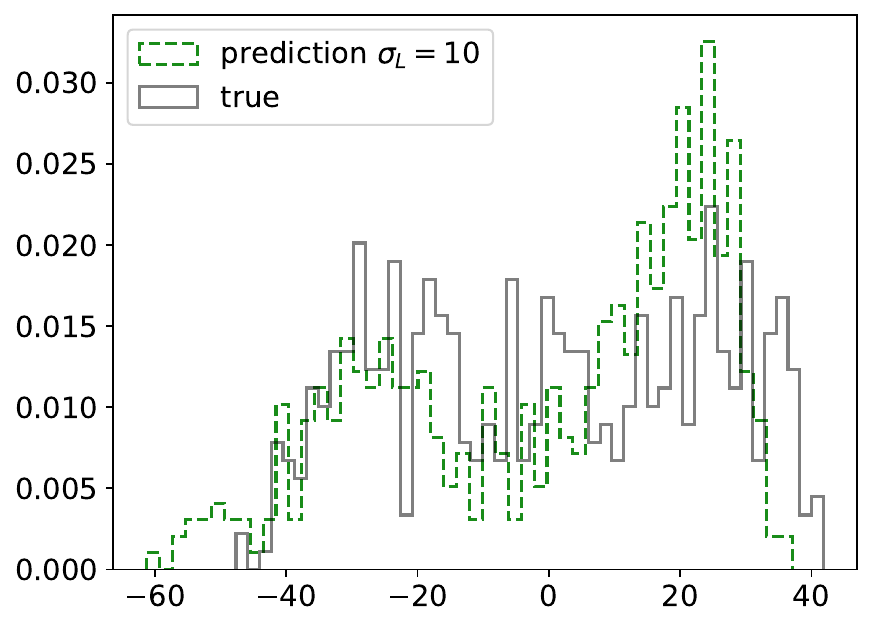}
         \caption{\scriptsize{$\alpha_{1,1}(t)$}}
     \end{subfigure}
     \begin{subfigure}[h]{0.4\textwidth}
         \centering
         \includegraphics[width=\textwidth]{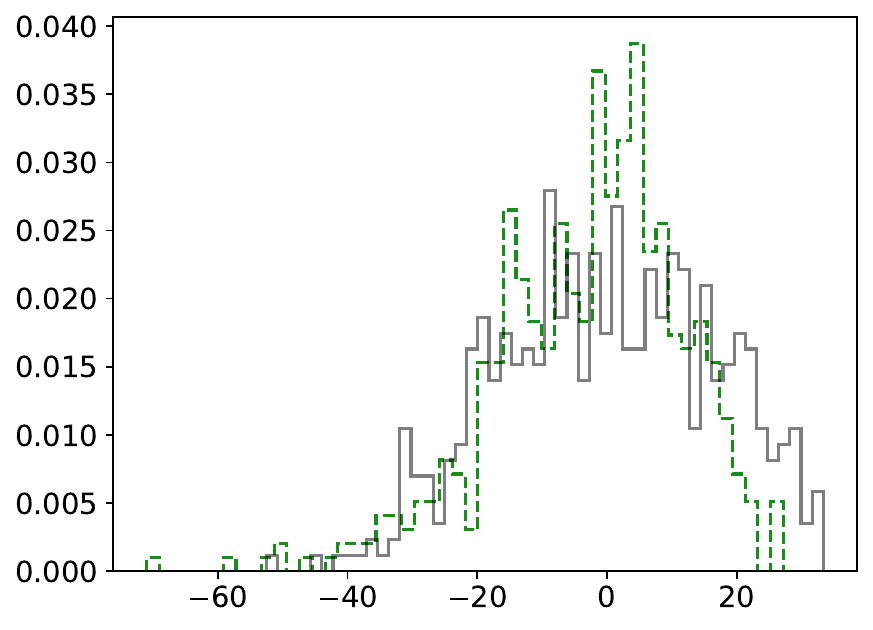}
         \caption{\scriptsize{$\alpha_{1,3}(t)$}}
     \end{subfigure}

     \begin{subfigure}[h]{0.4\textwidth}
         \centering
         \includegraphics[width=\textwidth]{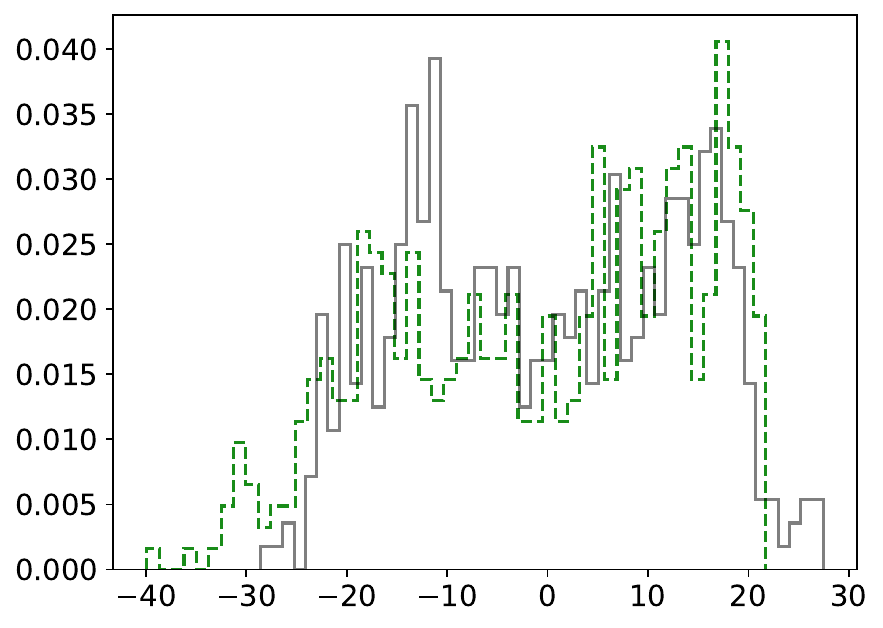}
         \caption{\scriptsize{$\alpha_{2,1}(t)$}}
     \end{subfigure}
     \begin{subfigure}[h]{0.4\textwidth}
         \centering
         \includegraphics[width=\textwidth]{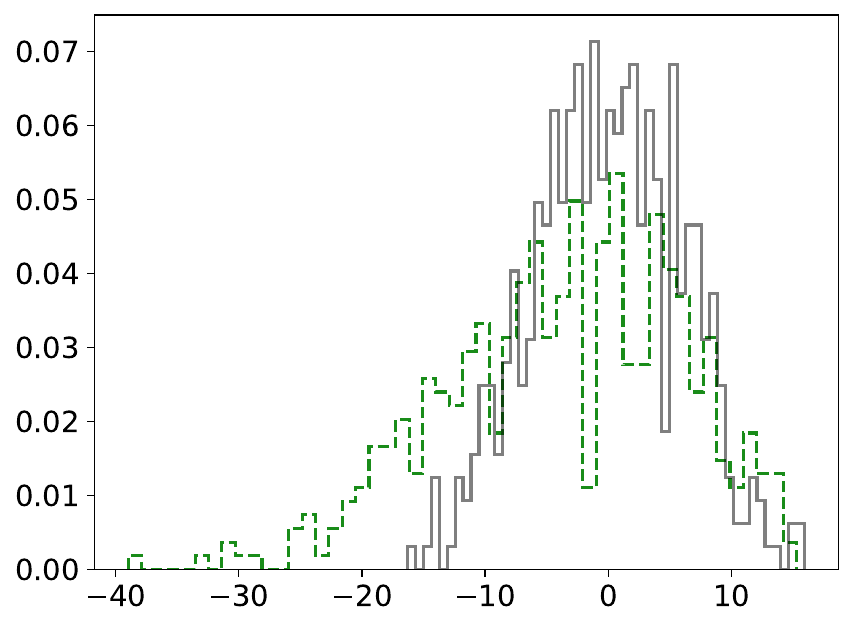}
         \caption{\scriptsize{$\alpha_{2,3}(t)$}}
     \end{subfigure}
\caption{Probability mass functions of the first and third modal coefficients of $q_1$ (top row) and $q_2$ (bottom row) over time interval $[50,100]$: 
comparison between true value, i.e.,  
for the FOM solution projected onto the appropriate reduced space, and the value computed by $\calmf_{q_1}$ and $\calmf_{q_2}$ for $\sigma_L=10$.
}
\label{fig:q-coeff_hist}
\end{figure}

\begin{figure}[htb!]
    \centering
    \begin{subfigure}[h]{0.47\textwidth}
         \centering
         \includegraphics[width=\textwidth]{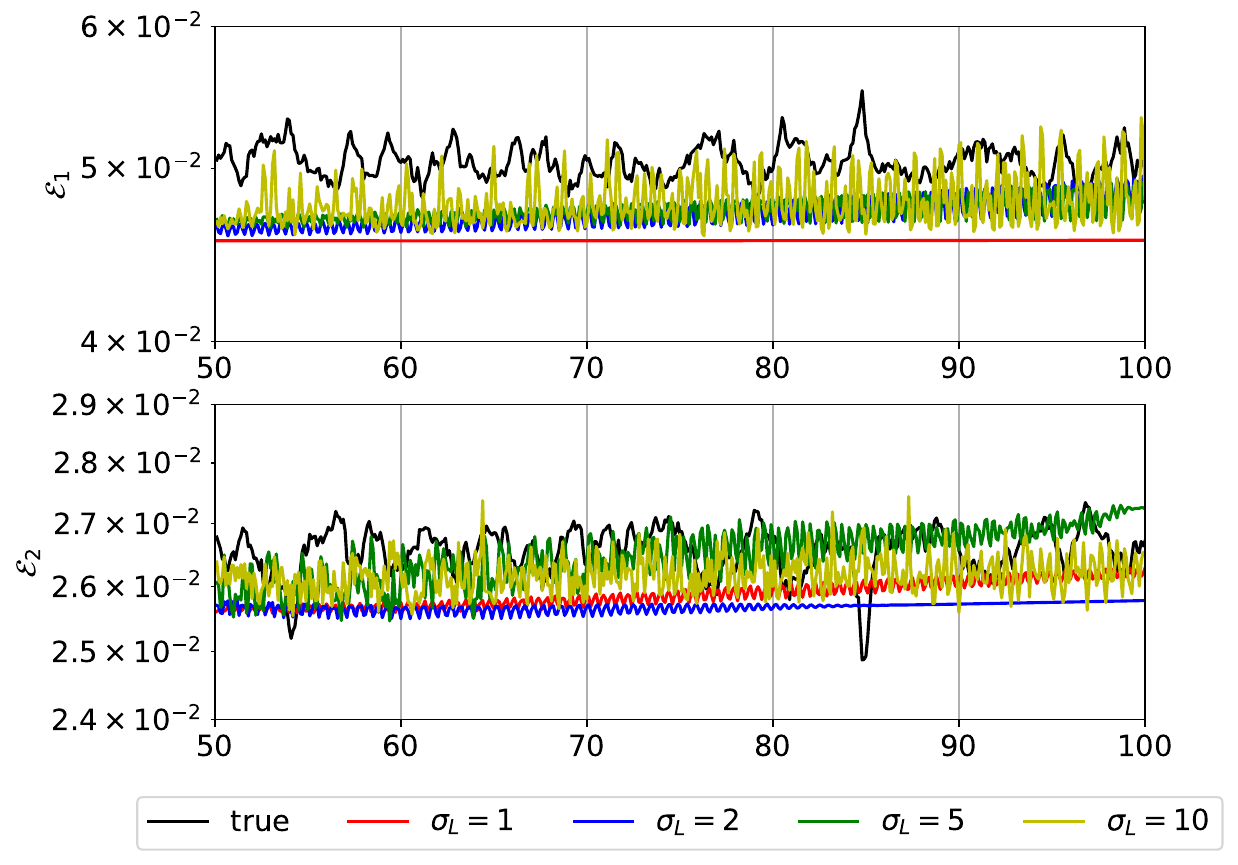}
     \end{subfigure}
     \begin{subfigure}[h]{0.47\textwidth}
         \centering
         \includegraphics[width=\textwidth]{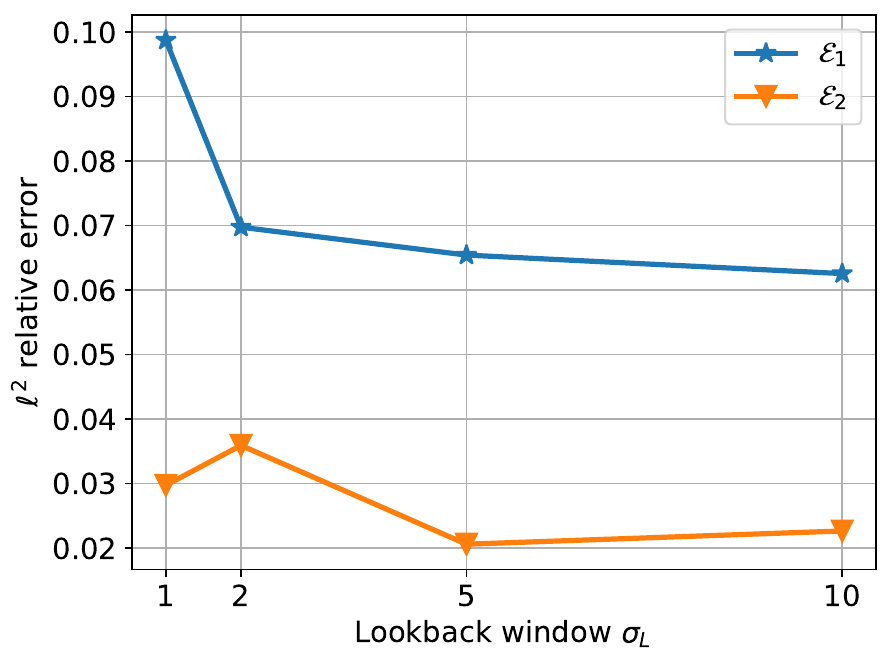}
     \end{subfigure}
\caption{
Time evolution of the enstrophy for the top layer $\mathcal{E}_1$ (top left) and the bottom layer $\mathcal{E}_2$ (bottom left) computed using the FOM and the POD-LSTM ROM with different lookback window $\sigma_L$. 
The corresponding relative $L^2$ error is shown on the right.}
\label{fig:ens-mu}
\end{figure}

Fig.~\ref{fig:q-coeff_hist}
compares the probability mass functions (PMF), i.e., a function that gives the probability of each possible value for a discrete random variable,
of the first and third modal coefficients 
over the unseen time interval $[50,100]$
for the  FOM solution projected onto the appropriate reduced space with the counterparts computed by $\calmf_{q_1}$ and $\calmf_{q_2}$ for 
$\sigma_L=10$. We observe a good agreement.
The comparison in terms of enstrophy evolution \eqref{eq:enstrophy} is shown in Fig.~\ref{fig:ens-mu} (left) for different values of $\sigma_L$.
For the top layer, $\sigma_L=1$ gives the worst prediction, which is consistent with our observation from Fig. \ref{fig:q1-coeff-mu}, 
while the best prediction (i.e., lowest relative $L^2$ error) is for $\sigma_L=10$. 
For $\sigma_L=10$, we obtain a good prediction 
also for the enstrophy in the bottom layer, although not the best in the relative $L^2$ norm.
See Fig.~\ref{fig:ens-mu} (right).

\begin{figure}[htb!]
\centering
\begin{tabular}{cccccc}
       & FOM & \hspace{-0.4cm} $N_{q_l}^r = 2$ & \hspace{-0.4cm} $N_{q_l}^r = 4$ & \hspace{-0.4cm} $N_{q_l}^r = 8$ & \hspace{-0.4cm} $N_{q_l}^r = 10$ \\
    $\tildeq_1$ & \includegraphics[align=c,scale = 0.3]{figs/q1/q1Mean_true.png} & \hspace{-0.4cm}\includegraphics[align=c,scale = 0.3]{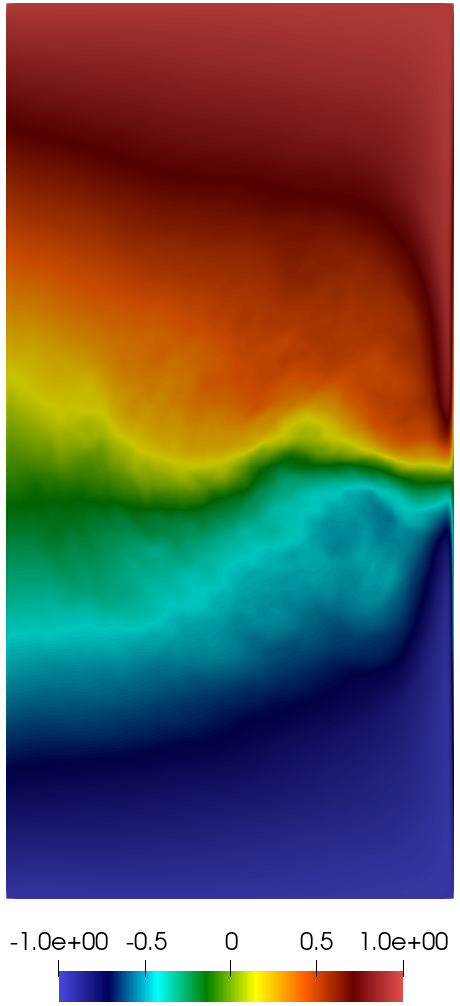} & \hspace{-0.4cm}\includegraphics[align=c,scale = 0.3]{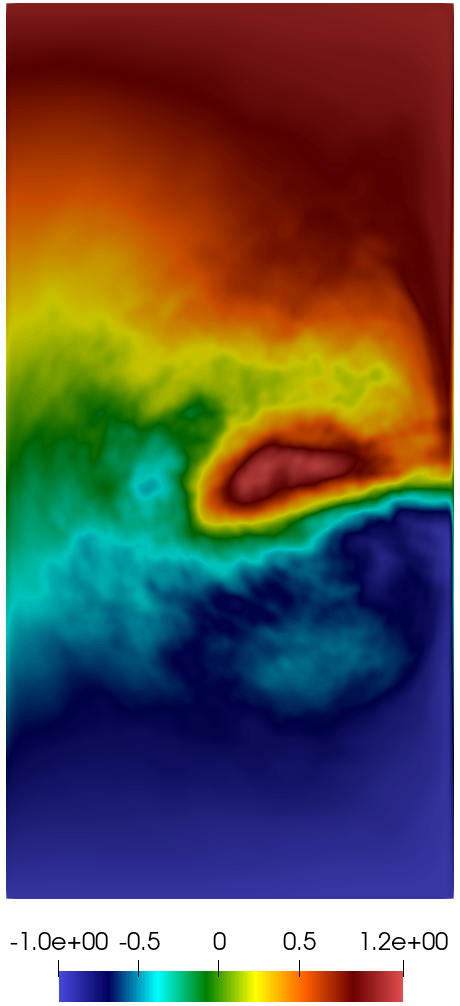} & \hspace{-0.4cm}\includegraphics[align=c,scale = 0.3]{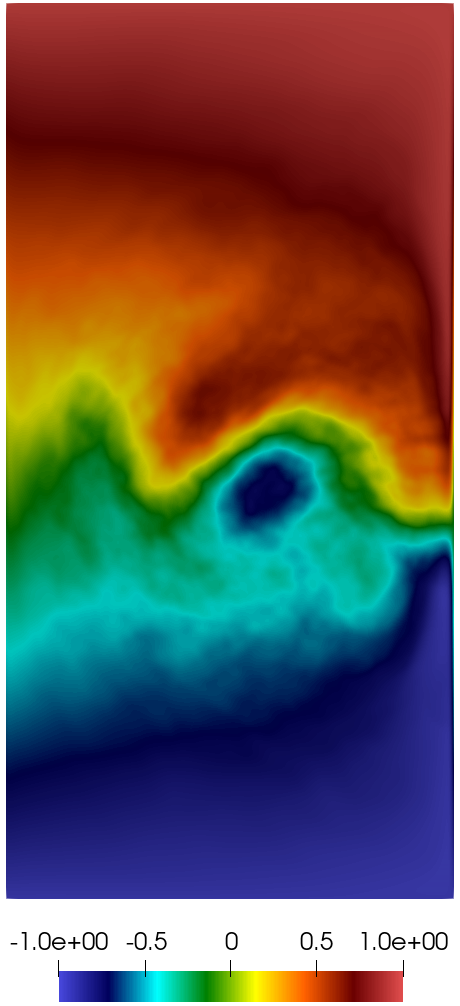} & \hspace{-0.4cm}\includegraphics[align=c,scale = 0.3]{figs/q1/q1Mean_mu10.png} \\
    $\tildeq_2$ & \includegraphics[align=c,scale = 0.3]{figs/q2/q2Mean_true.png} & \hspace{-0.4cm}\includegraphics[align=c,scale = 0.3]{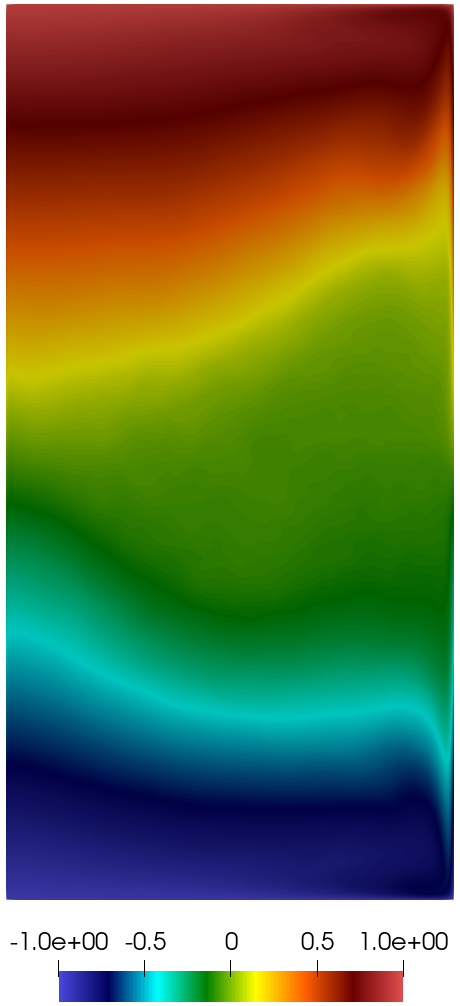} & \hspace{-0.4cm}\includegraphics[align=c,scale = 0.3]{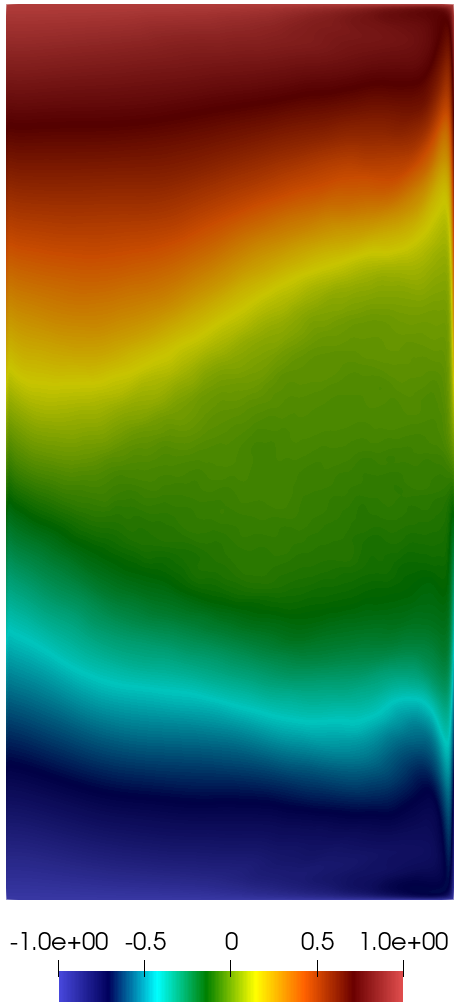} & \hspace{-0.4cm}\includegraphics[align=c,scale = 0.3]{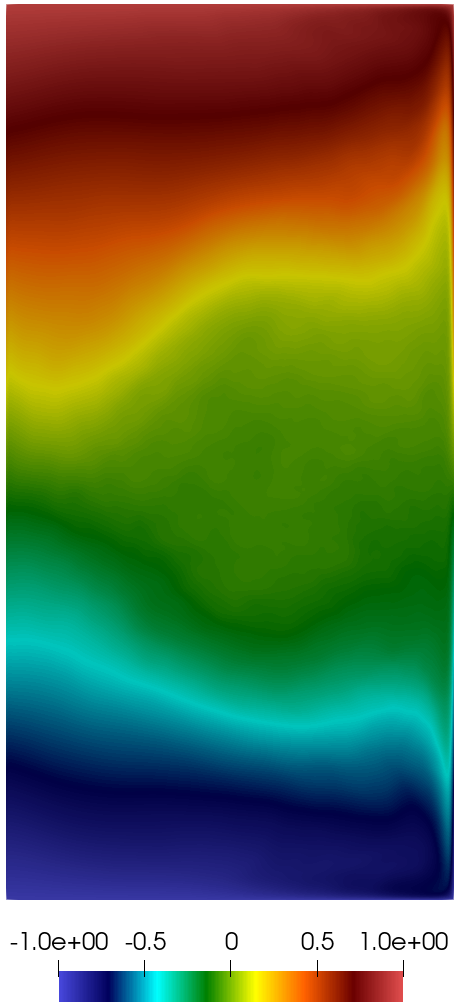} & \hspace{-0.4cm}\includegraphics[align=c,scale = 0.3]{figs/q2/q2Mean_mu10.png} \\
\end{tabular}
\caption{Time-averaged vorticities $\widetilde{q}_1$ (top) and $\widetilde{q}_2$ (bottom) over the predictive time interval $[50,100]$ computed by FOM (first column) and POD-LSTM ROM with $N_{q_l}^r=2,4,8,10$ (second to fifth columns).
}
\label{fig:q_mod}
\end{figure}

In light of the results presented so far, 
we fix the lookback window to $\sigma_L=10$ for both $\calmf_{q_1}$ and $\calmf_{q_2}$. 
Next, we vary the number of modes retained. Fig.~\ref{fig:q_mod} shows the time-averaged potential vorticities $\tildeq_1$ and $\tildeq_2$ computed using the FOM and the POD-LSTM ROM with different numbers of modal basis functions retained. 
We see that $\tildeq_1$ given by the ROM with $N_{q_1}^r = 4,8$ cannot correctly capture the solution around the center of the basin.
The ROM with $N_{q_1}^r = 10$ improves the solution, however the most accurate ROM solution
is obtained with $N_{q_1}^r = 2$. 
It could be argued that also for $\tildeq_2$
the most accurate ROM solution is obtained with 
$N_{q_2}^r = 2$. This is rather remarkable since
with $N_{q_1}^r = 2$ (resp., $N_{q_2}^r = 2$)
we retain only 5\% (resp., 6\%) of the cumulative eigenvalue energy \eqref{eq:energy}. 
To complement the qualitative assessment in 
Fig.~\ref{fig:q_mod}, 
Table \ref{tab:error-q-mod} reports the errors $\varepsilon_{q_l}^{(1)}$ \eqref{eq:rmse} and $\varepsilon_{q_l}^{(2)}$ \eqref{eq:l2-error}, for $l = 1,2 $, and the relative $L^2$ error
of the enstrophy in both layers. 

\begin{figure}[htb!]
    \centering
    \begin{subfigure}[h]{0.45\textwidth}
         \centering
         \includegraphics[width=\textwidth]{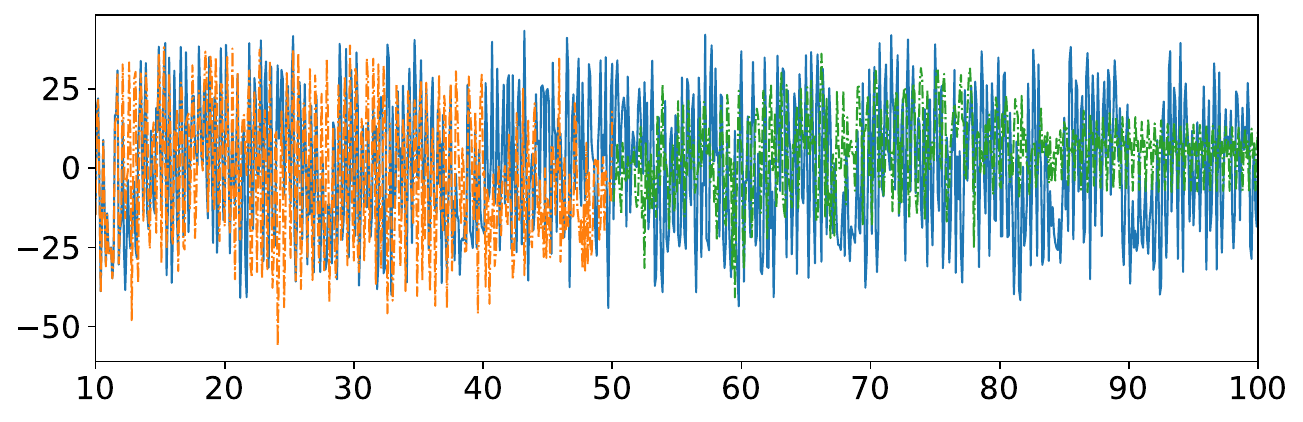}
         \caption{\scriptsize{$\alpha_{1,2}(t),N_{q_1}^r=2$}}
     \end{subfigure}
     \begin{subfigure}[h]{0.45\textwidth}
         \centering
         \includegraphics[width=\textwidth]{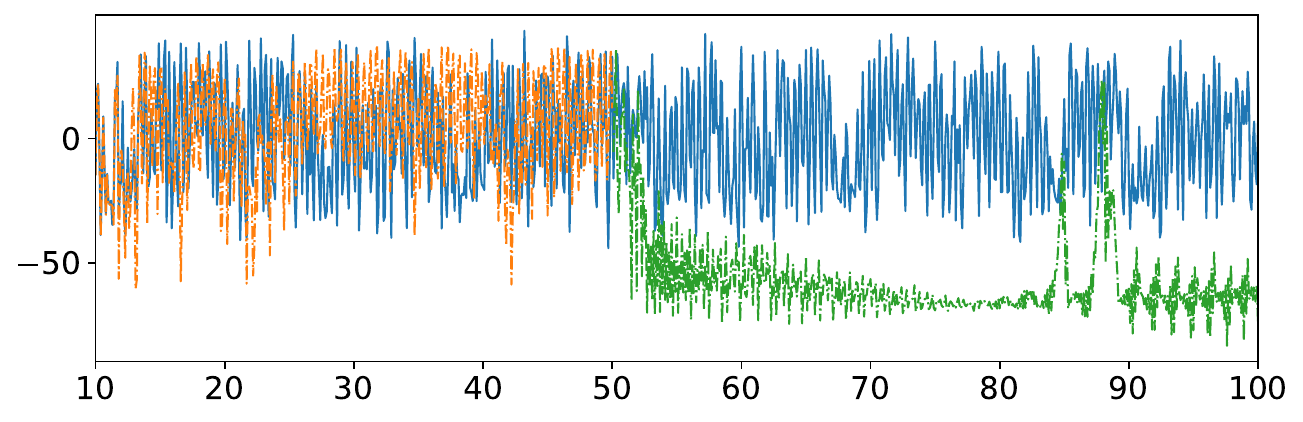}
         \caption{\scriptsize{$\alpha_{1,2}(t),N_{q_1}^r=4$}}
     \end{subfigure}
    \begin{subfigure}[h]{0.45\textwidth}
         \centering
         \includegraphics[width=\textwidth]{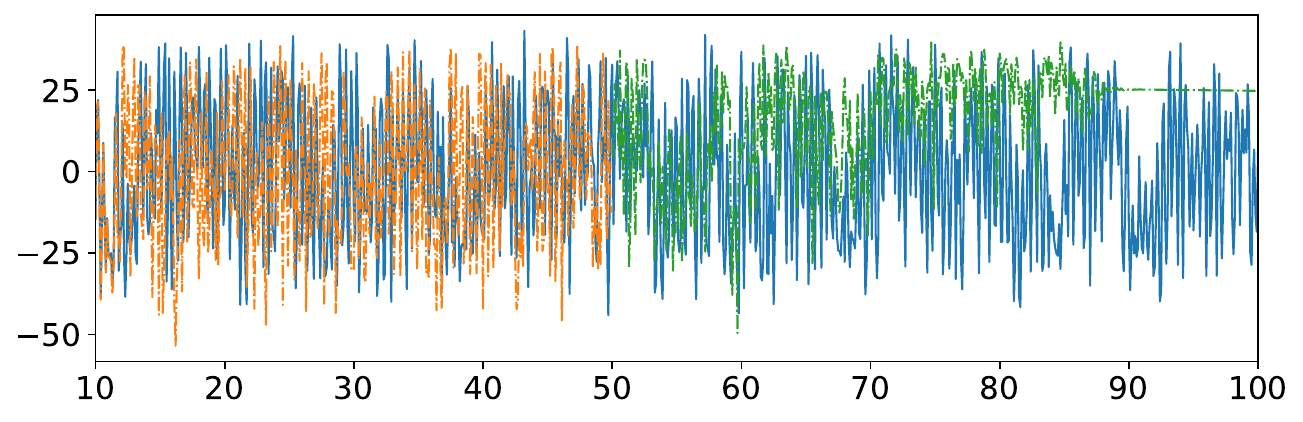}
         \caption{\scriptsize{$\alpha_{1,2}(t),N_{q_1}^r=8$}}
     \end{subfigure}
     \begin{subfigure}[h]{0.45\textwidth}
         \centering
         \includegraphics[width=\textwidth]{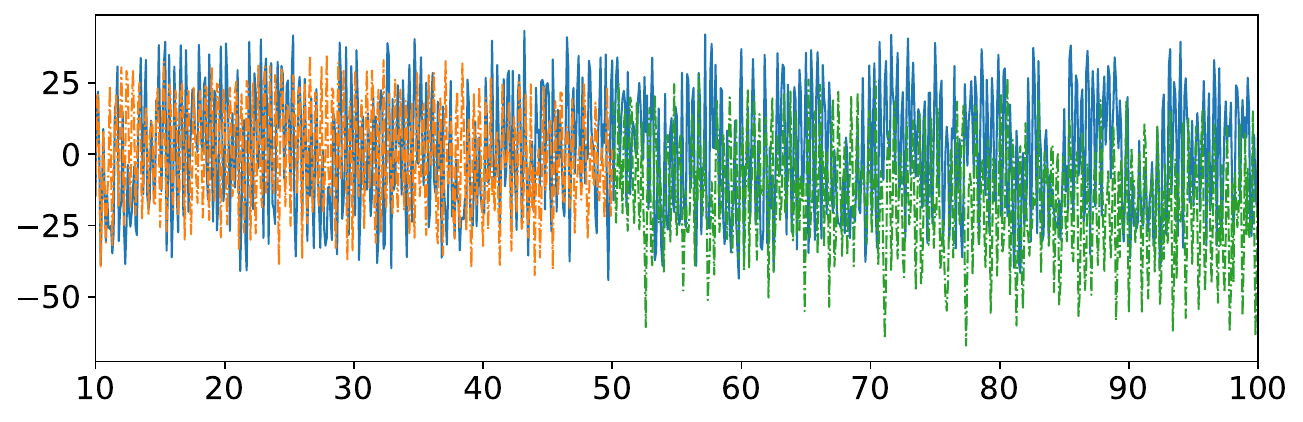}
         \caption{\scriptsize{$\alpha_{1,2}(t),N_{q_1}^r=10$}}
     \end{subfigure}

    \begin{subfigure}[h]{0.45\textwidth}
         \centering
         \includegraphics[width=\textwidth]{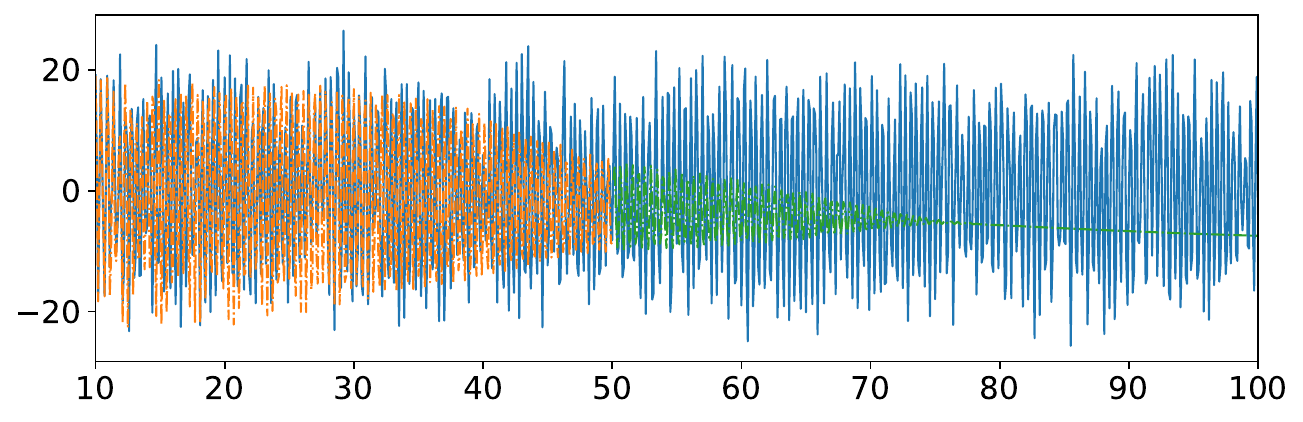}
         \caption{\scriptsize{$\alpha_{2,2}(t),N_{q_2}^r=2$}}
    \end{subfigure}
    \begin{subfigure}[h]{0.45\textwidth}
         \centering
         \includegraphics[width=\textwidth]{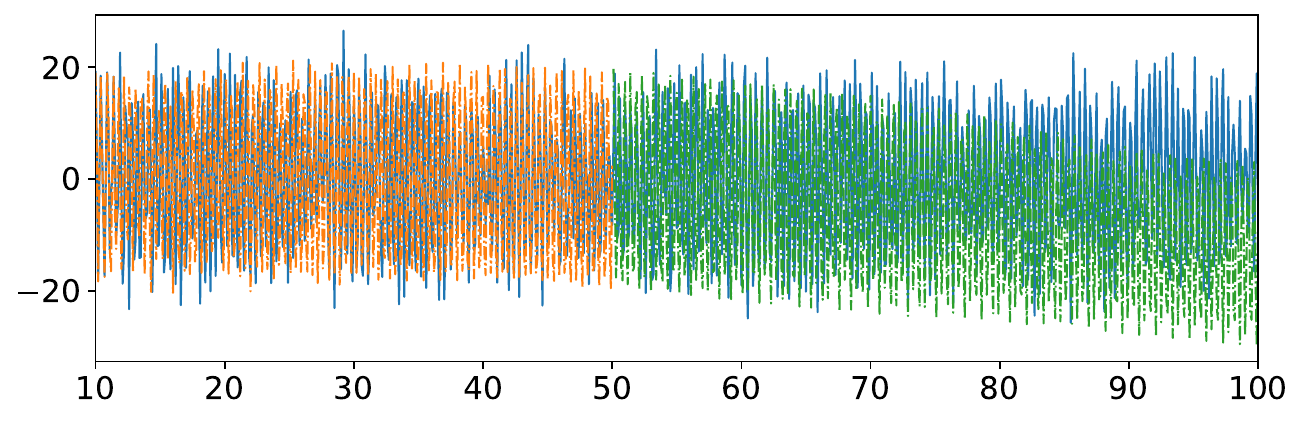}
         \caption{\scriptsize{$\alpha_{2,2}(t),N_{q_2}^r=4$}}
    \end{subfigure}
    \begin{subfigure}[h]{0.45\textwidth}
         \centering
         \includegraphics[width=\textwidth]{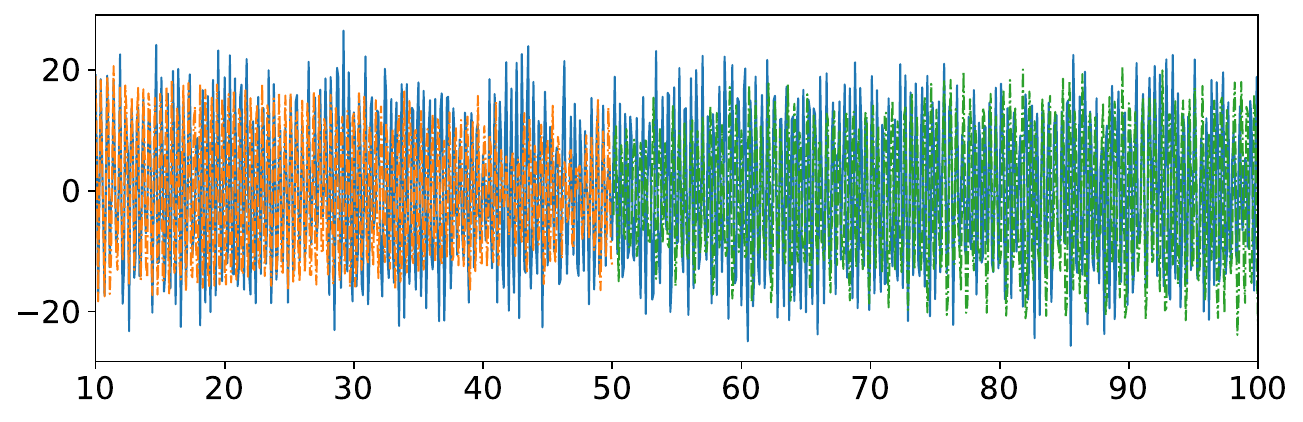}
         \caption{\scriptsize{$\alpha_{2,2}(t),N_{q_2}^r=8$}}
    \end{subfigure}
    \begin{subfigure}[h]{0.45\textwidth}
         \centering
         \includegraphics[width=\textwidth]{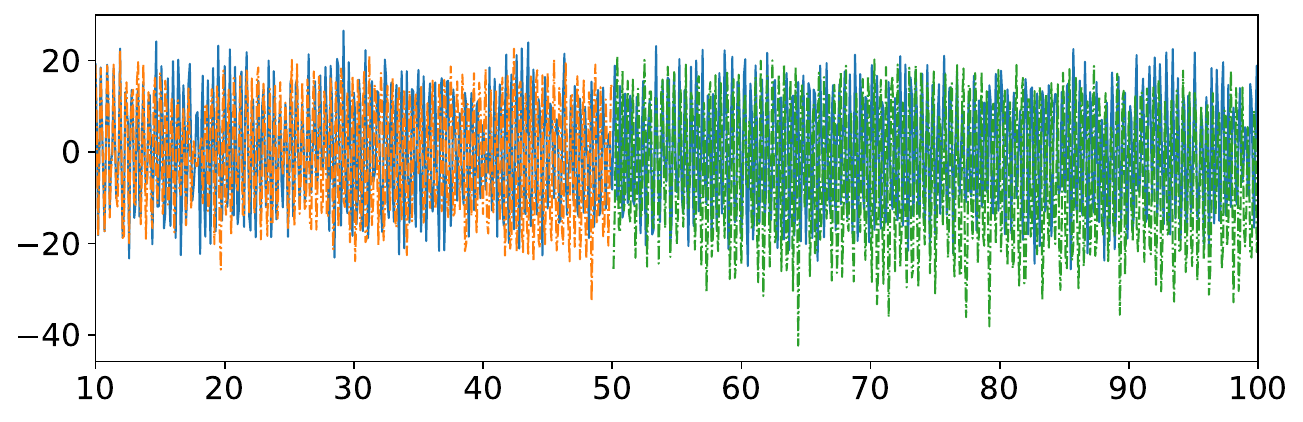}
         \caption{\scriptsize{$\alpha_{2,2}(t),N_{q_2}^r=10$}}
    \end{subfigure}

    \begin{subfigure}[h]{0.45\textwidth}
         \centering
         \vspace{0.3cm}
         \includegraphics[width=\textwidth]{figs/coeff-legend.png}
     \end{subfigure}
\caption{Time evolution of the second modal coefficient $\alpha_{l,2}$ over the time interval $[10,100]$ for different 
numbers of retained POD basis function $N_{q_l}^r$, $l = 1,2 $.}
\label{fig:q-coeff-mod}
\end{figure}

Next, we look at the evolution in time of 
second $\alpha_{l,2}$ modal coefficient computed using $\calmf_{q_l}$ with different
$N_{q_l}^r$, $l = 1,2$. 
This evolution, shown in Fig.~\ref{fig:q-coeff-mod}, is chosen to be representative of the evolutions of all other modal coefficients. 
We see that, even though $N_{q_1}^r=2$ gives lower errors in the metrics in Table~\ref{tab:error-q-mod}, both the reconstruction and prediction of $\alpha_{1,2}$ are inferior to those obtained with $N_{q_1}^r=10$. 
Compare Fig.~\ref{fig:q-coeff-mod} (a) and (d).
This, combined with the  fact that the $L^2$ relative error of $\mathcal{E}_1$ for $N_{q_1}^r=10$ is about 1.5 times smaller than the error $N_{q_1}^r=2$ (see Table~\ref{tab:error-q-mod}), leads us to believe that $\calmf_{q_1}$ with $N_{q_1}^r=10$ can 
predict time-dependent quantities with more accuracy than $\calmf_{q_1}$ with $N_{q_1}^r=2$. 
Further, Table \ref{tab:error-q-mod} and Fig. \ref{fig:q-coeff-mod} (e)-(h) indicate that
$\calmf_{q_2}$ with $N_{q_2}^r=10$ results in the most accurate prediction of the time-averaged field $\tildeq_2$, time evolution reconstruction of $\mathcal{E}_2$, and the second modal coefficient $\alpha_{2,2}$.

\begin{table}[h!]
    \centering
    \begin{tabular}{|c|cccc|}
    \hline
        $N_{q_1}^r$ & $\delta_{q_1}$ & $\varepsilon_{q_1}^{(1)}$ & $\varepsilon_{q_1}^{(2)}$ & $\mathcal{E}_1$ rel. err.\\
        \hline
        2 & 0.05  & 3.999E-02 & 5.922E-02 & 9.288E-02 \\
        \hline
        4 & 0.07 & 1.824E-01 & 2.705E-01 & 4.528E-02 \\
        \hline
        8 & 0.12 & 9.429E-02 & 1.398E-01 & 5.705E-02 \\
        \hline
        10 & 0.13 & 7.465E-02 & 1.107E-01 & 6.256E-02 \\
        \hline
         $N_{q_2}^r$ & $\delta_{q_2}$ & $\varepsilon_{q_2}^{(1)}$ & $\varepsilon_{q_2}^{(2)}$ & $\mathcal{E}_2$ rel. err.\\
        \hline
        2 & 0.06 & 2.428E-02 & 4.798E-02 & 3.880E-02 \\
        \hline
        4 & 0.09 & 2.727E-02 & 5.385E-02 & 3.006E-02 \\
        \hline
        8 & 0.14 & 2.605E-02 & 5.146E-02 & 2.831E-01 \\
        \hline
        10 & 0.16 & 2.235E-02 & 4.417E-02 & 2.266E-02 \\
        \hline
    \end{tabular}
    \caption{Fraction of retained eigenvalue energy $\delta_{q_l}$, 
    error metrics $\varepsilon_{q_l}^{(1)}$ \eqref{eq:rmse} and $\varepsilon_{q_l}^{(2)}$ \eqref{eq:l2-error}, 
    and relative $L^2$ error for enstrophy  $\mathcal{E}_l$ \eqref{eq:enstrophy},
    $l = 1,2 $, for different numbers of POD modes retained.}
    \label{tab:error-q-mod}
\end{table}


\begin{figure}[htb!]
\centering
\begin{tabular}{cccccc}
       & FOM & \hspace{-0.4cm} $\sigma_L = 1$ & \hspace{-0.4cm} $\sigma_L = 2$ & \hspace{-0.4cm} $\sigma_L = 5$ & \hspace{-0.4cm} $\sigma_L=10$ \\
    $\tildep_1$ & \includegraphics[align=c,scale = 0.3]{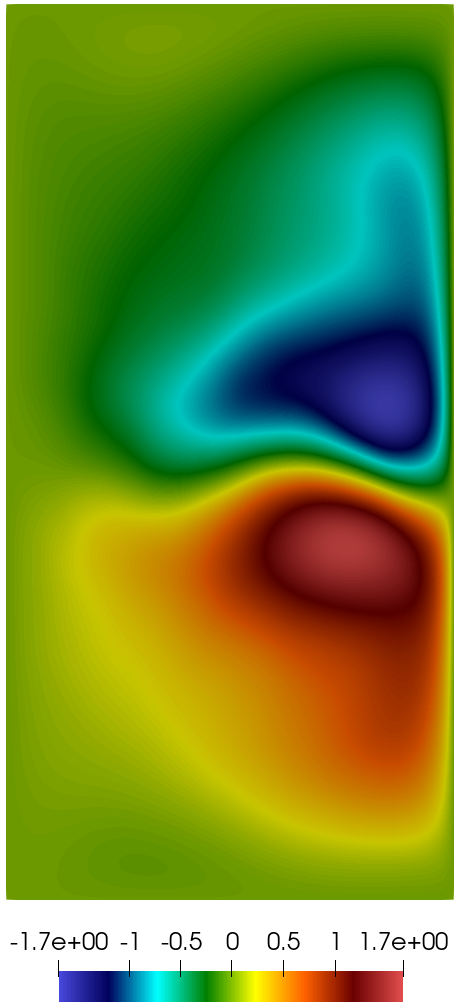} & \hspace{-0.4cm}\includegraphics[align=c,scale = 0.3]{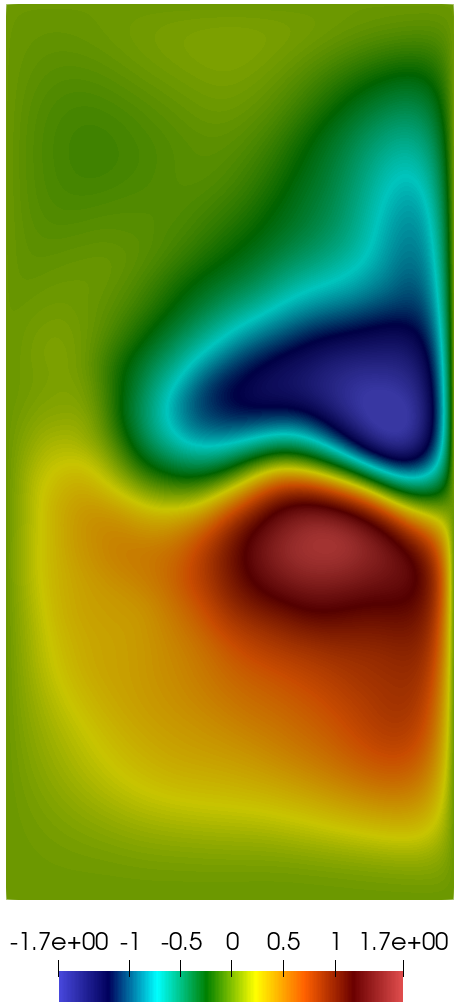} & \hspace{-0.4cm}\includegraphics[align=c,scale = 0.3]{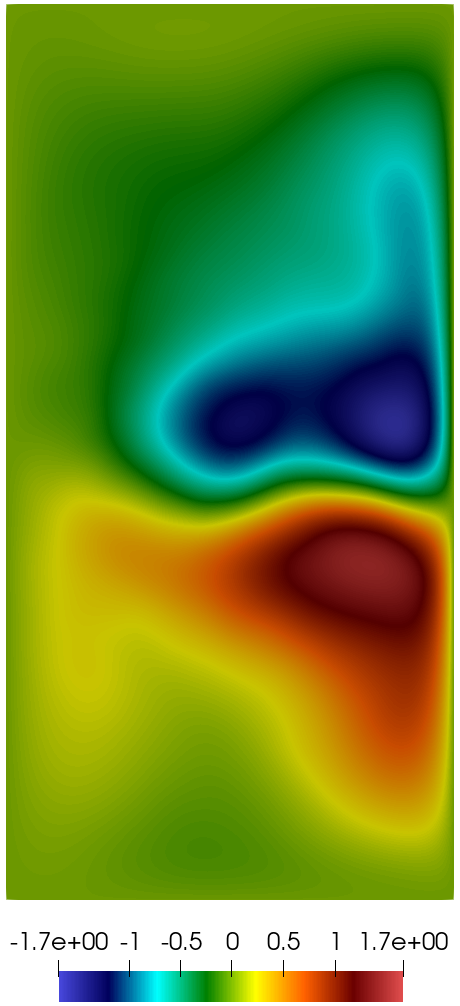} & \hspace{-0.4cm}\includegraphics[align=c,scale = 0.3]{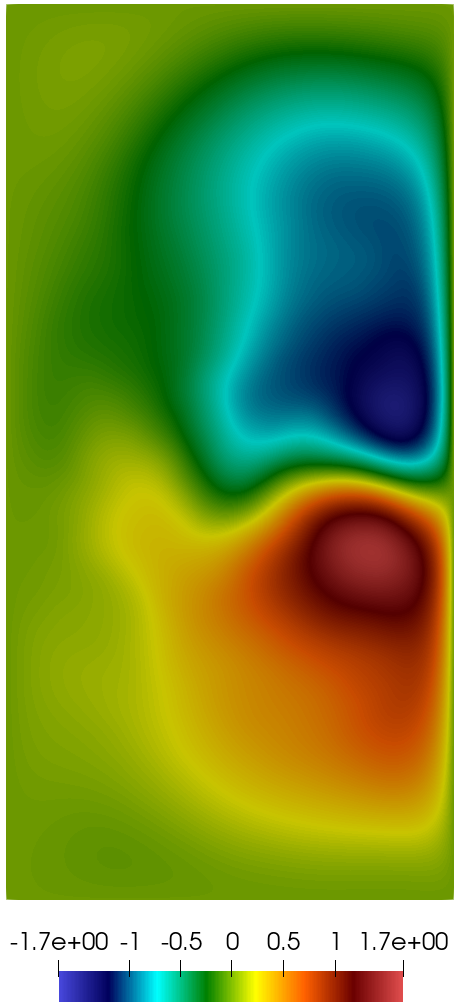} & \hspace{-0.4cm}\includegraphics[align=c,scale = 0.3]{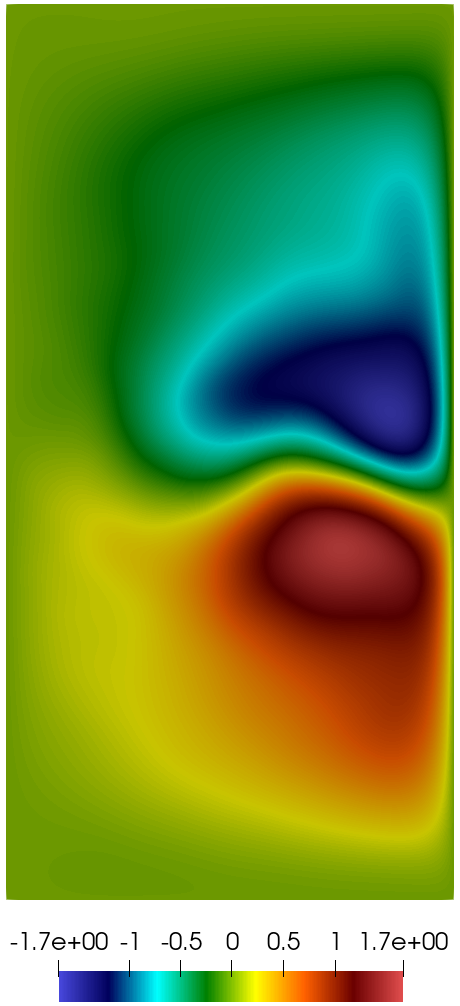} \\
    $\tildep_2$ & \includegraphics[align=c,scale = 0.3]{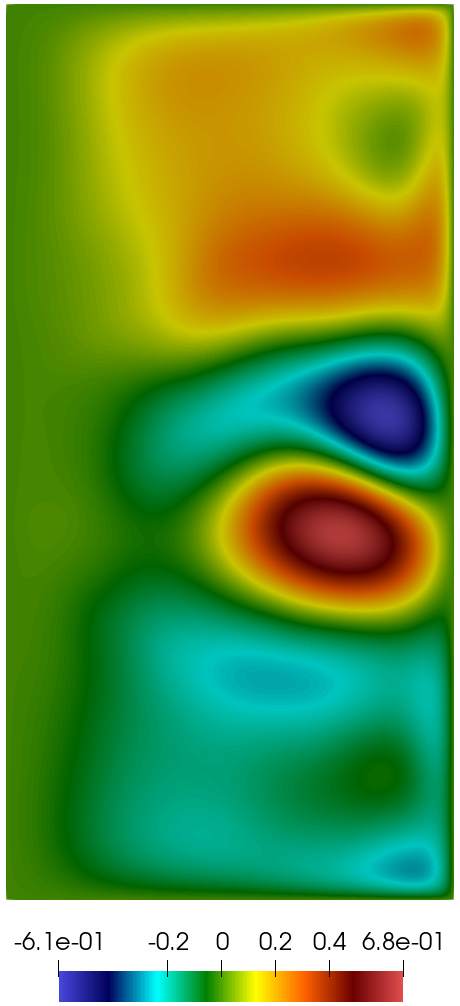} & \hspace{-0.4cm}\includegraphics[align=c,scale = 0.3]{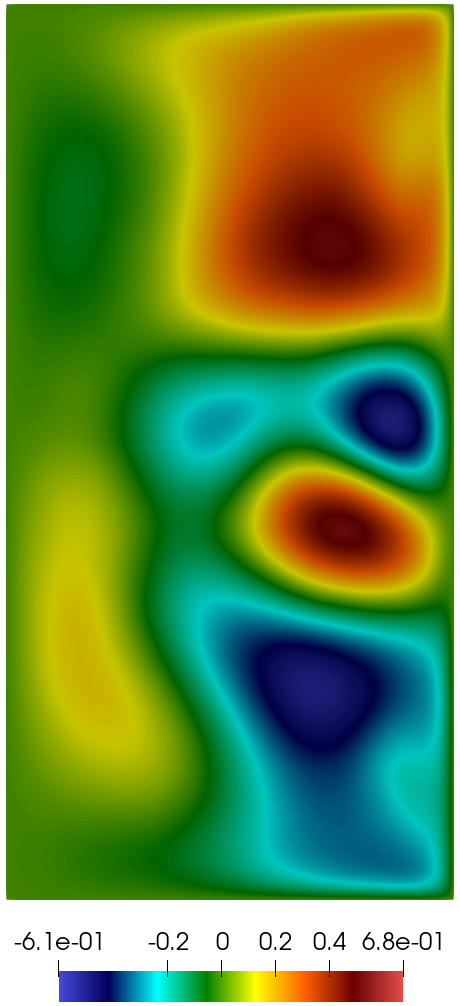} & \hspace{-0.4cm}\includegraphics[align=c,scale = 0.3]{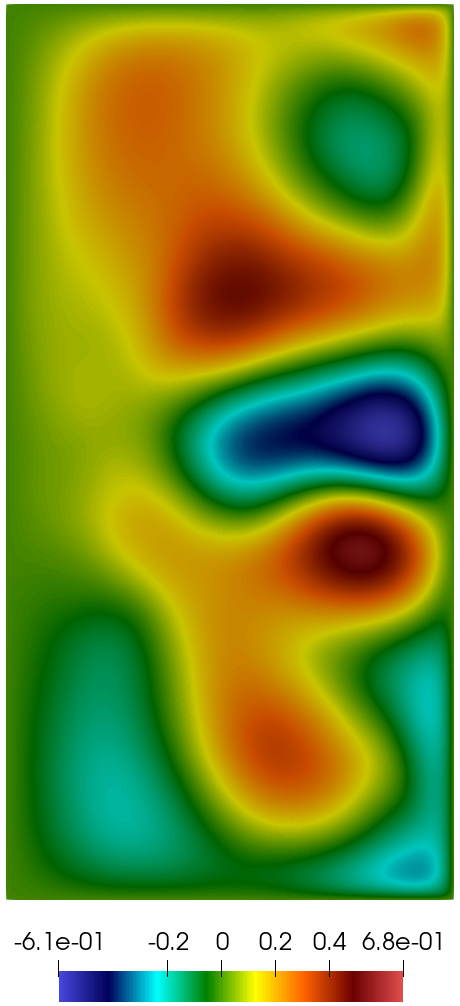} & \hspace{-0.4cm}\includegraphics[align=c,scale = 0.3]{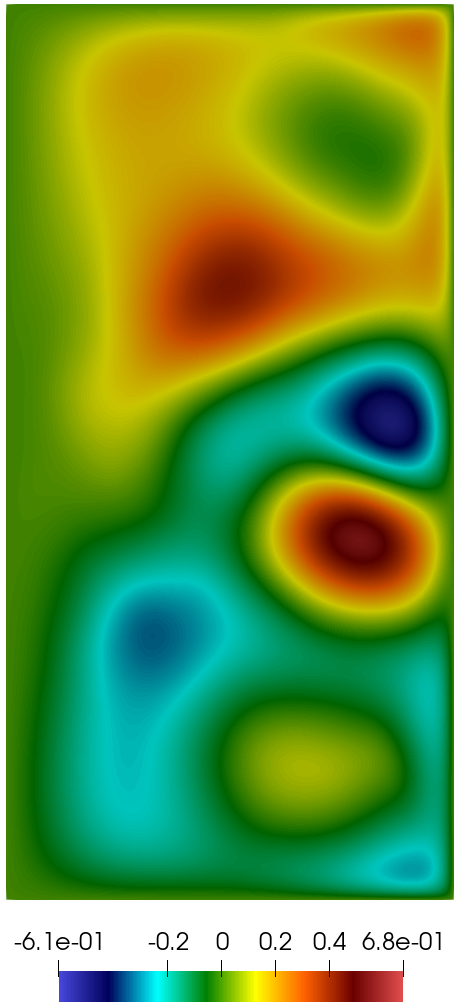} & \hspace{-0.4cm}\includegraphics[align=c,scale = 0.3]{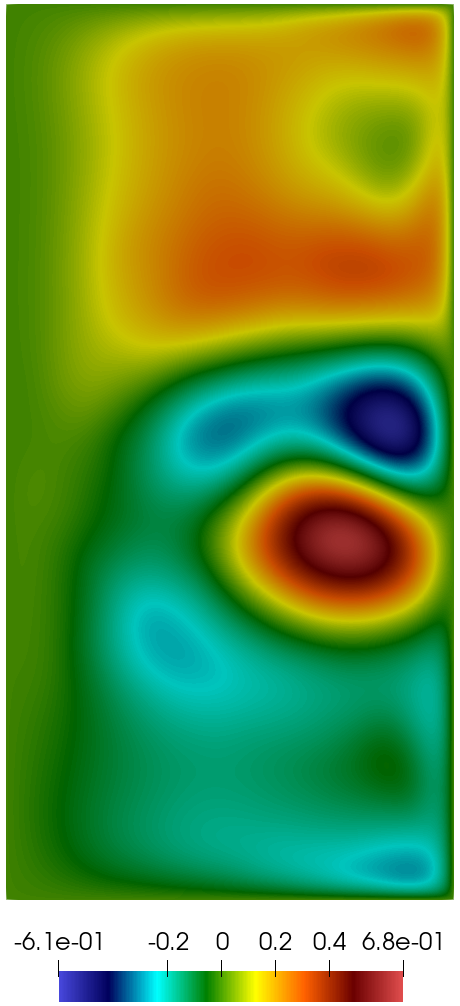} \\
\end{tabular}
\caption{$\tildep_1$ and $\tildep_2$ computed over the predictive time interval $[50,100]$ by the FOM (first column) and the POD-LSTM  $\calmf_{\psi_1}$ (first row) and $\calmf_{\psi_2}$ (second row) for different values of the lookback window $\sigma_L$ (second to fifth columns).}
\label{fig:psi_mu}
\end{figure}

Let us shift our focus to the predictive performance of the POD-LSTM ROM for the stream functions $\psi_1$ and $\psi_2$. We recall that the hyperparameters 
of the LSTM network are reported in Table \ref{tab:hyperparam}, third column. Following 
the analysis for $\calmf_{q_1}$ and $\calmf_{q_2}$, we first fix the number of retained POD basis functions to $N_{\psi_1}^r = N_{\psi_2}^r = 10$, which means retaining $54\%$ and $64\%$ of the cumulative singular value energy for $\psi_1'$ and $\psi_2'$, respectively.
Fig.~\ref{fig:psi_mu} compares the true $\tildep_1$ and $\tildep_2$ (first column) with $\tildep_1$ and $\tildep_2$ computed by the POD-LSTM ROM with different values of $\sigma_L$ (second to fifth columns). While the POD-LSTM ROM  manages to capture the correct number of gyres in $\tildep_1$
for all values of $\sigma_L$ and gets a good approximation of the gyres shape for $\sigma_L=1,10$, it cannot reliably predict $\tildep_2$ for lower values of $\sigma_L$.
However, for $\sigma_L=10$ the POD-LSTM ROM 
produces the correct shape of the gyres in $\tildep_2$, especially the outer gyres 
that can be harder to predict \cite{Besabe2024}.

To quantify the comparison in Fig.~\ref{fig:psi_mu},
Table \ref{tab:error_psi_mu} reports the errors $\varepsilon_{\psi_l}^{(1)}$ \eqref{eq:rmse} and $\varepsilon_{\psi_l}^{(2)}$ \eqref{eq:l2-error} for $l = 1,2$. In both metrics, $\sigma_L=10$ produces the lowest errors.

\begin{table}[h!]
    \centering
    \begin{tabular}{|c|c|c|c|c|}
        \hline
        $\sigma_L$ & $\varepsilon_{\psi_1}^{(1)}$ & $\varepsilon_{\psi_1}^{(2)}$ & $\varepsilon_{\psi_2}^{(1)}$ & $\varepsilon_{\psi_2}^{(2)}$ \\
        \hline
        1 & 1.102E-01 & 1.911E-01 & 1.206E-01 & 6.579E-01 \\
        \hline
        2 & 1.671E-01 & 2.900E-01 & 1.630E-01 & 8.893E-01 \\
        \hline
        5 & 1.654E-01 & 2.870E-01 & 9.909E-02 & 5.406E-01 \\
        \hline
        10 & 6.041E-02 & 1.048E-01 & 4.543E-02 & 2.479E-01 \\
        \hline
    \end{tabular}
    \caption{Error metrics $\varepsilon_{\psi_l}^{(1)}$ \eqref{eq:rmse} and $\varepsilon_{\psi_l}^{(2)}$ \eqref{eq:l2-error}, $l = 1,2 $, for different
    values of the lookback window $\sigma_L$.}
    \label{tab:error_psi_mu}
\end{table}

\begin{figure}[htb!]
    \centering
    \begin{subfigure}[h]{0.45\textwidth}
         \centering
         \includegraphics[width=\textwidth]{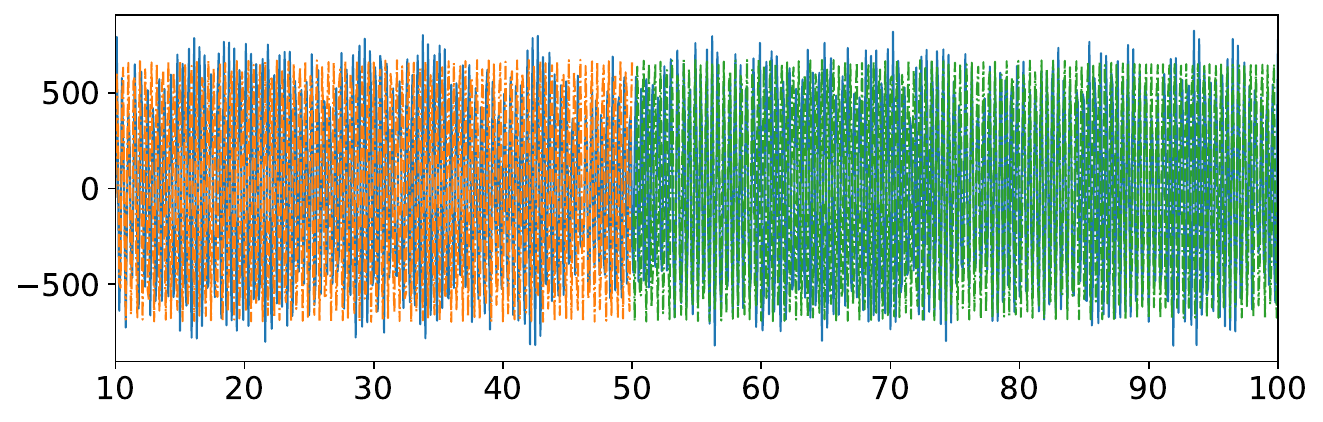}
         \caption{\scriptsize{$\beta_{1,1}(t),\sigma_L=1$}}
     \end{subfigure}
     \begin{subfigure}[h]{0.45\textwidth}
         \centering
         \includegraphics[width=\textwidth]{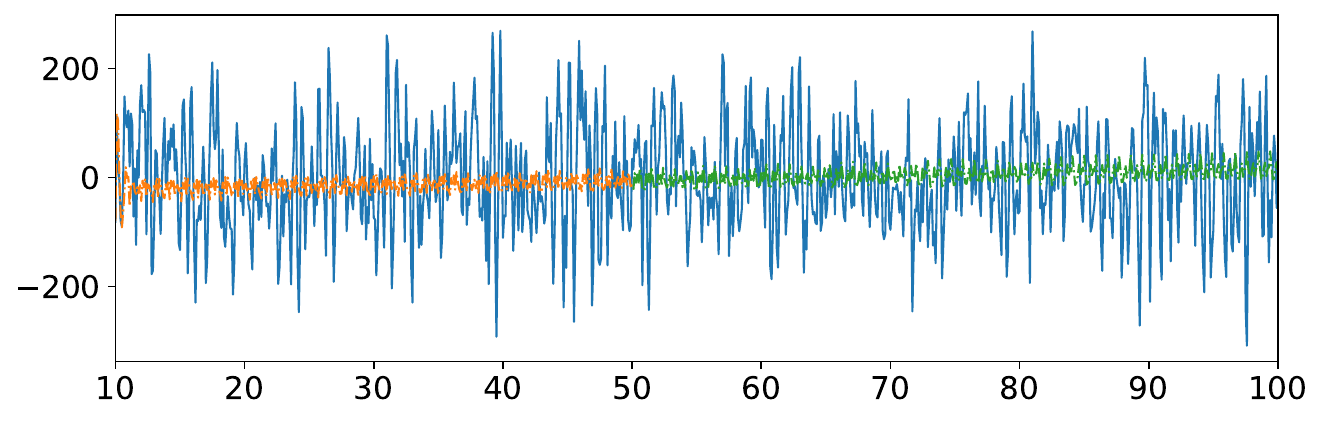}
         \caption{\scriptsize{$\beta_{1,3}(t),\sigma_L=1$}}
     \end{subfigure}

     \begin{subfigure}[h]{0.45\textwidth}
         \centering
         \includegraphics[width=\textwidth]{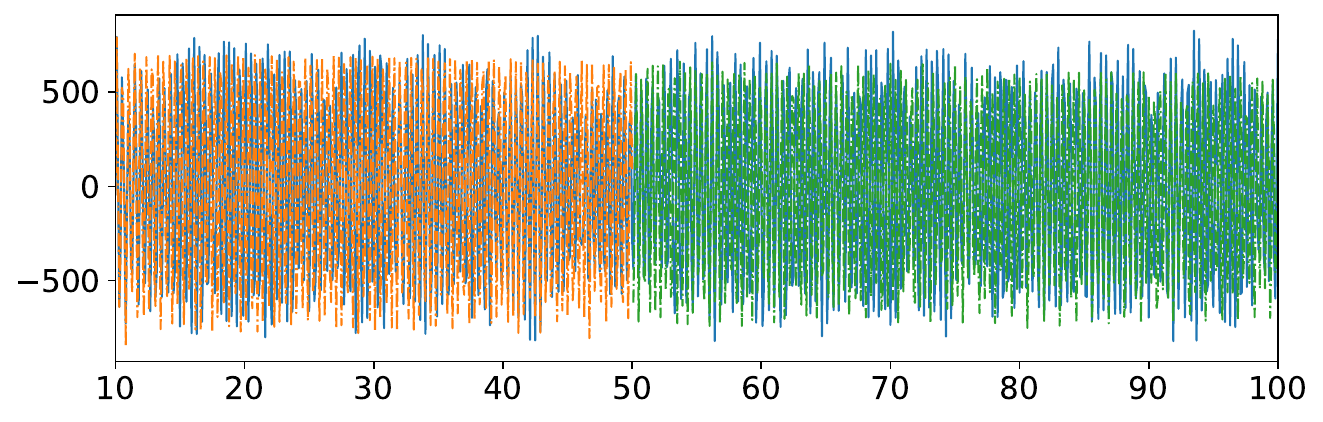}
         \caption{\scriptsize{$\beta_{1,1}(t),\sigma_L=2$}}
     \end{subfigure}
     \begin{subfigure}[h]{0.45\textwidth}
         \centering
         \includegraphics[width=\textwidth]{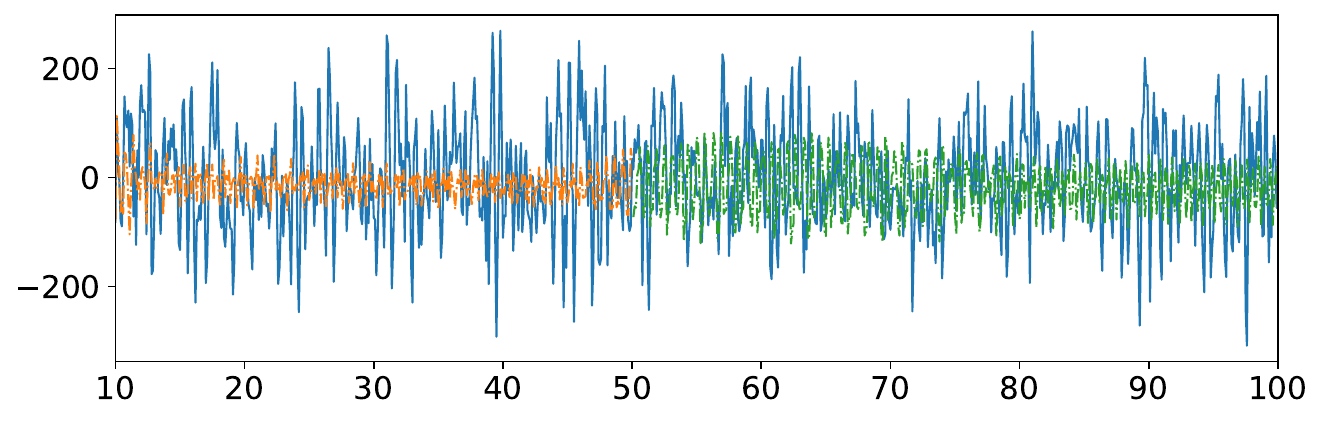}
         \caption{\scriptsize{$\beta_{1,3}(t),\sigma_L=2$}}
     \end{subfigure}

     \begin{subfigure}[h]{0.45\textwidth}
         \centering
         \includegraphics[width=\textwidth]{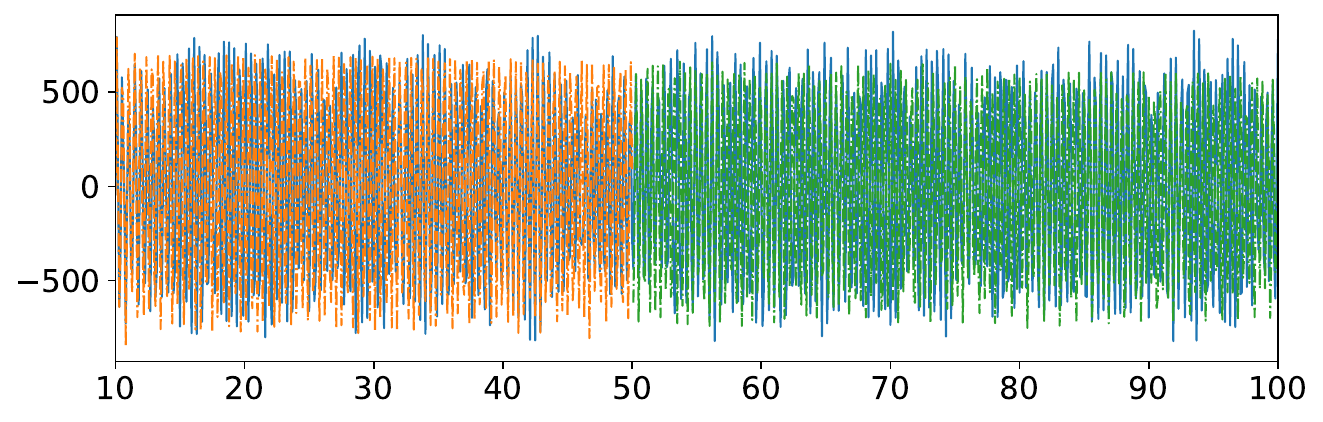}
         \caption{\scriptsize{$\beta_{1,1}(t),\sigma_L=5$}}
     \end{subfigure}
     \begin{subfigure}[h]{0.45\textwidth}
         \centering
         \includegraphics[width=\textwidth]{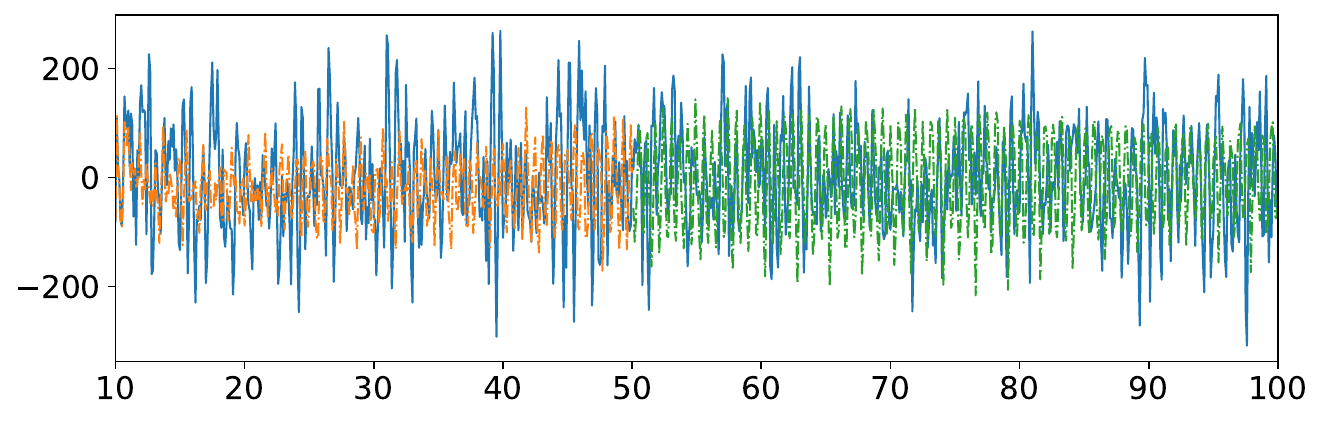}
         \caption{\scriptsize{$\beta_{1,3}(t),\sigma_L=5$}}
     \end{subfigure}

     \begin{subfigure}[h]{0.45\textwidth}
         \centering
         \includegraphics[width=\textwidth]{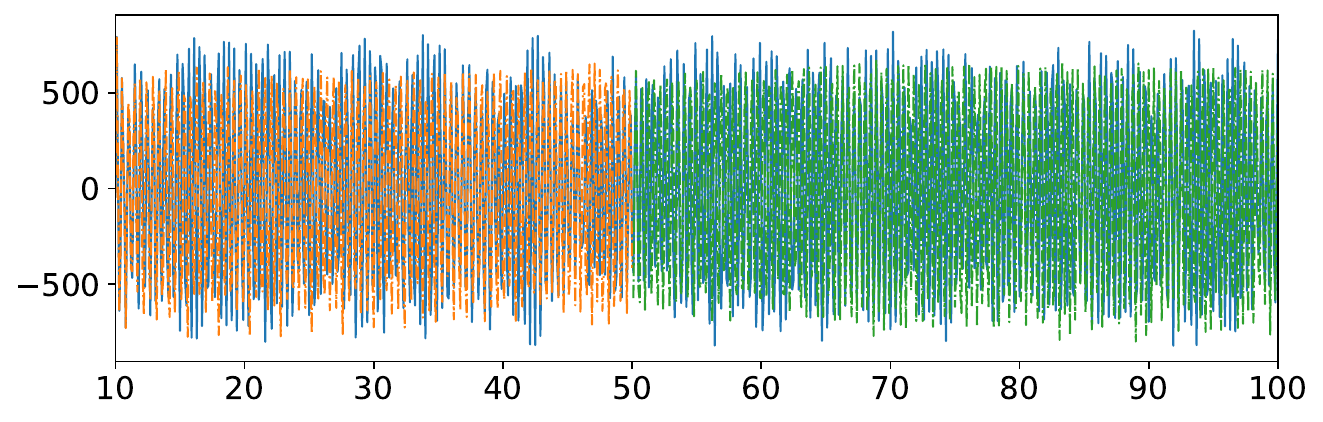}
         \caption{\scriptsize{$\beta_{1,1}(t),\sigma_L=10$}}
     \end{subfigure}
     \begin{subfigure}[h]{0.45\textwidth}
         \centering
         \includegraphics[width=\textwidth]{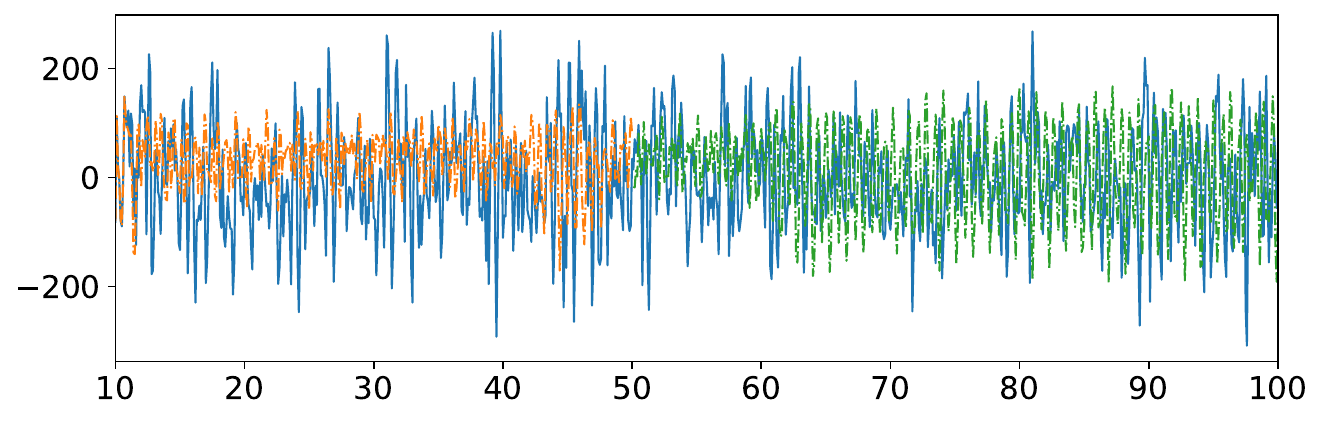}
         \caption{\scriptsize{$\beta_{1,3}(t),\sigma_L=10$}}
     \end{subfigure}

     \begin{subfigure}[h]{0.45\textwidth}
         \centering
         \vspace{0.3cm}
         \includegraphics[width=\textwidth]{figs/coeff-legend.png}
     \end{subfigure}
\caption{Time evolution of the first $\beta_{1,1}$ (left) and third $\beta_{1,3}$ (right) modal coefficients for $\psi_1$ over the time interval $[10,100]$ for different values of the lookback window $\sigma_L$. 
}
\label{fig:psi1-coeff-mu}
\end{figure}

\begin{figure}[htb!]
    \centering
    \begin{subfigure}[h]{0.45\textwidth}
         \centering
         \includegraphics[width=\textwidth]{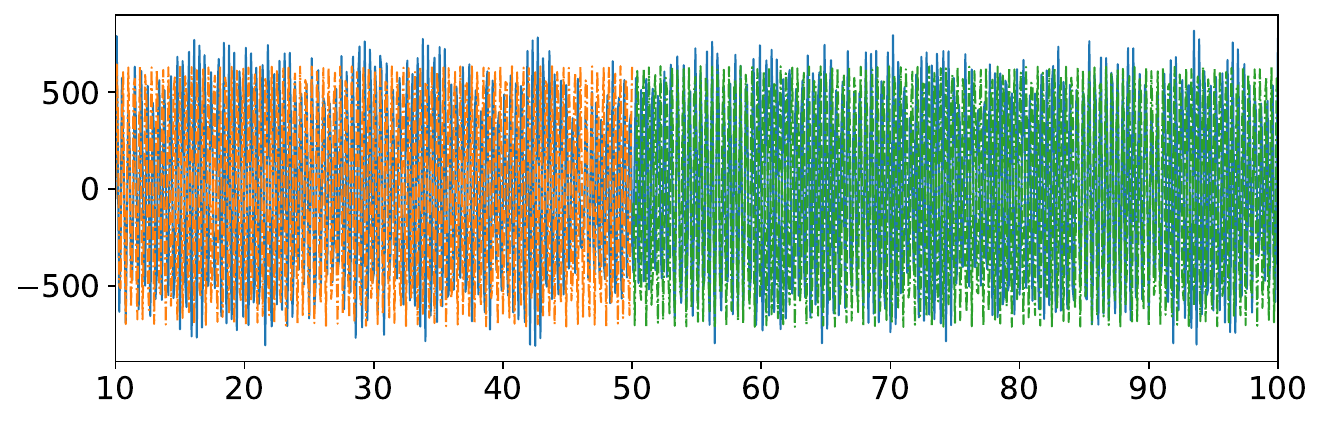}
         \caption{\scriptsize{$\beta_{2,1}(t),\sigma_L=1$}}
     \end{subfigure}
     \begin{subfigure}[h]{0.45\textwidth}
         \centering
         \includegraphics[width=\textwidth]{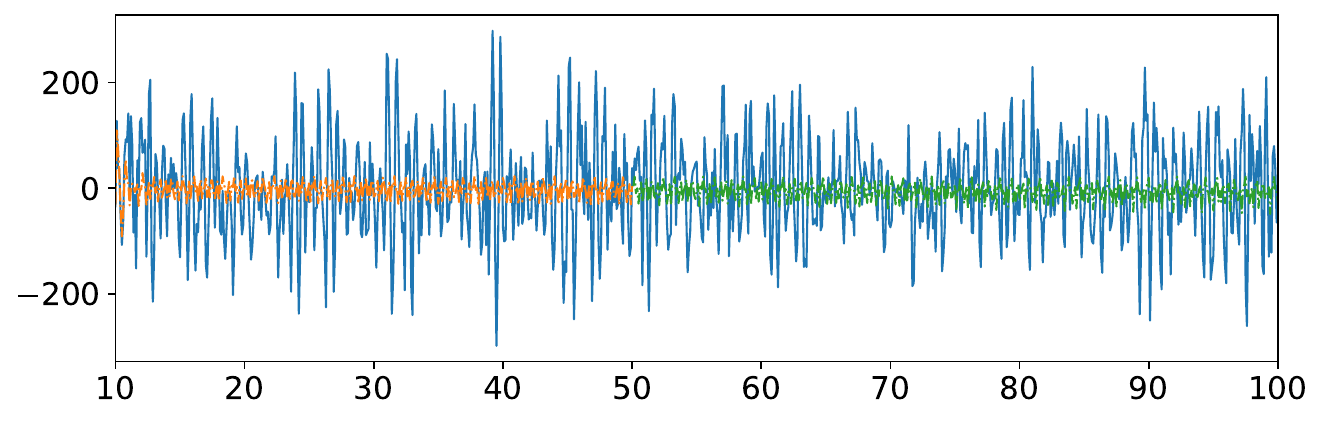}
         \caption{\scriptsize{$\beta_{2,3}(t),\sigma_L=1$}}
     \end{subfigure}

     \begin{subfigure}[h]{0.45\textwidth}
         \centering
         \includegraphics[width=\textwidth]{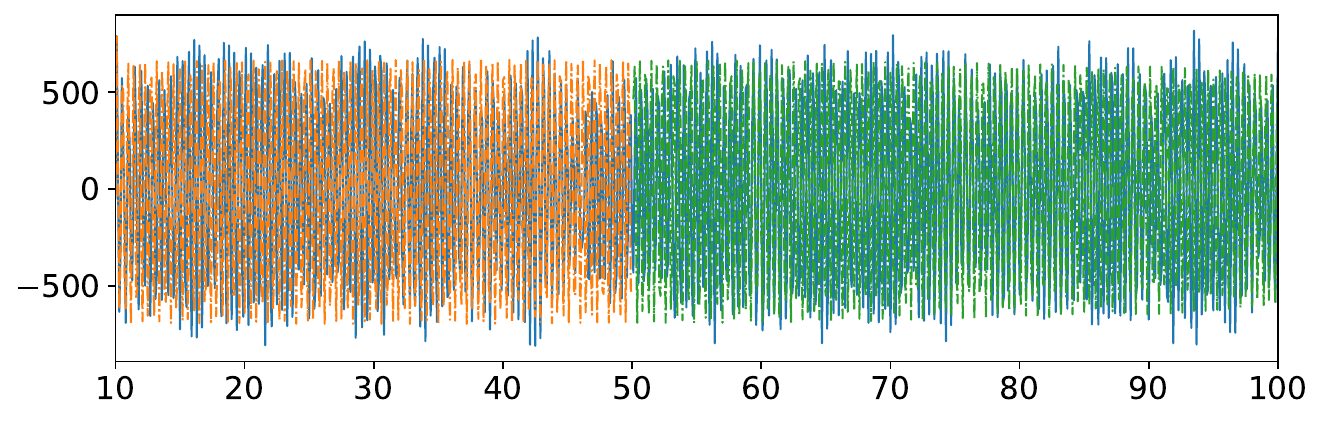}
         \caption{\scriptsize{$\beta_{2,1}(t),\sigma_L=2$}}
     \end{subfigure}
     \begin{subfigure}[h]{0.45\textwidth}
         \centering
         \includegraphics[width=\textwidth]{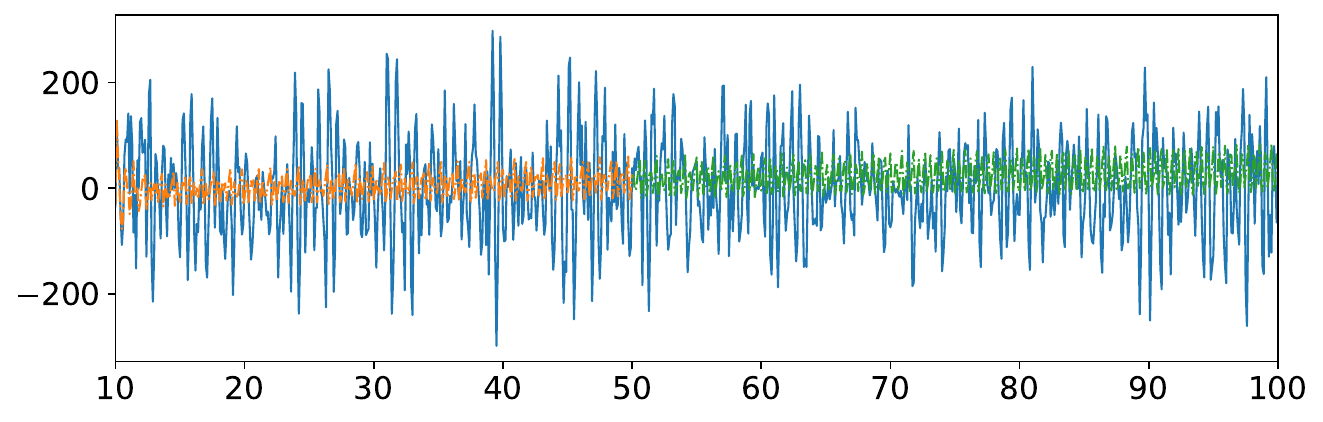}
         \caption{\scriptsize{$\beta_{2,3}(t),\sigma_L=2$}}
     \end{subfigure}

     \begin{subfigure}[h]{0.45\textwidth}
         \centering
         \includegraphics[width=\textwidth]{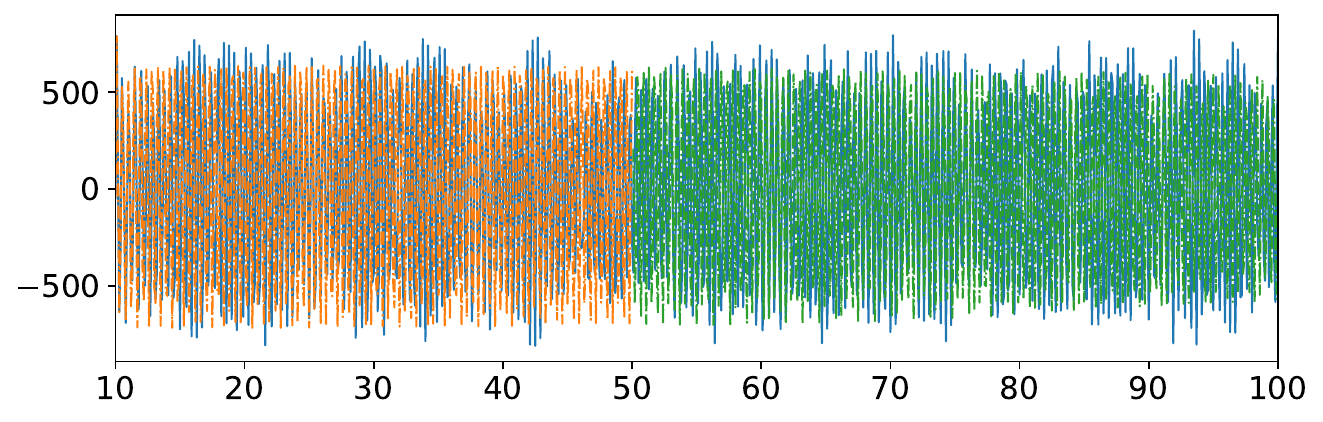}
         \caption{\scriptsize{$\beta_{2,1}(t),\sigma_L=5$}}
     \end{subfigure}
     \begin{subfigure}[h]{0.45\textwidth}
         \centering
         \includegraphics[width=\textwidth]{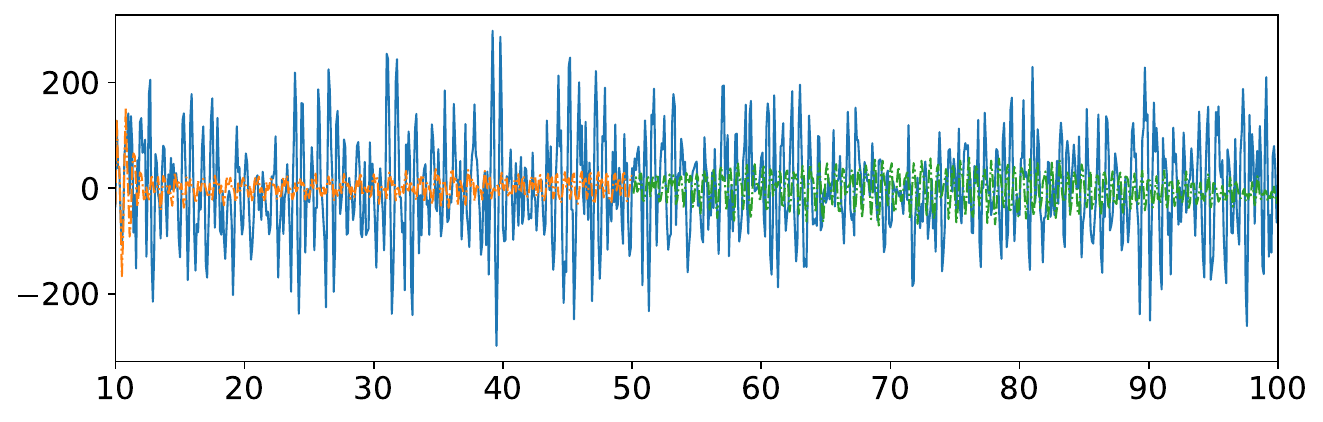}
         \caption{\scriptsize{$\beta_{2,3}(t),\sigma_L=5$}}
     \end{subfigure}

     \begin{subfigure}[h]{0.45\textwidth}
         \centering
         \includegraphics[width=\textwidth]{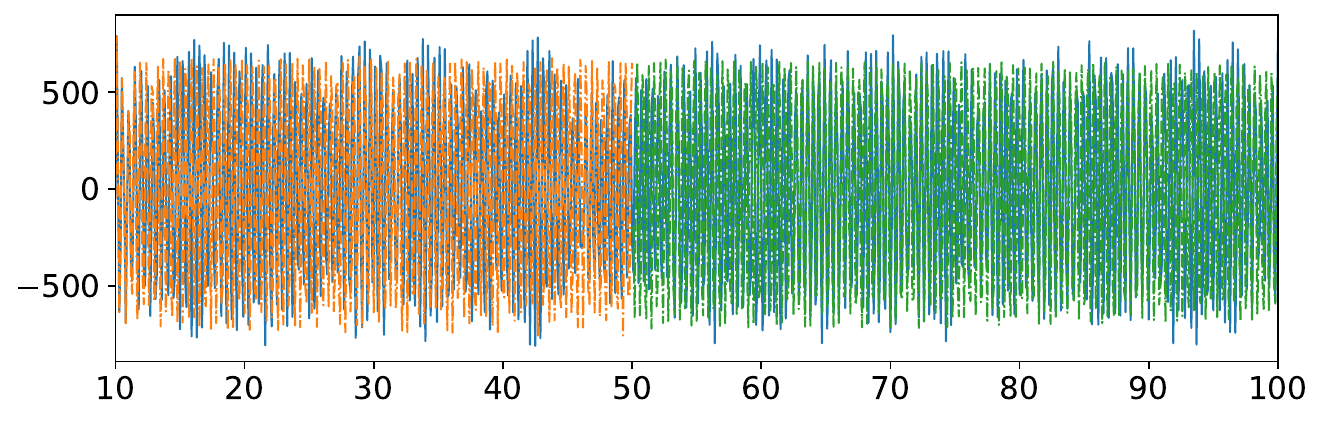}
         \caption{\scriptsize{$\beta_{2,1}(t),\sigma_L=10$}}
     \end{subfigure}
     \begin{subfigure}[h]{0.45\textwidth}
         \centering
         \includegraphics[width=\textwidth]{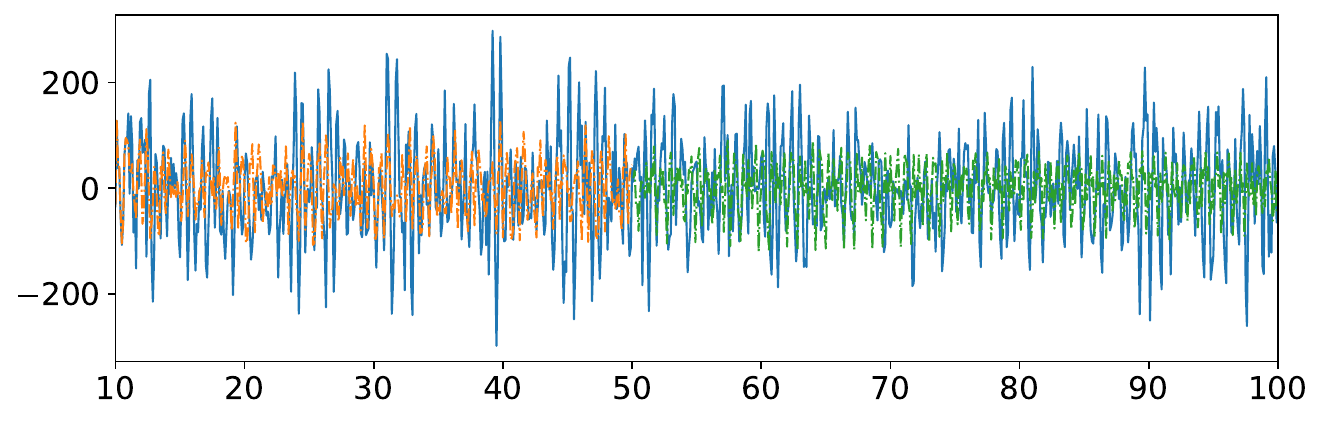}
         \caption{\scriptsize{$\beta_{2,3}(t),\sigma_L=10$}}
     \end{subfigure}

     \begin{subfigure}[h]{0.45\textwidth}
         \centering
         \vspace{0.3cm}
         \includegraphics[width=\textwidth]{figs/coeff-legend.png}
     \end{subfigure}
\caption{Time evolution of the first $\beta_{2,1}$ (left) and third $\beta_{2,3}$ (right) modal coefficients for $\psi_2$ over the time interval $[10,100]$ for different values of the lookback window $\sigma_L$. 
}
\label{fig:psi2-coeff-mu}
\end{figure}

\begin{figure}[htb!]
    \centering
    \begin{subfigure}[h]{0.4\textwidth}
         \centering
         \includegraphics[width=\textwidth]{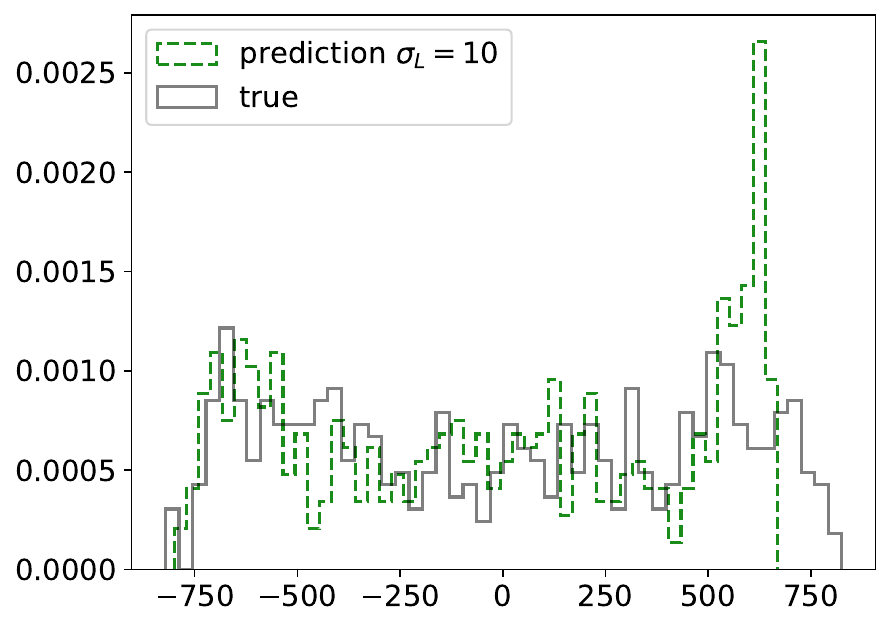}
         \caption{\scriptsize{$\beta_{1,1}(t)$}}
     \end{subfigure}
     \begin{subfigure}[h]{0.4\textwidth}
         \centering
         \includegraphics[width=\textwidth]{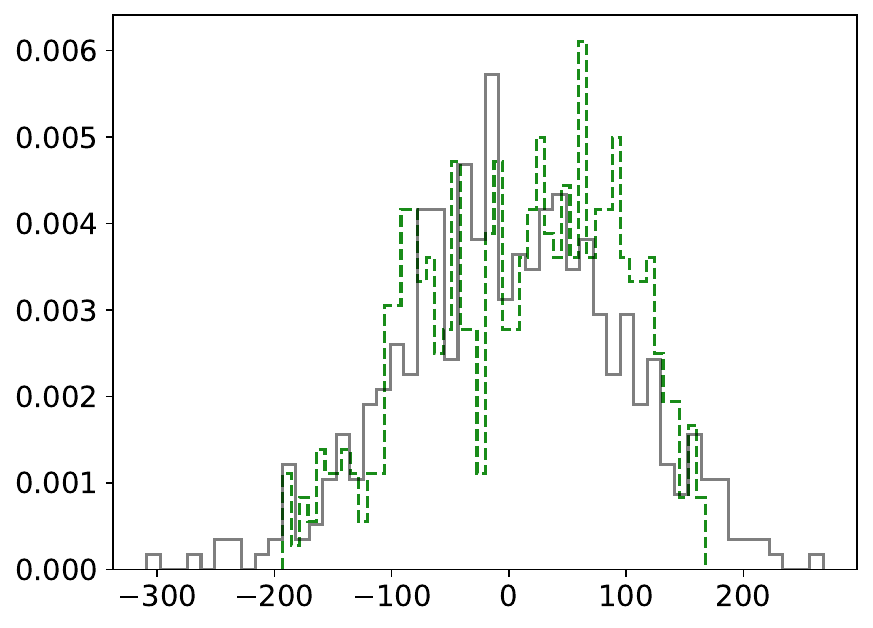}
         \caption{\scriptsize{$\beta_{1,3}(t)$}}
     \end{subfigure}

     \begin{subfigure}[h]{0.4\textwidth}
         \centering
         \includegraphics[width=\textwidth]{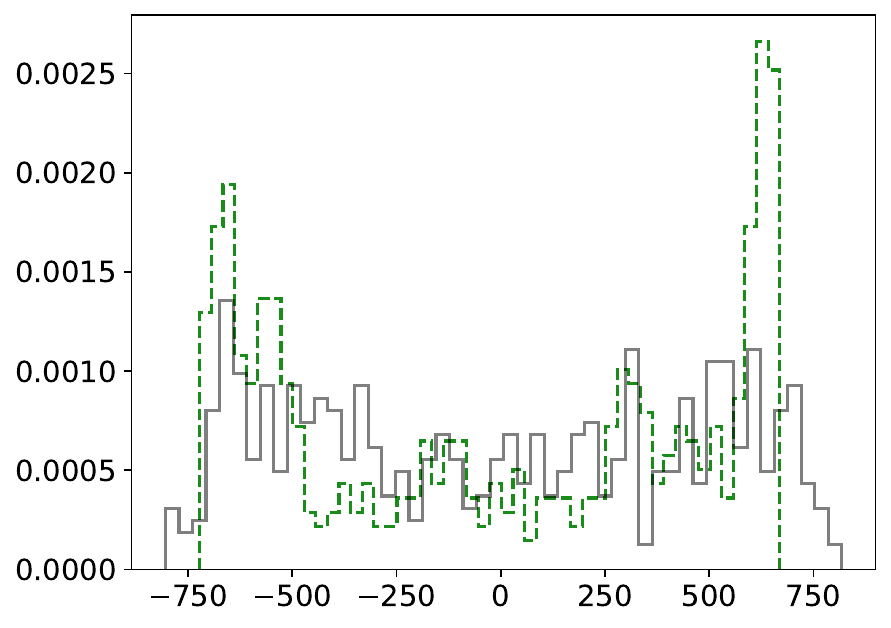}
         \caption{\scriptsize{$\beta_{2,1}(t)$}}
     \end{subfigure}
     \begin{subfigure}[h]{0.4\textwidth}
         \centering
         \includegraphics[width=\textwidth]{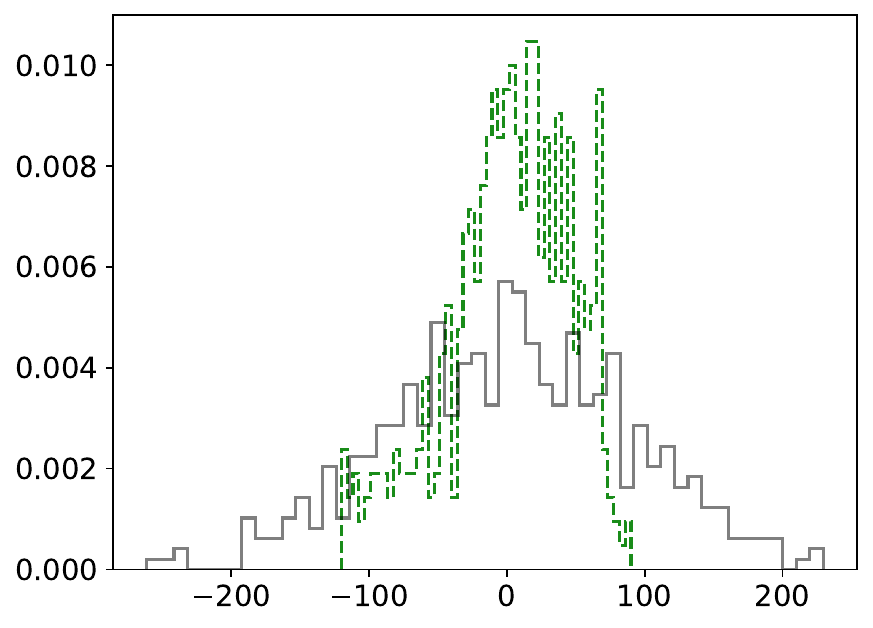}
         \caption{\scriptsize{$\beta_{2,3}(t)$}}
     \end{subfigure}
\caption{
Probability mass functions of the first and third modal coefficients of $\psi_1$ (top row) and $\psi_2$ (bottom row) over time interval $[50,100]$: 
comparison between true value, i.e.,  
for the FOM solution projected onto the appropriate reduced space, and the value computed by $\calmf_{\psi_1}$ and $\calmf_{\psi_2}$ for $\sigma_L=10$.
}
\label{fig:psi-coeff_hist}
\end{figure}

\begin{figure}[htb!]
    \centering
    \begin{subfigure}[h]{0.47\textwidth}
         \centering
         \includegraphics[width=\textwidth]{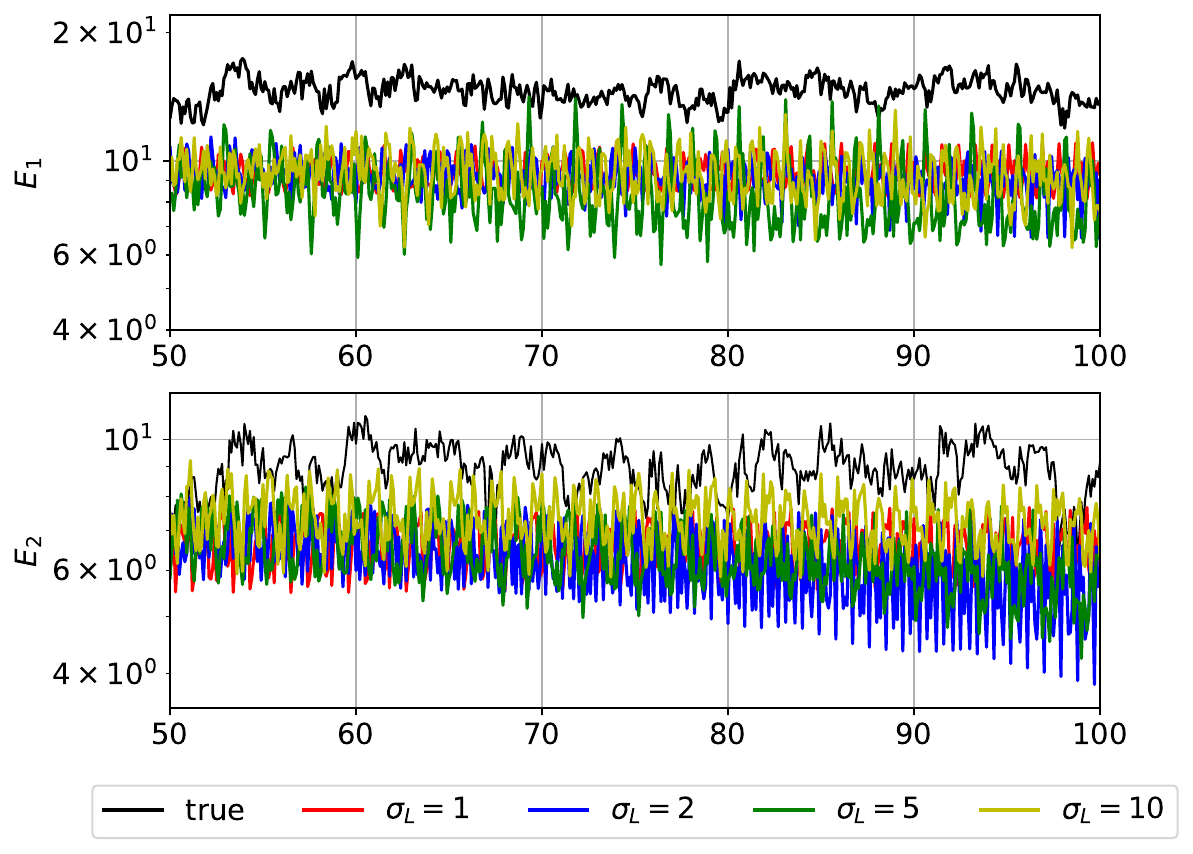}
     \end{subfigure}
     \begin{subfigure}[h]{0.47\textwidth}
         \centering
         \includegraphics[width=\textwidth]{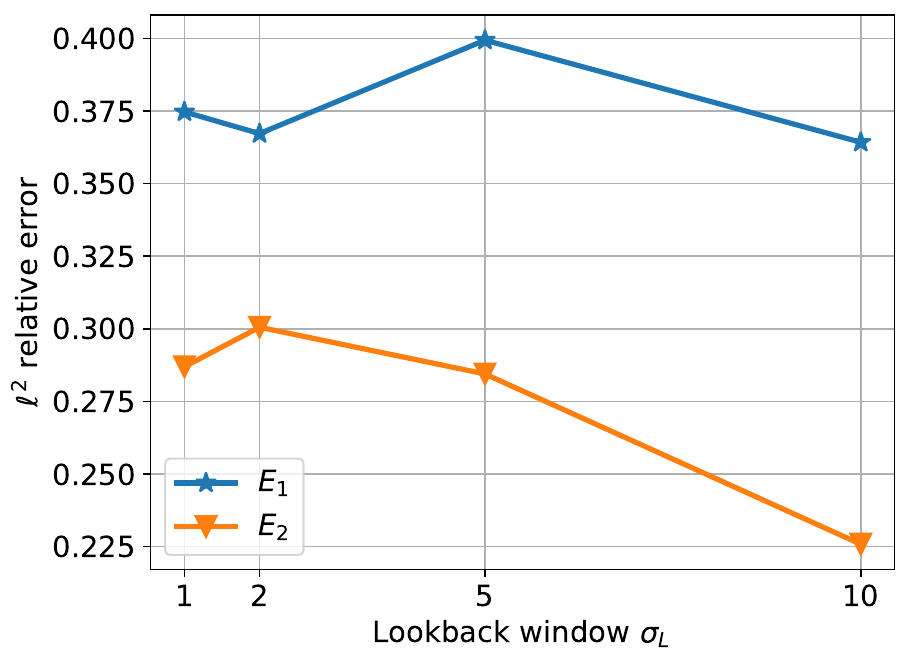}
     \end{subfigure}
\caption{
Time evolution of the kinetic energy for the top layer $\mathcal{E}_1$ (top left) and the bottom layer $\mathcal{E}_2$ (bottom left) computed using the FOM and the POD-LSTM ROM with different lookback window $\sigma_L$. 
The corresponding relative $L^2$ error is shown on the right.}
\label{fig:ke-mu}
\end{figure}

In order to understand the time-dependent predictive performance of POD-LSTM ROM, we show the first $\beta_{l,1}$ and third $\beta_{l,3}$ modal coefficients computed using $\calmf_{\psi_l}$, $l = 1,2$, with different $\sigma_L$
in Fig.~\ref{fig:psi1-coeff-mu}-\ref{fig:psi2-coeff-mu}. 
While $\beta_{1,1}$ and $\beta_{2,1}$ are rather well approximated regardless of the value of 
$\sigma_L$, for $\beta_{1,3}$ and $\beta_{2,3}$ we see that $\sigma_L=10$ gives the most accurate reconstructions and predictions. 
Fig.~\ref{fig:psi-coeff_hist}
compares the probability mass functions 
of the first and third modal coefficients 
over the unseen time interval $[50,100]$
for the  FOM solution projected onto 
the appropriate reduced space with the counterparts computed by $\calmf_{\psi_1}$ and $\calmf_{\psi_2}$ for 
$\sigma_L=10$. The true and predicted PMFs are in good agreement, except for $\beta_{2,3}$.
The comparison in terms of kinetic energy evolution \eqref{eq:kin-energy} is shown in Fig.~\ref{fig:ke-mu} (left) for different values of $\sigma_L$.
We see that the POD-LSTM ROM underestimates the kinetic energy for either layer, which we suspect is caused by the retention of only roughly half of
the cumulative singular value energy. 
The lowest relative $L^2$ errors with respect to the true kinetic energy of the system is given by 
$\sigma_L=10$. See the right panel of Fig.~\ref{fig:ke-mu} right.

\begin{figure}[htb!]
\centering
\begin{tabular}{cccccc}
       & FOM & \hspace{-0.4cm}$N_{\psi_l}^r = 2$ & \hspace{-0.4cm}$N_{\psi_l}^r = 4$ & \hspace{-0.4cm}$N_{\psi_l}^r = 8$ & \hspace{-0.4cm}$N_{\psi_l}^r = 10$ \\
    $\tildep_1$ & \includegraphics[align=c,scale = 0.3]{figs/psi1/psi1Mean_true.png} & \hspace{-0.4cm}\includegraphics[align=c,scale = 0.3]{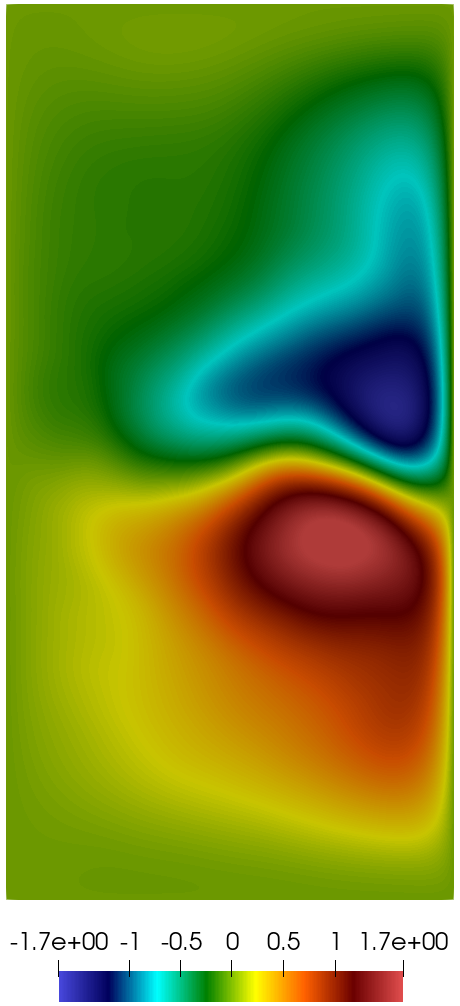} & \hspace{-0.4cm}\includegraphics[align=c,scale = 0.3]{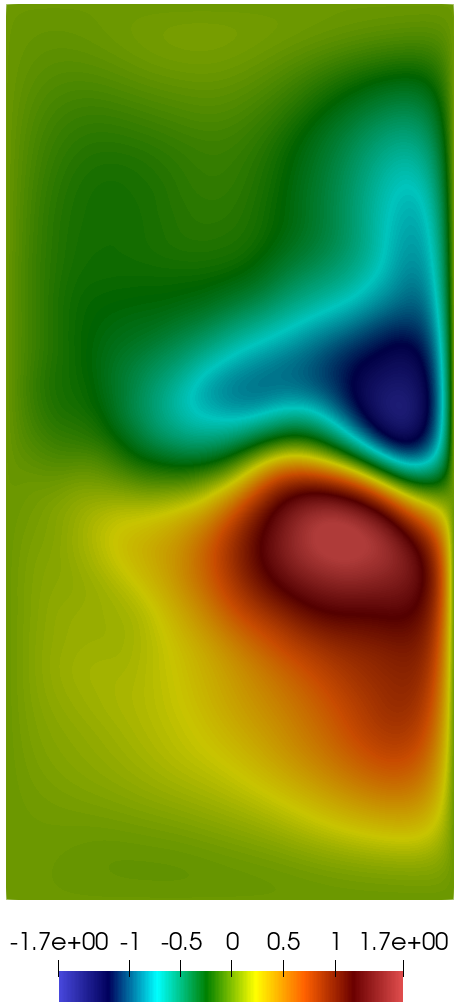} & \hspace{-0.4cm}\includegraphics[align=c,scale = 0.3]{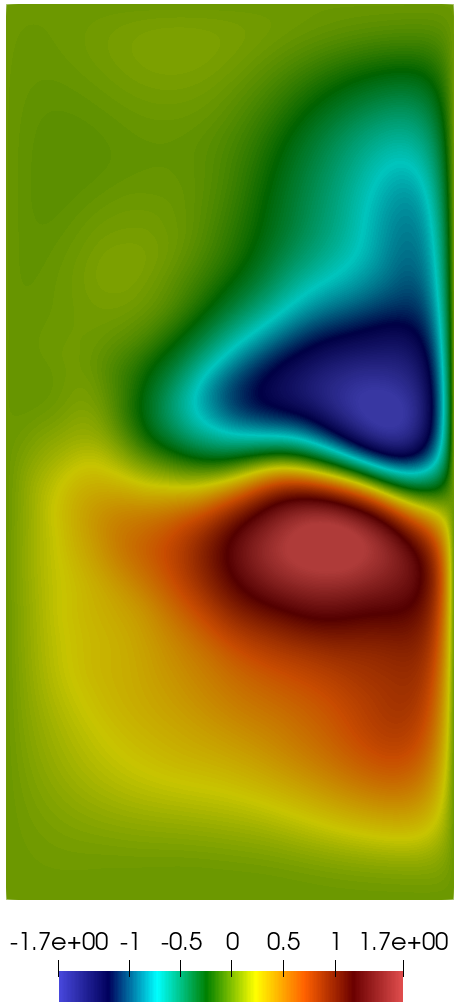} & \hspace{-0.4cm}\includegraphics[align=c,scale = 0.3]{figs/psi1/psi1Mean_mu10.png} \\
    $\tildep_2$ & \includegraphics[align=c,scale = 0.3]{figs/psi2/psi2Mean_true.png} & \hspace{-0.4cm}\includegraphics[align=c,scale = 0.3]{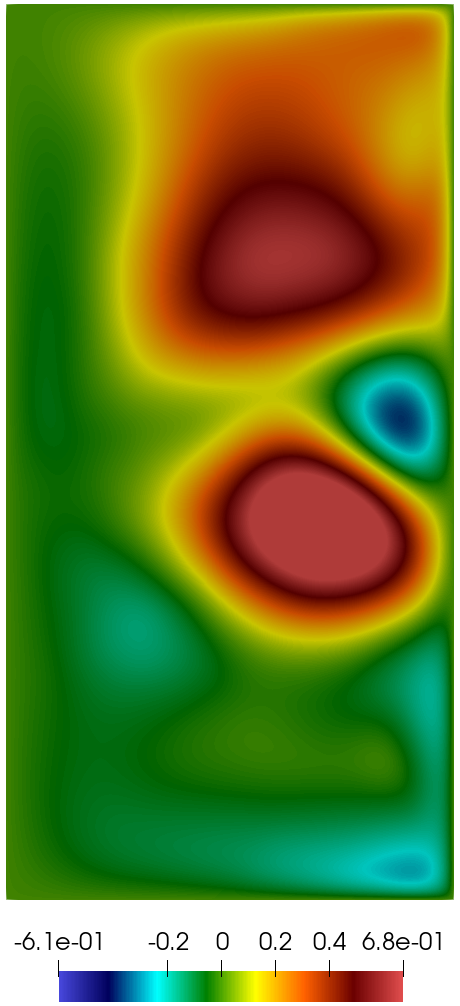} & \hspace{-0.4cm}\includegraphics[align=c,scale = 0.3]{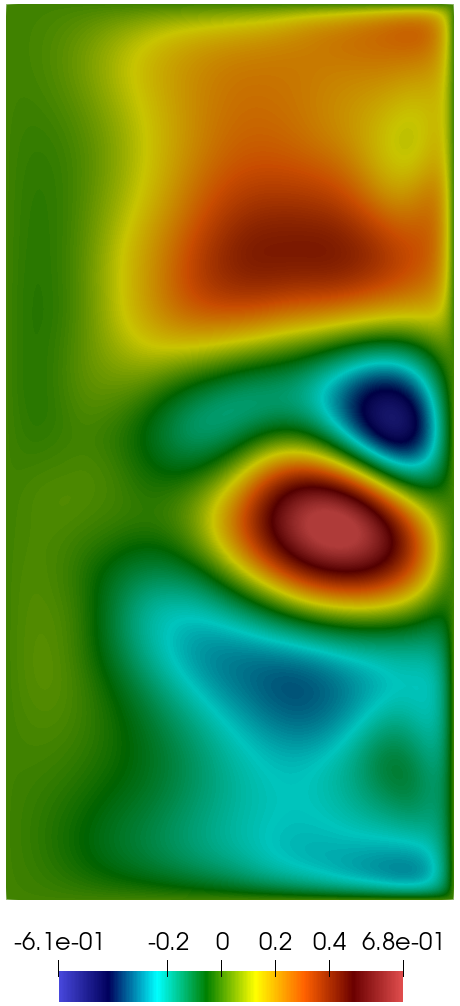} & \hspace{-0.4cm}\includegraphics[align=c,scale = 0.3]{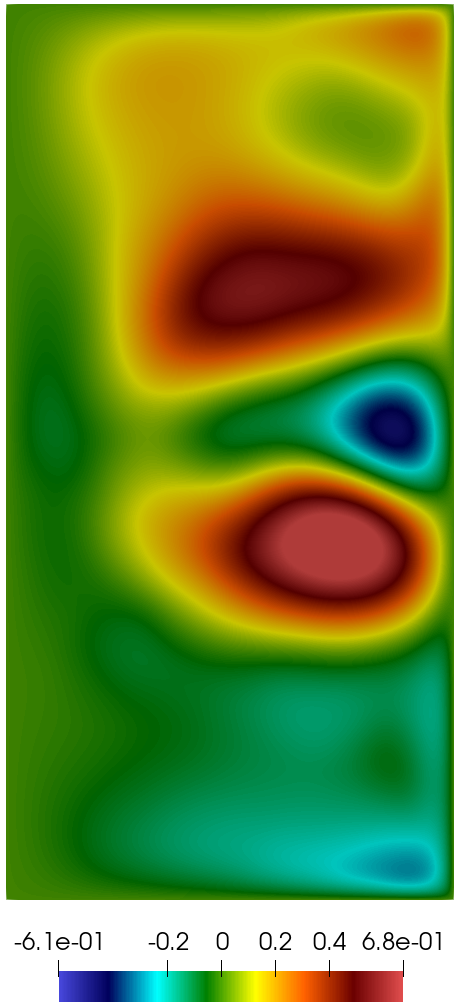} & \hspace{-0.4cm}\includegraphics[align=c,scale = 0.3]{figs/psi2/psi2Mean_mu10.png} \\
\end{tabular}
\caption{
Time-averaged stream functions $\widetilde{\psi}_1$ (top) and $\widetilde{\psi}_2$ (bottom) over the predictive time interval $[50,100]$ computed by FOM (first column) and POD-LSTM ROM with $N_{\psi_l}^r=2,4,8,10$ (second to fifth columns).
}
\label{fig:psi_mod}
\end{figure}

Next, we fix the lookback window to $\sigma_L=10$ for both $\calmf_{\psi_1}$ and $\calmf_{\psi_2}$
and vary $N_{\psi_l}$.
Fig.~\ref{fig:psi_mod} shows the time-averaged potential vorticities $\widetilde{\psi}_1$ and $\widetilde{\psi}_2$ computed using the FOM and the POD-LSTM ROM with different numbers of modal basis functions retained. 
We observe that even when retaining just $N_{\psi_1}^r=1$, the POD-LSTM ROM is able to capture the gyre shapes in $\widetilde{\psi}_1$.
In the case of $\widetilde{\psi}_2$ instead, 
a very low numbers of modal basis functions retained
leads to large regions of overshoot (bright red in the northern gyre) and incorrect gyre shapes. 
As $N_{\psi_l}^r$ increases, $\widetilde{\psi}_2$ computed by the ROM becomes more accurate, 
with $N_{\psi_l}^r = 10$ providing the most accurate solution.

Error metrics $\varepsilon_{\psi_l}^{(1)}$ \eqref{eq:rmse} and $\varepsilon_{\psi_l}^{(2)}$ \eqref{eq:l2-error} for the stream functions, together with the relative $L^2$ error for the kinetic energy, are presented in Table \ref{tab:error-psi1-mod}. The lowest errors are for  $N_{\psi_l}^r = 10$, with the exception of the relative $L^2$ error for $E_1$, which is lowest for $N_{\psi_1}^r=2$.

\begin{table}[h!]
    \centering
    \begin{tabular}{|c|cccc|}
    \hline
        $N_{\psi_1}^r$ & $\delta_{\psi_1}$ & $\varepsilon_{\psi_1}^{(1)}$ & $\varepsilon_{\psi_1}^{(2)}$ & $E_1$ rel. err.\\
        \hline
        2 & 0.34 & 7.936E-02 & 1.377E-01 & 3.166E-01 \\
        \hline
        4 & 0.41 & 1.082E-01 & 1.879E-01 & 4.228E-01 \\
        \hline
        8 & 0.50 & 1.249E-01 & 2.169E-01 & 4.130E-01 \\
        \hline
        10 & 0.54 & 6.041E-02 & 1.048E-01 & 3.642E-01 \\
        \hline
        $N_{\psi_2}^r$ & $\delta_{\psi_2}$ & $\varepsilon_{\psi_2}^{(1)}$ & $\varepsilon_{\psi_2}^{(2)}$ & $E_2$ rel. err.\\
        \hline
        2 & 0.42 & 1.708E-01 & 9.319E-01 & 6.699E-01 \\
        \hline
        4 & 0.49 & 6.127E-02 & 3.343E-01 & 8.616E-01 \\
        \hline
        8 & 0.60 & 1.223E-01 & 6.672E-01 & 9.000E-01 \\
        \hline
        10 & 0.64 & 4.543E-02 & 2.479E-01 & 2.258E-01 \\
        \hline
    \end{tabular}
    \caption{
    Fraction of retained eigenvalue energy $\delta_{\psi_l}$, 
    error metrics $\varepsilon_{\psi_l}^{(1)}$ \eqref{eq:rmse} and $\varepsilon_{\psi_l}^{(2)}$ \eqref{eq:l2-error}, 
    and relative $L^2$ error for the kinetic energy  ${E}_l$ \eqref{eq:enstrophy},
    $l = 1,2 $, for different numbers of POD modes retained.}
    \label{tab:error-psi1-mod}
\end{table}

\begin{figure}[htb!]
    \centering
    \begin{subfigure}[h]{0.45\textwidth}
         \centering
         \includegraphics[width=\textwidth]{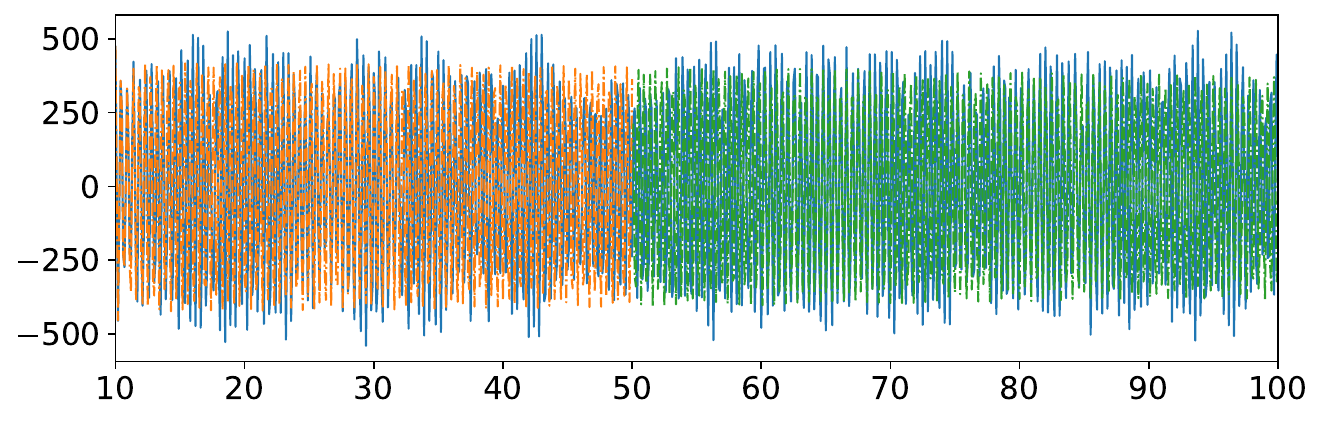}
         \caption{\scriptsize{$\beta_{1,2}(t),N_{\psi_1}^r=2$}}
     \end{subfigure}
     \begin{subfigure}[h]{0.45\textwidth}
         \centering
         \includegraphics[width=\textwidth]{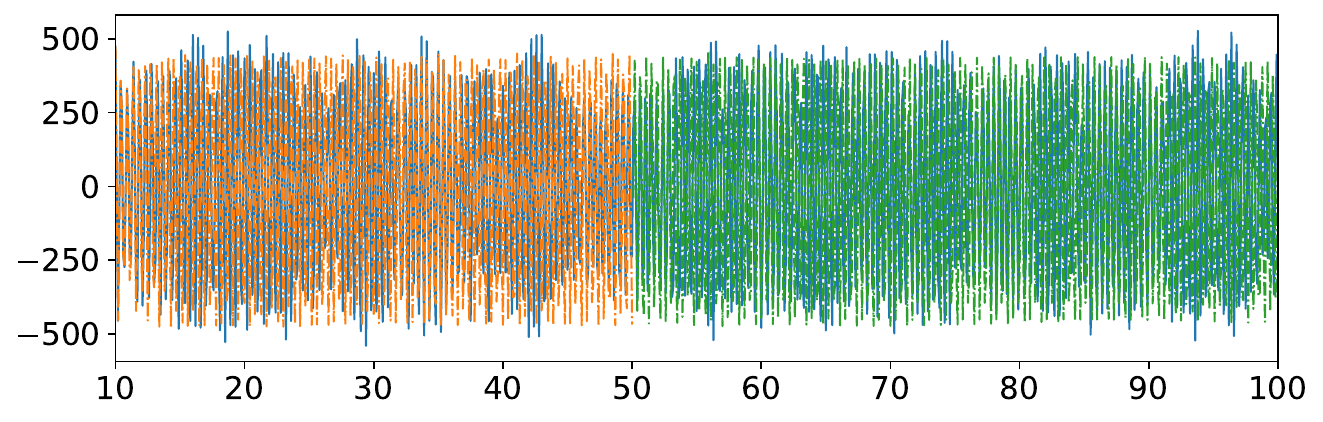}
         \caption{\scriptsize{$\beta_{1,2}(t),N_{\psi_1}^r=4$}}
     \end{subfigure}
    \begin{subfigure}[h]{0.45\textwidth}
         \centering
         \includegraphics[width=\textwidth]{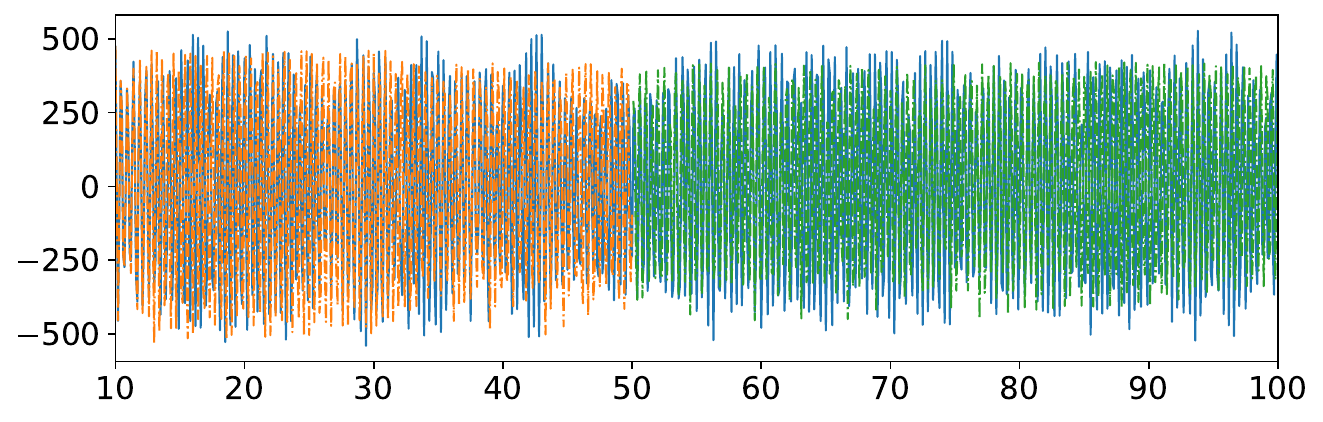}
         \caption{\scriptsize{$\beta_{1,2}(t),N_{\psi_1}^r=8$}}
     \end{subfigure}
     \begin{subfigure}[h]{0.45\textwidth}
         \centering
         \includegraphics[width=\textwidth]{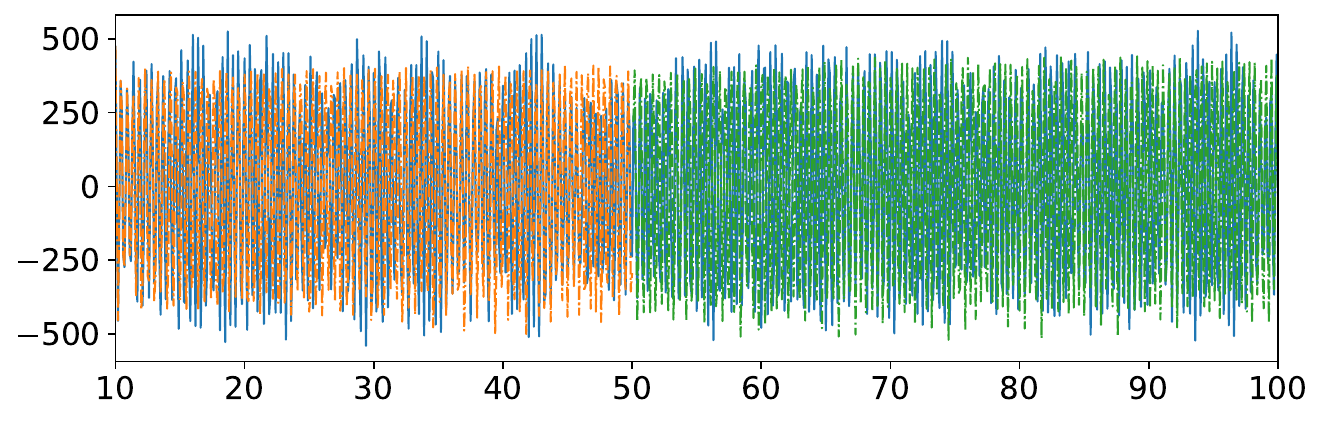}
         \caption{\scriptsize{$\beta_{1,2}(t),N_{\psi_1}^r=10$}}
     \end{subfigure}

    \begin{subfigure}[h]{0.45\textwidth}
         \centering
         \includegraphics[width=\textwidth]{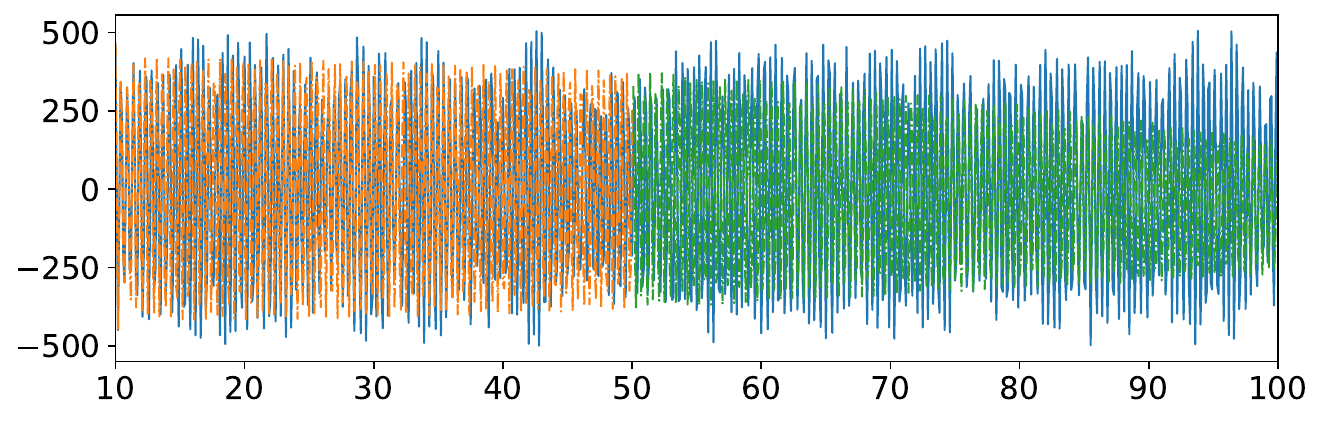}
         \caption{\scriptsize{$\beta_{2,2}(t),N_{\psi_2}^r=2$}}
    \end{subfigure}
    \begin{subfigure}[h]{0.45\textwidth}
         \centering
         \includegraphics[width=\textwidth]{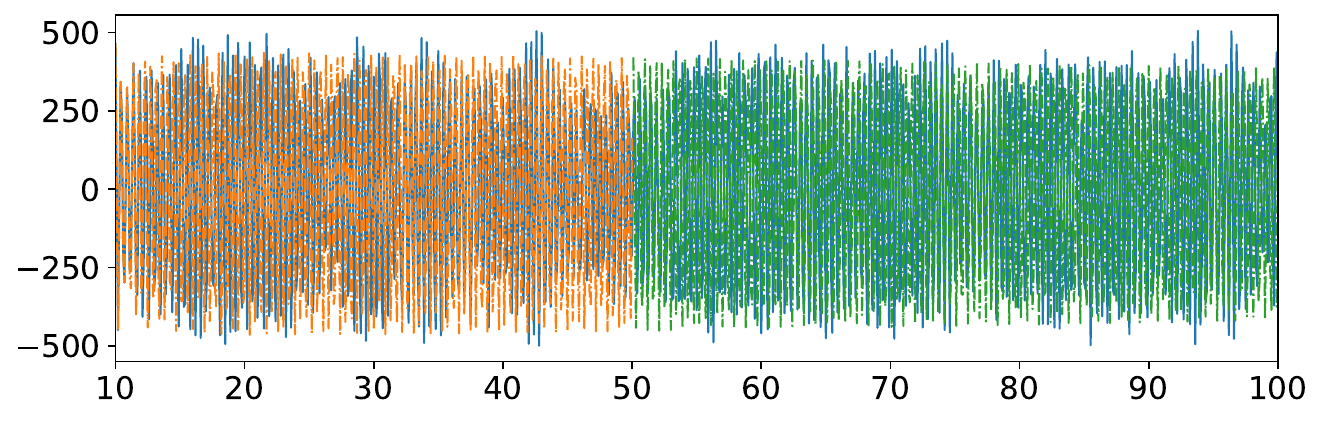}
         \caption{\scriptsize{$\beta_{2,2}(t),N_{\psi_2}^r=4$}}
    \end{subfigure}
    \begin{subfigure}[h]{0.45\textwidth}
         \centering
         \includegraphics[width=\textwidth]{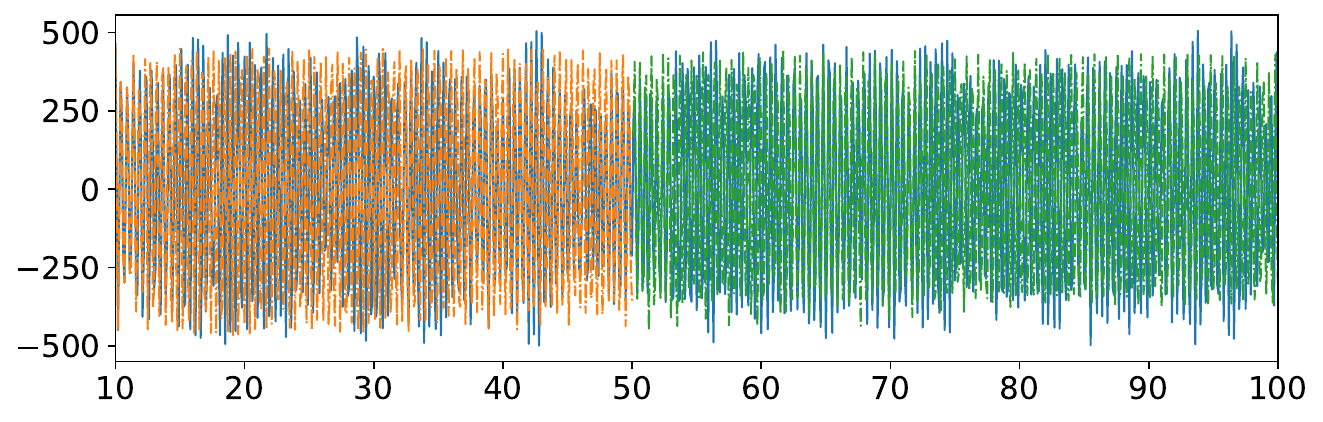}
         \caption{\scriptsize{$\beta_{2,2}(t),N_{\psi_2}^r=8$}}
    \end{subfigure}
    \begin{subfigure}[h]{0.45\textwidth}
         \centering
         \includegraphics[width=\textwidth]{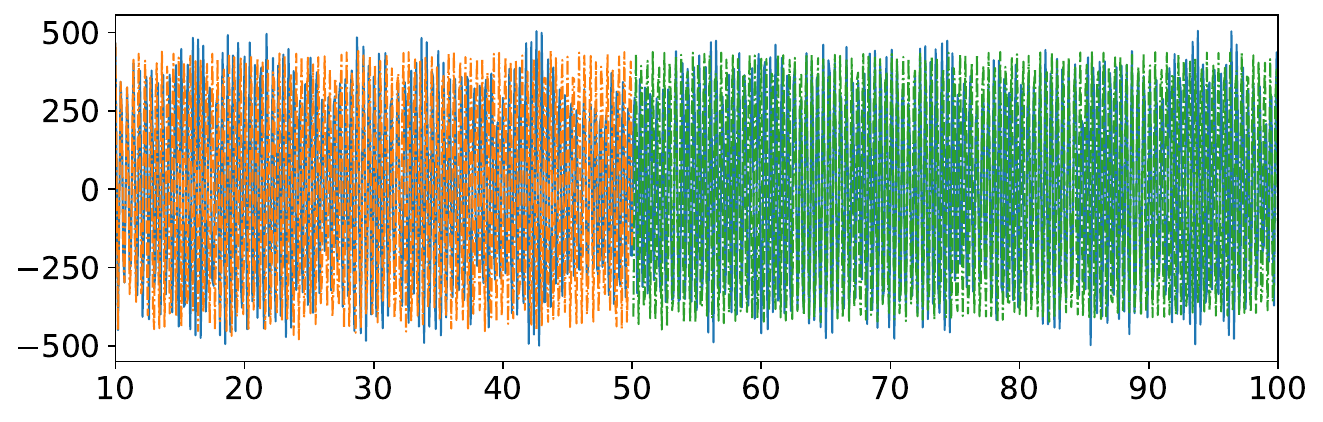}
         \caption{\scriptsize{$\beta_{2,2}(t),N_{\psi_2}^r=10$}}
    \end{subfigure}

    \begin{subfigure}[h]{0.45\textwidth}
         \centering
         \vspace{0.3cm}
         \includegraphics[width=\textwidth]{figs/coeff-legend.png}
     \end{subfigure}
\caption{
Time evolution of the second modal coefficient $\beta_{l,2}$ over the time interval $[10,100]$ for different 
numbers of retained POD basis function $N_{\psi_l}^r$, $l = 1,2 $.
}
\label{fig:psi-coeff-mod}
\end{figure}

The evolution in time of 
second $\beta_{l,2}$ modal coefficient computed using $\calmf_{\psi_l}$ with different
$N_{q_l}^r$, $l = 1,2$, is shown in Fig.~\ref{fig:psi-coeff-mod}. We can see that no coefficient blows up and the reconstruction/prediction given by the POD-LSTM ROM is reasonably
good, with the exception of $\beta_{2,2}$ for
$N_{\psi_2}^r=2$ during the predictive phase $[50,100]$. See Fig. \ref{fig:psi-coeff-mod} (e). 

\begin{figure}[htb!]
    \centering
    \begin{subfigure}[h]{0.15\textwidth}
         \centering
         \includegraphics[width=\textwidth]{figs/psi1/psi1Mean_true.png}
         \caption{\scriptsize{$\tildep_1$, FOM}}
     \end{subfigure}
     \begin{subfigure}[h]{0.15\textwidth}
         \centering
         \includegraphics[width=\textwidth]{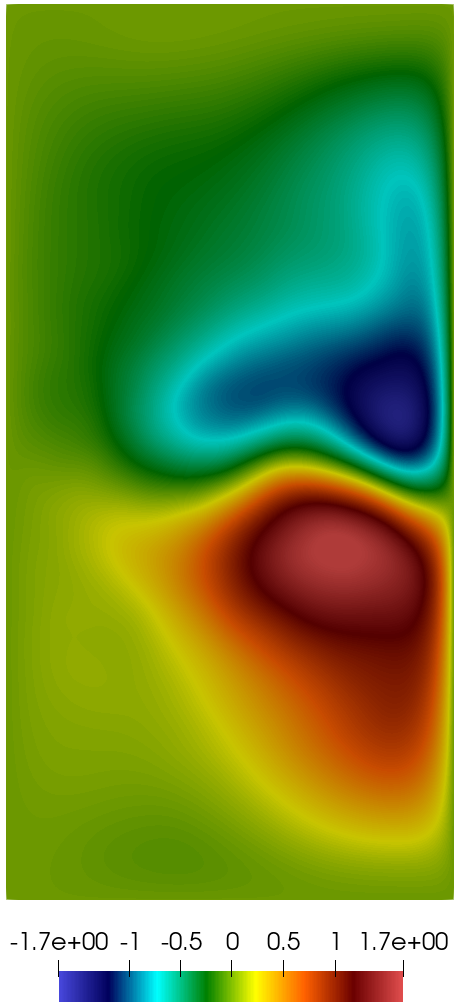}
         \caption{\scriptsize{$\tildep_1, \calmf_{q_1}$}}
     \end{subfigure}
    \begin{subfigure}[h]{0.15\textwidth}
         \centering
         \includegraphics[width=\textwidth]{figs/psi1/psi1Mean_mu10.png}
         \caption{\scriptsize{$\tildep_1, \calmf_{\psi_1}$}}
     \end{subfigure}\\
     \begin{subfigure}[h]{0.15\textwidth}
         \centering
         \includegraphics[width=\textwidth]{figs/psi2/psi2Mean_true.png}
         \caption{\scriptsize{$\tildep_2$, FOM}}
     \end{subfigure}
     \begin{subfigure}[h]{0.15\textwidth}
         \centering
         \includegraphics[width=\textwidth]{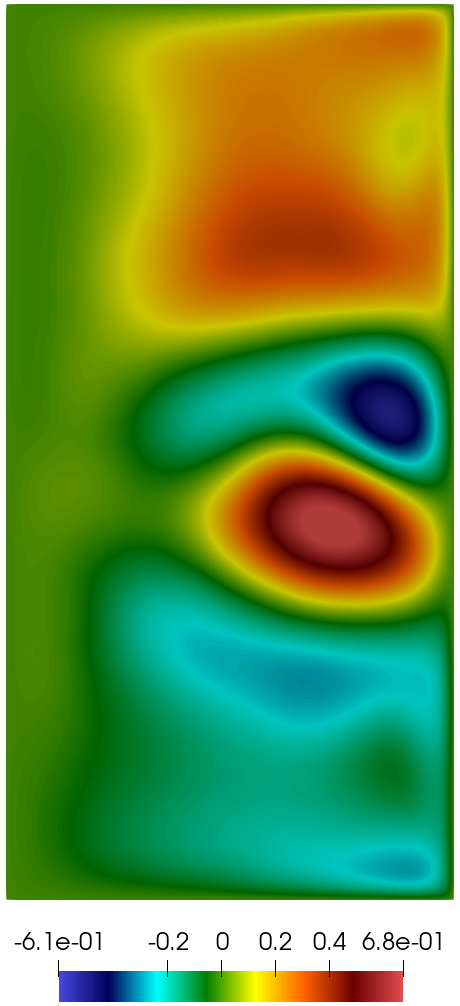}
         \caption{\scriptsize{$\tildep_2, \calmf_{q_2}$}}
     \end{subfigure}
    \begin{subfigure}[h]{0.15\textwidth}
         \centering
         \includegraphics[width=\textwidth]{figs/psi2/psi2Mean_mu10.png}
         \caption{\scriptsize{$\tildep_2, \calmf_{\psi_2}$}}
     \end{subfigure}
\caption{$\widetilde{\psi}_1$ (top) and $\widetilde{\psi}_2$ (bottom) computed over the predictive time interval $[50,100]$ by the FOM (first column) and the POD-LSTM ROM with modal coefficients from $\calmf_{q_l}$ (second column) and $\calmf_{\psi_l}$ (third column).}
\label{fig:psi-q}
\end{figure}

To conclude, we present results obtained with the alternative strategy suggested in Remark \ref{rem1}.
We fix the number of modes to $N_{q_1}^r = N_{q_2}^r = 10$ 
(which we recall implies $N_{\psi_1}^r = N_{\psi_2}^r = 10$) 
and the lookback window to $\sigma_L=10$.
Fig.~\ref{fig:psi-q} compares the time-averaged fields $\tildep_l$ computed by $\calmf_{q_l}$ and with those computed by $\calmf_{\psi_l}$ for $l = 1,2 $. We observe that the predictions of $\tildep_1$ and $\tildep_2$ given by $\calmf_{q_1}$ and $\calmf_{q_2}$ capture the general shapes of the gyres, however $\tildep_1$ and $\tildep_2$ computed by
$\calmf_{\psi_1}$ and $\calmf_{\psi_2}$ 
compare more favorably with the FOM counterparts, 
especially in the case of the gyres in the northern part of the basin. To quantitate the
agreement in Fig.~\ref{fig:psi-q}, we list in 
Table \ref{tab:error-psi-q}
error metrics $\varepsilon_{\psi_l}^{(1)}$ \eqref{eq:rmse} and $\varepsilon_{\psi_l}^{(2)}$ \eqref{eq:l2-error} for the stream functions.
The values in Table \ref{tab:error-psi-q} 
confirm that both $\calmf_{q_l}$ and $\calmf_{\psi_l}$ provide an accurate
reconstruction of $\tildep_l$, with the solutions 
given by $\calmf_{\psi_l}$ slightly more accurate. Table \ref{tab:error-psi-q} 
reports also the relative $L^2$ error for the kinetic energy \eqref{eq:kin-energy}.
We see that such errors are 3 times larger when
using $\calmf_{q_l}$, indicating that 
the time-dependent behavior is better captured by
$\calmf_{\psi_l}$.
We speculate that this is due to the fact that the modal coefficients computed by $\calmf_{q_l}$, $l = 1, 2$, are an order smaller than those given by $\calmf_{\psi_l}$. Since these models capture the fluctuations from the mean, $\calmf_{q_l}$ does well in predicting 
time-averaged fields but not time-dependent quantities because it predicts smaller fluctuations from the time-averaged fields. 

\begin{table}[h!]
    \centering
    \begin{tabular}{|cc|c|c|c|}
    \hline
    \multicolumn{2}{|c|}{Model} & $\varepsilon_{\psi_l}^{(1)}$ & $\varepsilon_{\psi_l}^{(2)}$ & $E_l$ rel. error \\ \hline
    \multicolumn{1}{|c|}{\multirow{2}{*}{$\psi_1$}} & $\calmf_{q_1}$ & 7.838E-02 & 1.360E-01 & 9.999E-01 \\ \cline{2-5} 
    \multicolumn{1}{|c|}{} & $\calmf_{\psi_1}$ & 6.041E-02 & 1.048E-01 & 3.642E-01 \\ \hline
    \multicolumn{1}{|c|}{\multirow{2}{*}{$\psi_2$}} & $\calmf_{q_2}$ & 4.393E-02 & 2.597E-01 & 9.524E-01 \\ \cline{2-5} 
    \multicolumn{1}{|c|}{} & $\calmf_{\psi_1}$ & 4.543E-02 & 2.479E-01 & 2.258E-01 \\ \hline
\end{tabular}
\caption{Error metrics for time-averaged fields $\tildep_1$ and $\tildep_2$, and kinetic energies $E_1$ and $E_2$ computed by POD-LSTM models $\calmf_{q_l}$ and $\calmf_{\psi_l}$ with respect to the FOM computed $\tildep_l$ and $E_l$.}
\label{tab:error-psi-q}
\end{table}

All the simulations presented in this section
were performed on a computing server running a 64-bit version of Linux with a Xeon\textsuperscript{\textregistered} processor and 128GB of RAM. A FOM simulation with the parameters specified in \eqref{eq:params1} over the time interval $[10,100]$ takes a total of $8.5$ days to complete. The online phase of our POD-LSTM ROM with $N_{q_l}^r = N_{\psi_l}^r = 10$ and $\sigma_L=10$ takes 0.5 s for $q_1$, 0.7 s for $q_2$, 0.9 s for $\psi_1$ and $\psi_2$. 
This is more than 1E+07 computational speed up compared to the FOM simulation. However, for the offline phase of the ROM, we need to collect the fluctuation field snapshots over the time interval $[10,50]$ from the FOM and this takes around 3.7 days. 
Additionally, the training of the LSTM networks $\calmf_{\Phi}$ takes 58 s for $q_1$, 59 s for $q_2$, 61 s for $\psi_1$, and 55 s for $\psi_2$. Both the snapshot collection and the LSTM training only need to be performed once. 
Lastly, we emphasize that the training of the networks and the online phase for each variable are independent of each other. This means that it is possible to compute the modal coefficients for the variables in parallel.

\subsection{Parametric study} \label{sec:param-study}

While the case of time as the only parameter is interesting in its own right and has allowed us to learn about the advantages and limitations of the POD-LSTM ROM, 
it is in a parametric study (i.e., a study where physical parameters are varied) that 
a ROM can show all of its potential. For this reason, we will consider $\delta$ as a variable parameter. All the other parameters are fixed and set as in \eqref{eq:params1}.

\begin{figure}[htb!]
    \centering
    \includegraphics[width = 0.45\textwidth]{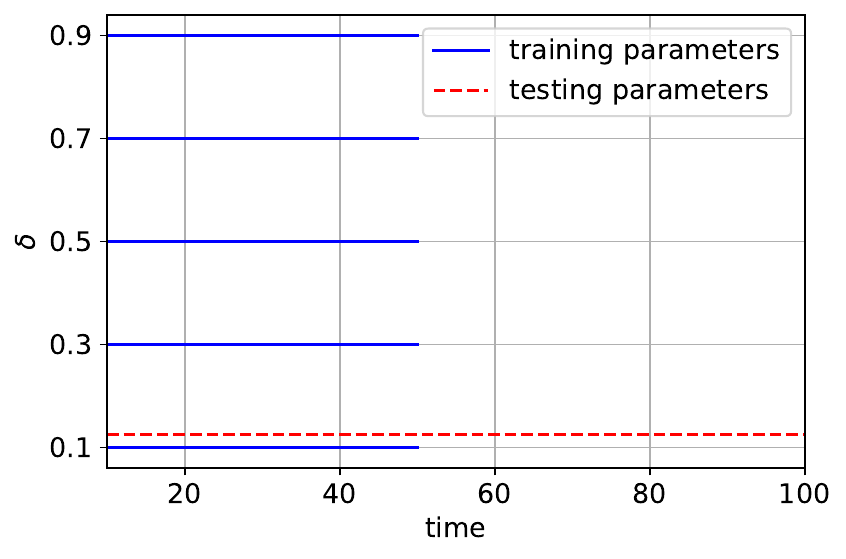}
    \caption{Parametric study: Training and testing parameters to assess the predictive capabilities of the POD-LSTM ROM.}
    \label{fig:train-test-param}
\end{figure}

We choose to train the POD-LSTM ROM with the following values: $\delta = 0.1, 0.3, 0.5, 0.7, 0.9$. Note that $\delta = 0.5$ is case
considered in Sec.~\ref{sec:asses-pod}. To test the prediction capability of the POD-LSTM ROM, we test it for $\delta = 0.125$ and
predict the time-averaged fields $\tildeq_l, \tildep_l$ over the time interval $[50,100]$, kinetic energies $E_l$ \eqref{eq:kin-energy} and enstrophies $\mathcal{E}_l$ \eqref{eq:enstrophy}, for both layers $l$. To build the POD basis, we collect 401 fluctuation snapshots for each training parameter set over the time interval $[10,50]$ for a total of 2005 fluctuation snapshots for each field of interest. See Fig.~\ref{fig:train-test-param} for the visualization of training and testing data used in this section.
For the LSTM network,
we keep the lookback window constant ($\sigma_L = 10$)  and use the same hyperparameters specified in Table~\ref{tab:hyperparam}.

\begin{figure}[htb!]
\centering
\begin{tabular}{cccccc}
       & FOM & \hspace{-0.4cm} $N_{q_l}^r = 2$ & \hspace{-0.4cm} $N_{q_l}^r = 4$ & \hspace{-0.4cm} $N_{q_l}^r = 8$ & \hspace{-0.4cm} $N_{q_l}^r = 10$ \\
    $\tildeq_1$ & \includegraphics[align=c,scale = 0.3]{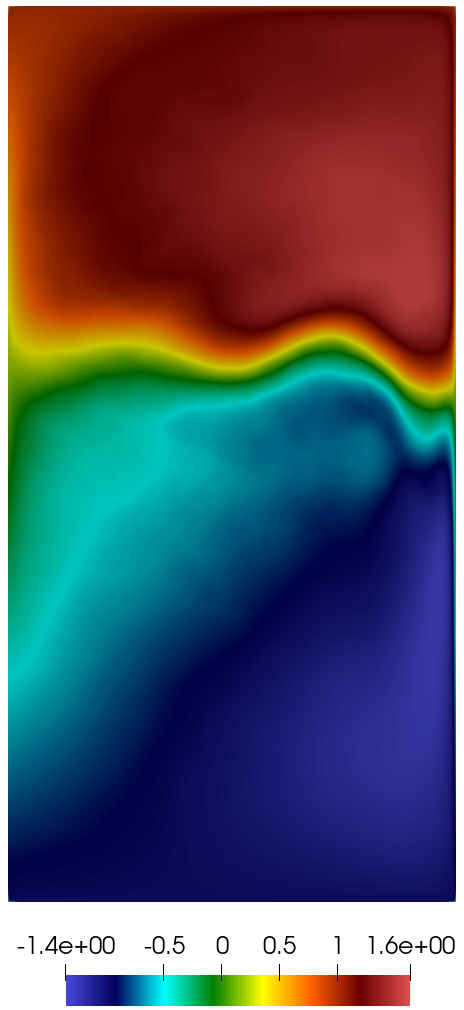} & \hspace{-0.4cm}\includegraphics[align=c,scale = 0.3]{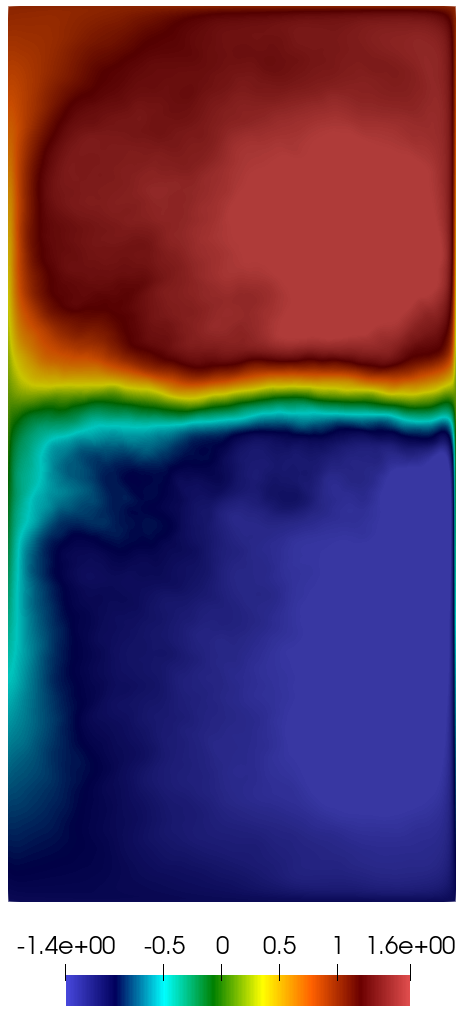} & \hspace{-0.4cm}\includegraphics[align=c,scale = 0.3]{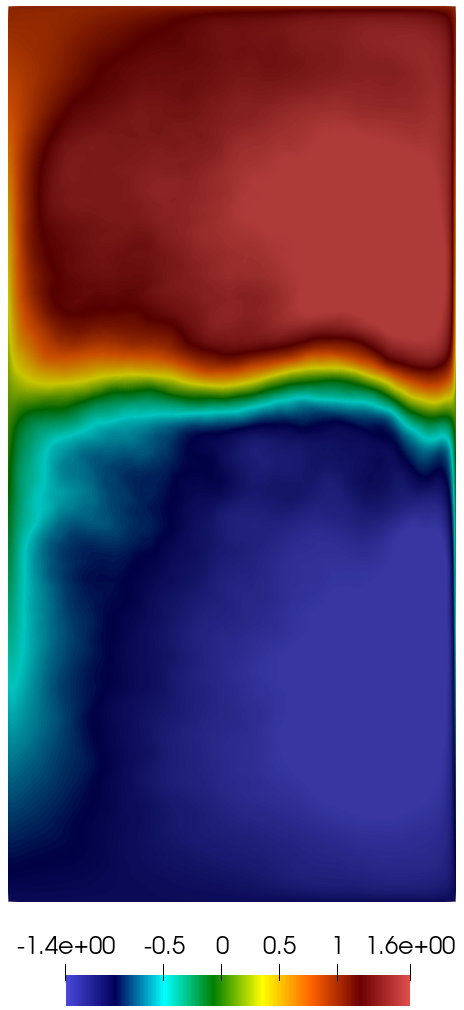} & \hspace{-0.4cm}\includegraphics[align=c,scale = 0.3]{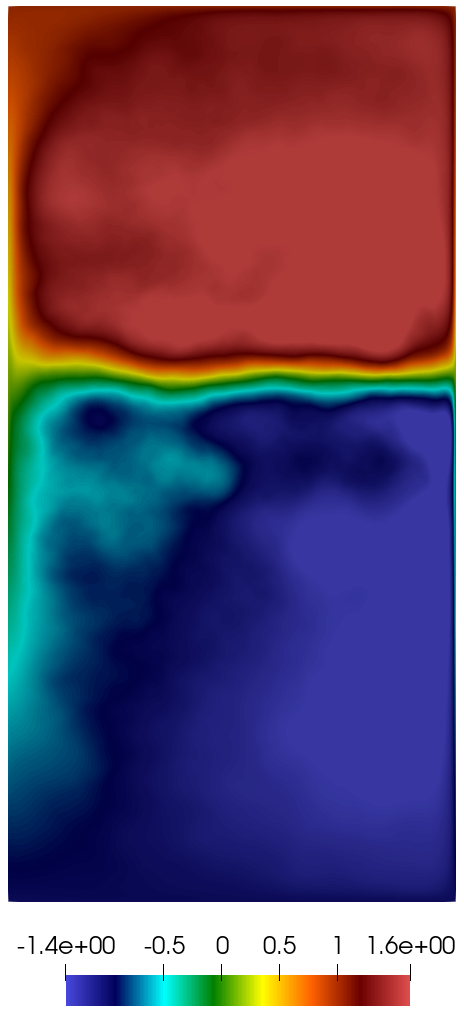} & \hspace{-0.4cm}\includegraphics[align=c,scale = 0.3]{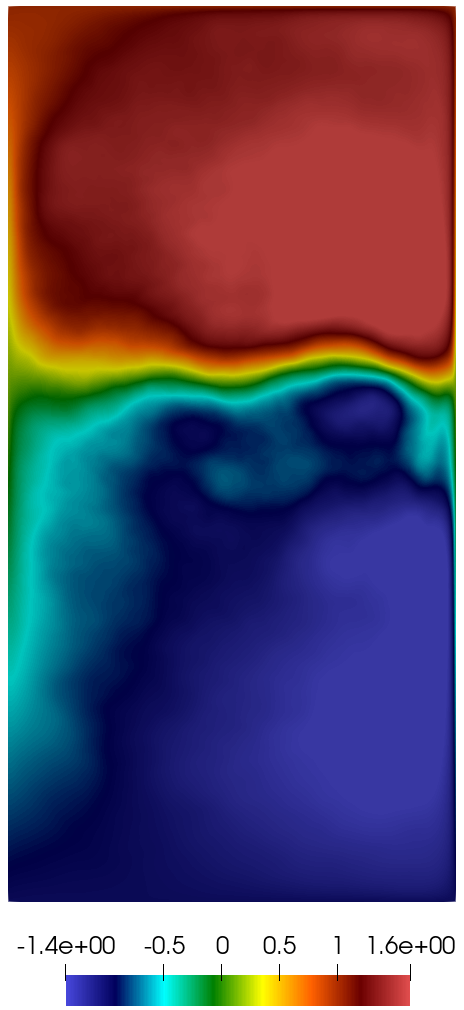} \\
    $\tildeq_2$ & \includegraphics[align=c,scale = 0.3]{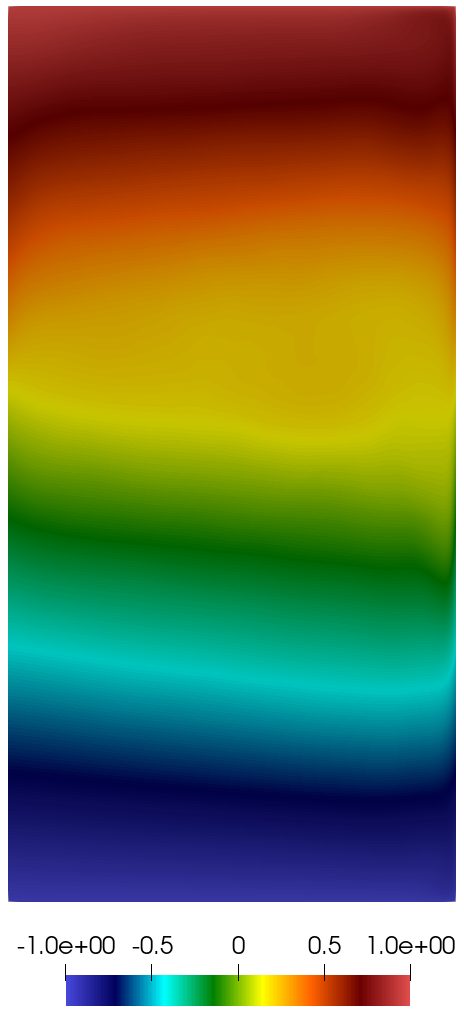} & \hspace{-0.4cm}\includegraphics[align=c,scale = 0.3]{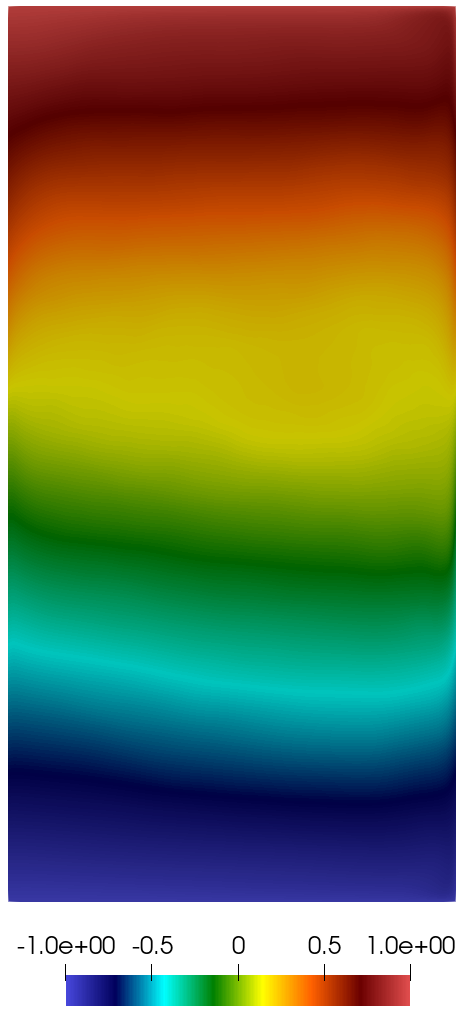} & \hspace{-0.4cm}\includegraphics[align=c,scale = 0.3]{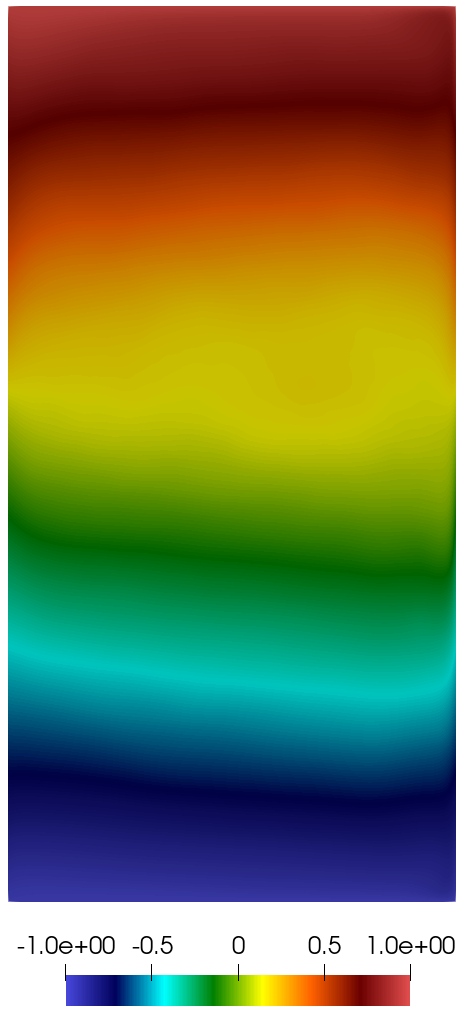} & \hspace{-0.4cm}\includegraphics[align=c,scale = 0.3]{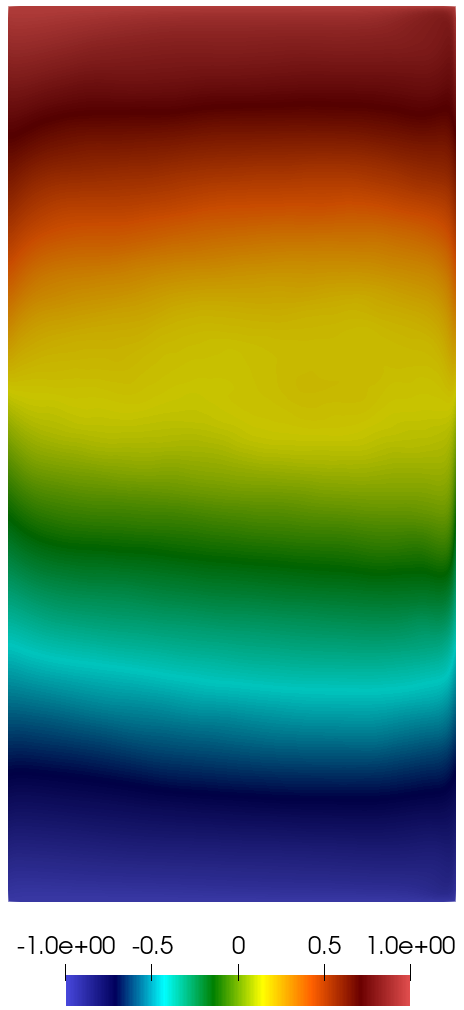} & \hspace{-0.4cm}\includegraphics[align=c,scale = 0.3]{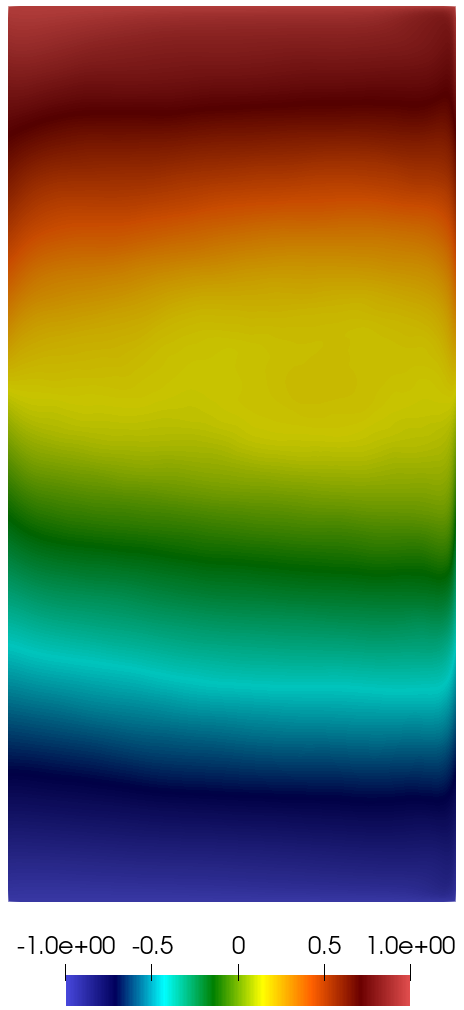} \\
\end{tabular}
\caption{Parametric study: Time-averaged vorticities $\widetilde{q}_1$ (top) and $\widetilde{q}_2$ (bottom) over the predictive time interval $[50,100]$ computed by FOM (first column) and POD-LSTM ROM with $N_{q_l}^r=2,4,8,10$ (second to fifth columns)}.
\label{fig:q_param}
\end{figure}

Let us start by investigating the effect of
$N_{\Phi}^r$ on the ROM prediction for the testing parameter $\delta = 0.125$. 
The closest sample to $\delta = 0.125$ is obviously $\delta = 0.1$, so we 
set ${\bm \mu}_c = 0.1$ in \eqref{eq:time_av_mu}. Fig.~\ref{fig:q_param}
compares $\tildeq_l$, $l = 1, 2$, computed by the FOM for $\delta = 0.125$
with its ROM counterpart for $N_{q_l}^r = 2,4,8,10$. We see that $\tildeq_2$
is rather well reconstructed by the ROM
for all $N_{q_2}^r$, while the quality of the reconstruction of $\tildeq_1$
varies as  $N_{q_1}^r$ changes.
Towards the basin center, $\tildeq_1$ exhibits a thin oscillating transition layer which is captured when $N_{q_1}^r=4$ and $N_{q_1}^r=10$, but not when $N_{q_1}^r=2$ or $N_{q_1}^r=8$.
We also observe that the small negative values (cyan region) is better captured as we increase $N_{q_1}^r$, with $N_{q_1}^r=10$ giving the solution that compares more
favorably with the FOM solution. 
This is also supported by the error metrics $\varepsilon_{q_1}^{(1)}$ \eqref{eq:rmse} and $\varepsilon_{q_1}^{(2)}$ \eqref{eq:l2-error}
reported in Table \ref{tab:error-q-param}.
Although $N_{q_2}^r=10$ does not provide the lowest error in either metric, all the values are comparable and it gives the smallest  
relative $L^2$ error for the 
enstrophy in the bottom layer. 

\begin{table}[h!]
    \centering
    \begin{tabular}{|c|cccc|}
    \hline
        $N_{q_1}^r$ & $\delta_{q_1}$ & $\varepsilon_{q_1}^{(1)}$ & $\varepsilon_{q_1}^{(2)}$ & $\mathcal{E}_1$ rel. err.\\
        \hline
        2 & 0.03 & 3.195E-01 & 3.099E-01 & 2.296E-01 \\
        \hline
        4 & 0.04 & 2.427E-01 & 2.354E-01 & 1.863E-01 \\
        \hline
        8 & 0.07 & 3.092E-01 & 2.998E-01 & 2.791E-01 \\
        \hline
        10 & 0.08 & 2.025E-01 & 1.965E-01 & 1.873E-01 \\
        \hline
         $N_{q_2}^r$ & $\delta_{q_2}$ & $\varepsilon_{q_2}^{(1)}$ & $\varepsilon_{q_2}^{(2)}$ & $\mathcal{E}_2$ rel. err.\\
        \hline
        2 & 0.04 & 1.437E-02 & 2.706E-02 & 3.421E-02 \\
        \hline
        4 & 0.06 & 1.421E-02 & 2.675E-02 & 2.932E-02 \\
        \hline
        8 & 0.09 & 1.496E-02 & 2.818E-02 & 2.819E-02 \\
        \hline
        10 & 0.10 & 1.524E-02 & 2.869E-02 & 2.599E-02 \\
        \hline
    \end{tabular}
    \caption{Parametric study: fraction of retained eigenvalue energy $\delta_{q_l}$, 
    error metrics $\varepsilon_{q_l}^{(1)}$ \eqref{eq:rmse} and $\varepsilon_{q_l}^{(2)}$ \eqref{eq:l2-error}, 
    and relative $L^2$ error for enstrophy $\mathcal{E}_l$ \eqref{eq:enstrophy},
    $l = 1,2 $, for different numbers of POD modes retained.}
    \label{tab:error-q-param}
\end{table}

\begin{figure}[htb!]
\centering
\begin{tabular}{cccccc}
       & FOM & \hspace{-0.4cm} $N_{\psi_l}^r = 2$ & \hspace{-0.4cm} $N_{\psi_l}^r = 4$ & \hspace{-0.4cm} $N_{\psi_l}^r = 8$ & \hspace{-0.4cm} $N_{\psi_l}^r = 10$ \\
    $\tildep_1$ & \includegraphics[align=c,scale = 0.3]{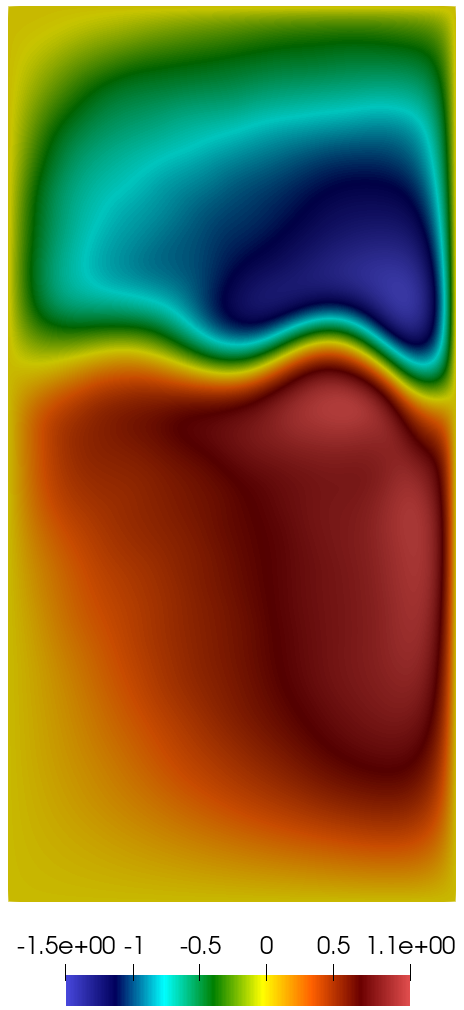} & \hspace{-0.4cm}\includegraphics[align=c,scale = 0.3]{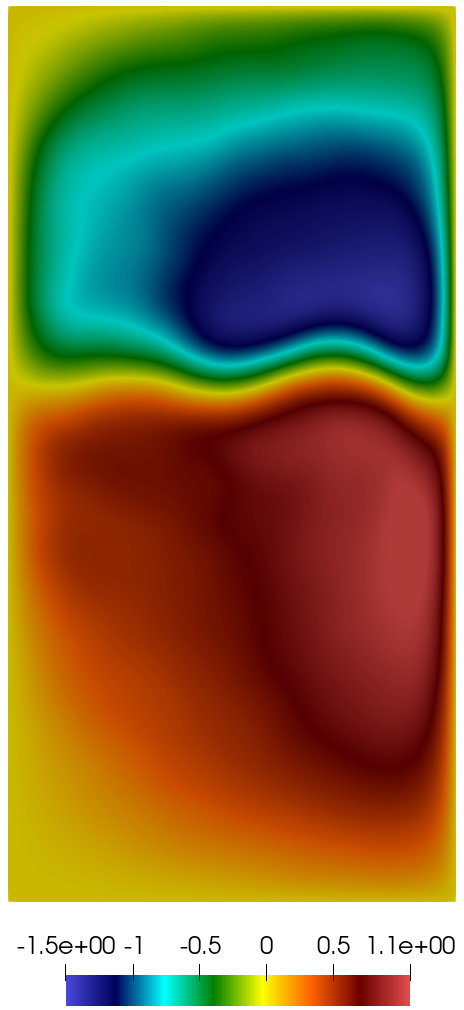} & \hspace{-0.4cm}\includegraphics[align=c,scale = 0.3]{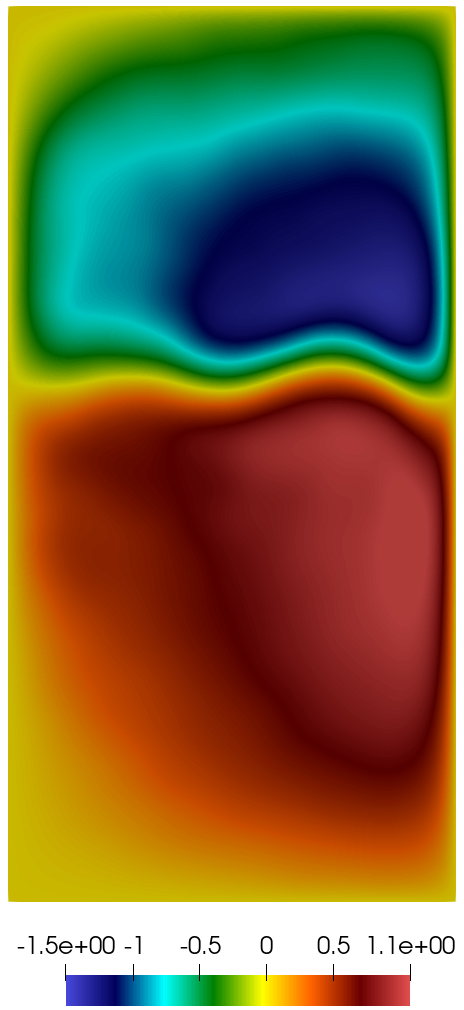} & \hspace{-0.4cm}\includegraphics[align=c,scale = 0.3]{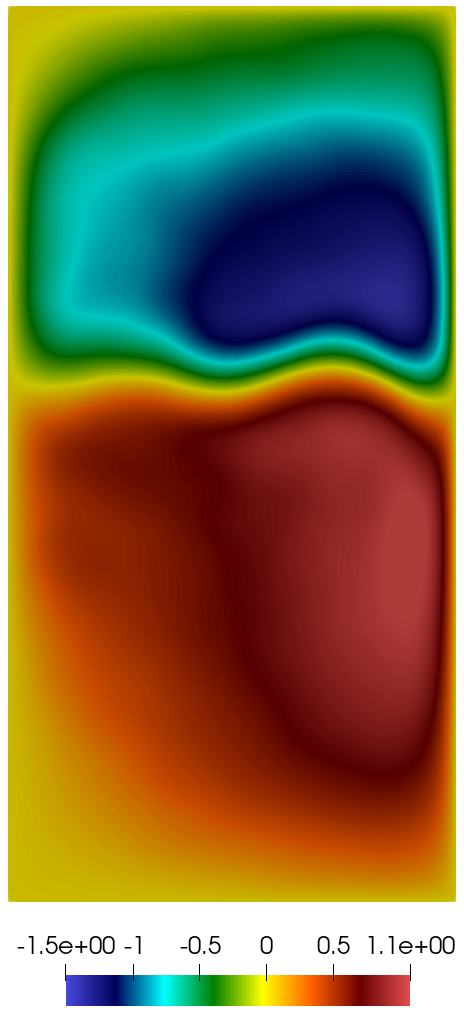} & \hspace{-0.4cm}\includegraphics[align=c,scale = 0.3]{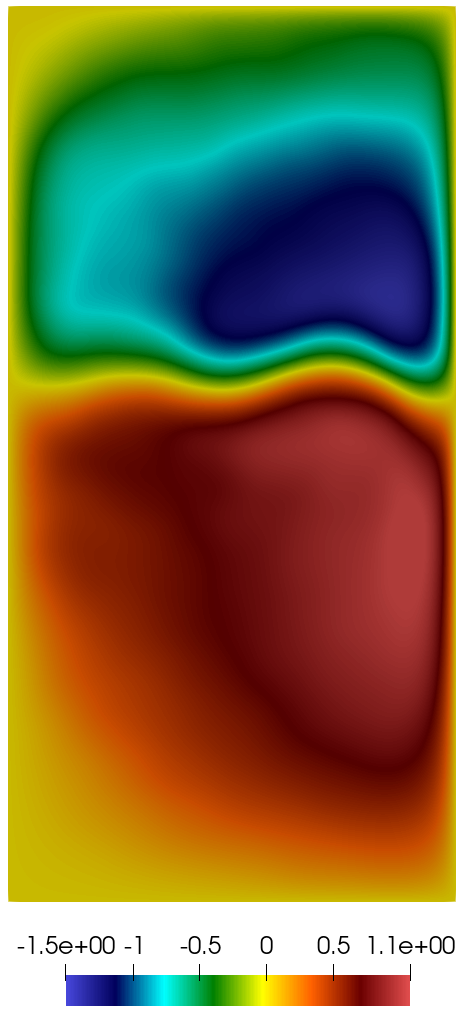} \\
    $\tildep_2$ & \includegraphics[align=c,scale = 0.3]{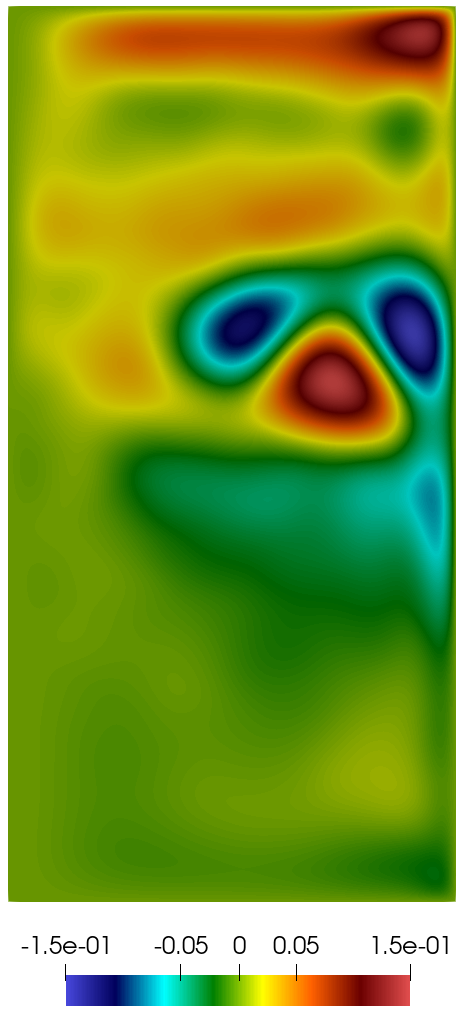} & \hspace{-0.4cm}\includegraphics[align=c,scale = 0.3]{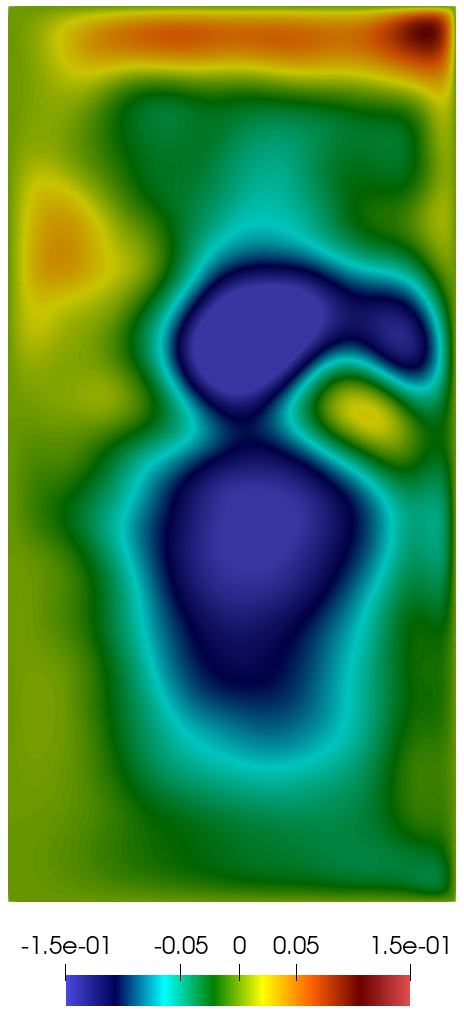} & \hspace{-0.4cm}\includegraphics[align=c,scale = 0.3]{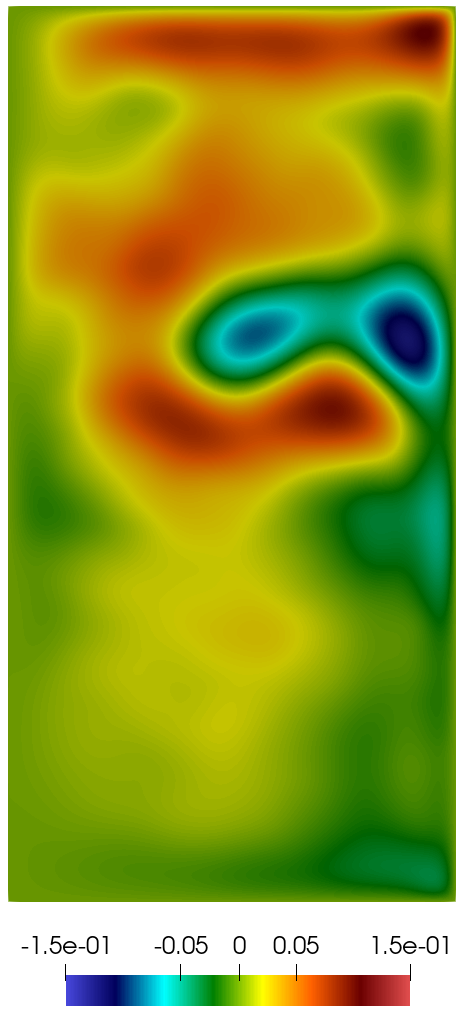} & \hspace{-0.4cm}\includegraphics[align=c,scale = 0.3]{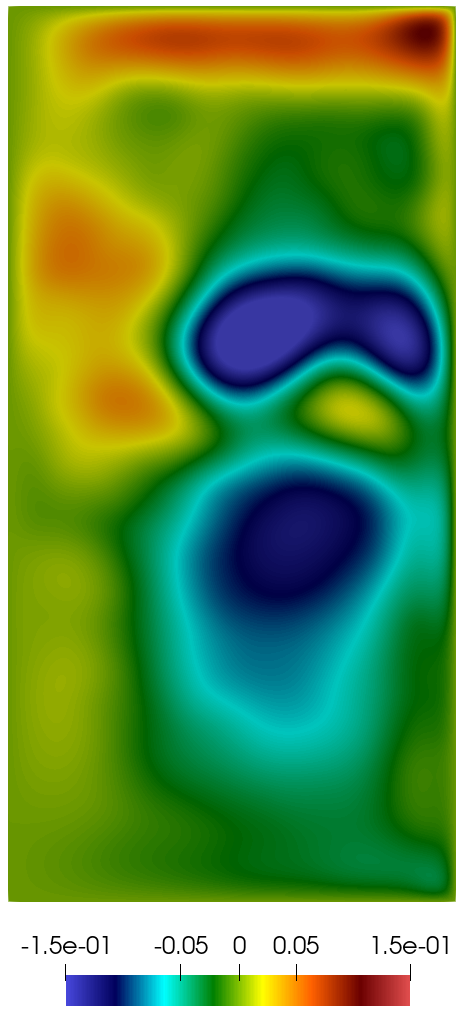} & \hspace{-0.4cm}\includegraphics[align=c,scale = 0.3]{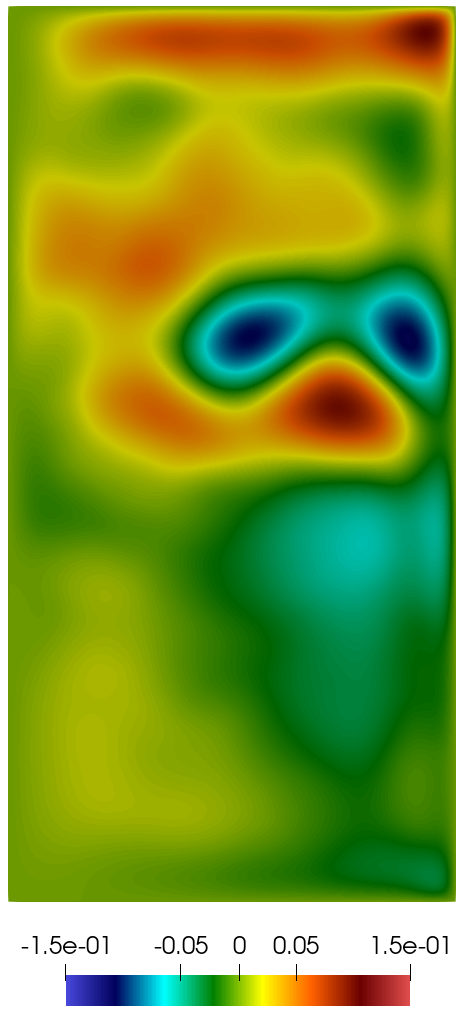} \\
\end{tabular}
\caption{Parametric study: Time-averaged stream functions $\widetilde{\psi}_1$ (top) and $\widetilde{\psi}_2$ (bottom) over the predictive time interval $[50,100]$ computed by FOM (first column) and POD-LSTM ROM with $N_{\psi_l}^r=2,4,8,10$ (second to fifth columns)}.
\label{fig:psi_param}
\end{figure}

Next, let us look at the reconstruction of 
$\tildep_l$, $l = 1, 2$, when $N_{\psi_l}^r$
varies in Fig.~\ref{fig:psi_param}. Note that in the top layer
there are two large and irregularly shaped gyes, while the large gyres are lost in the bottom layer.
The ROM predictions of $\tildep_1$ match well with the $\tildep_1$ given by the FOM for all $N_{\psi_1}^r$.
We see that $N_{\psi_1}^r = 2$ predicts the shape of the red gyre (large positive values) in the Southern region of the ocean basin most 
accurately but $N_{\psi_1}^r = 10$ is able to capture the blue gyre (large negative values) in the Northern region better than the rest of the cases considered. 
In fact, $N_{\psi_1}^r = 10$ results into the lowest error metrics $\varepsilon_{\psi_l}^{(1)}$ \eqref{eq:rmse} and $\varepsilon_{\psi_l}^{(2)}$ \eqref{eq:l2-error}, as reported 
in the third and fourth columns of Table \ref{tab:error-psi-param}.
The results for $\tildep_2$ differ significantly when $N_{\psi_2}^r$ is varied. 
We see that $N_{\psi_2}^r = 2$ and $N_{\psi_2}^r = 8$ cannot accurately capture the time-averaged
dynamics in the lower layer, especially
around the center of the basin and in the Southern region.
On the other hand, the ROM prediction improves when $N_{\psi_2}^r = 4$, and even more when $N_{\psi_2}^r = 10$. 
These observations are reflected in the root mean square error $\varepsilon_{\psi_2}^{(1)}$ \eqref{eq:rmse} and the $L^2$ relative error $\varepsilon_{\psi_2}^{(2)}$ \eqref{eq:l2-error} in the third and fourth columns of Table \ref{tab:error-psi-param}. Table \ref{tab:error-psi-param} also reports
the relative error for the predicted time-evolution of the kinetic energy $E_l$. We see that the errors are comparable across all values of $N_{\psi_l}^r$.

\begin{table}[h!]
    \centering
    \begin{tabular}{|c|cccc|}
    \hline
        $N_{\psi_1}^r$ & $\delta_{\psi_1}$ & $\varepsilon_{\psi_1}^{(1)}$ & $\varepsilon_{\psi_1}^{(2)}$ & $E_1$ rel. err.\\
        \hline
        2 & 0.31 & 1.398E-01 & 2.027E-01 & 7.151E-01 \\
        \hline
        4 & 0.40 & 1.336E-01 & 1.937E-01 & 7.129E-01 \\
        \hline
        8 & 0.50 & 1.339E-01 & 1.940E-01 & 7.171E-01 \\
        \hline
        10 & 0.54 & 1.282E-01 & 1.855E-01 & 7.189E-01 \\
        \hline
         $N_{\psi_2}^r$ & $\delta_{\psi_2}$ & $\varepsilon_{\psi_2}^{(1)}$ & $\varepsilon_{\psi_2}^{(2)}$ & $E_2$ rel. err.\\
        \hline
        2 & 0.37 & 5.923E-02 & 1.731E+00 & 9.561E-01 \\
        \hline
        4 & 0.49 & 2.561E-02 & 7.485E-01 & 9.643E-01 \\
        \hline
        8 & 0.61 & 4.253E-02 & 1.243E+00 & 9.571E-01 \\
        \hline
        10 & 0.65 & 1.887E-02 & 5.513E-01 & 9.647E-01 \\
        \hline
    \end{tabular}
    \caption{Parametric study: fraction of retained eigenvalue energy $\delta_{\psi_l}$, 
    error metrics $\varepsilon_{\psi_l}^{(1)}$ \eqref{eq:rmse} and $\varepsilon_{\psi_l}^{(2)}$ \eqref{eq:l2-error}, 
    and relative $L^2$ error for kinetic energy $E_l$ \eqref{eq:enstrophy},
    $l = 1,2 $, for different numbers of POD modes retained.}
    \label{tab:error-psi-param}
\end{table}

\begin{figure}[htb!]
    \centering
    \begin{subfigure}[h]{0.15\textwidth}
         \centering
         \includegraphics[width=\textwidth]{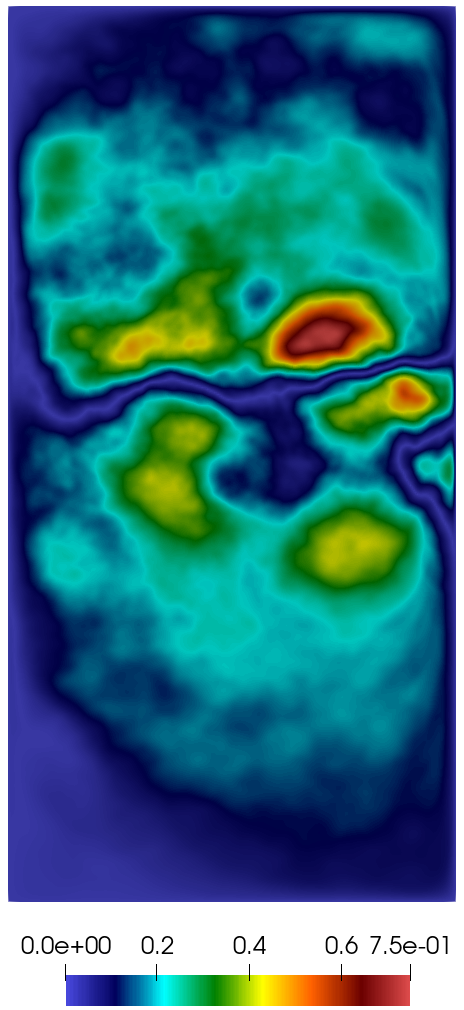}
         \caption{\tiny{$\abs{\tildeq_1^\text{FOM} - \tildeq_1^\text{ROM}}$}}
     \end{subfigure}
     \begin{subfigure}[h]{0.15\textwidth}
         \centering
         \includegraphics[width=\textwidth]{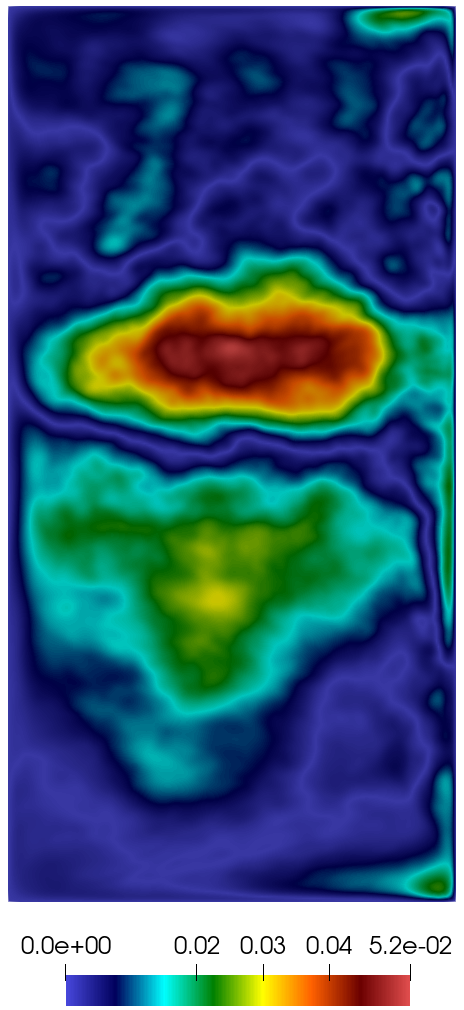}
         \caption{\tiny{$\abs{\tildeq_2^\text{FOM} - \tildeq_2^\text{ROM}}$}}
     \end{subfigure}
    \begin{subfigure}[h]{0.15\textwidth}
         \centering
         \includegraphics[width=\textwidth]{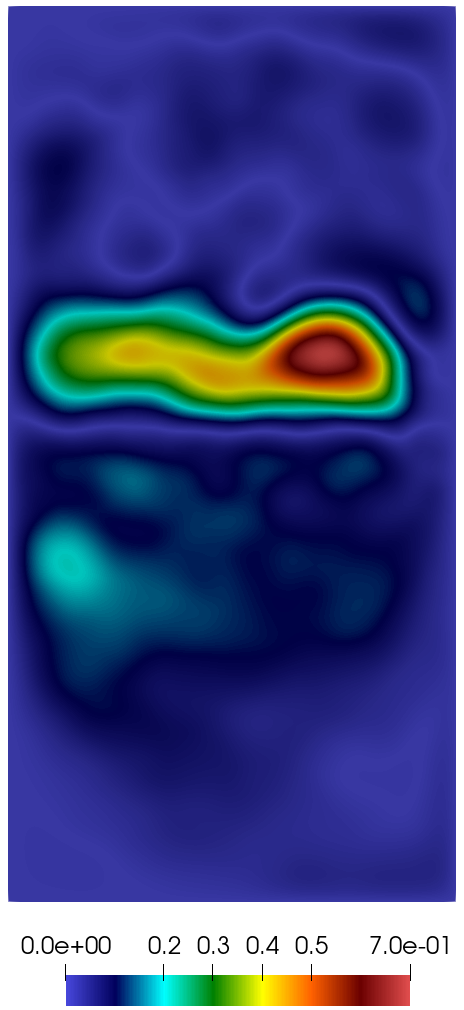}
         \caption{\tiny{$\abs{\tildep_1^\text{FOM} - \tildep_1^\text{ROM}}$}}
     \end{subfigure}
     \begin{subfigure}[h]{0.15\textwidth}
         \centering
         \includegraphics[width=\textwidth]{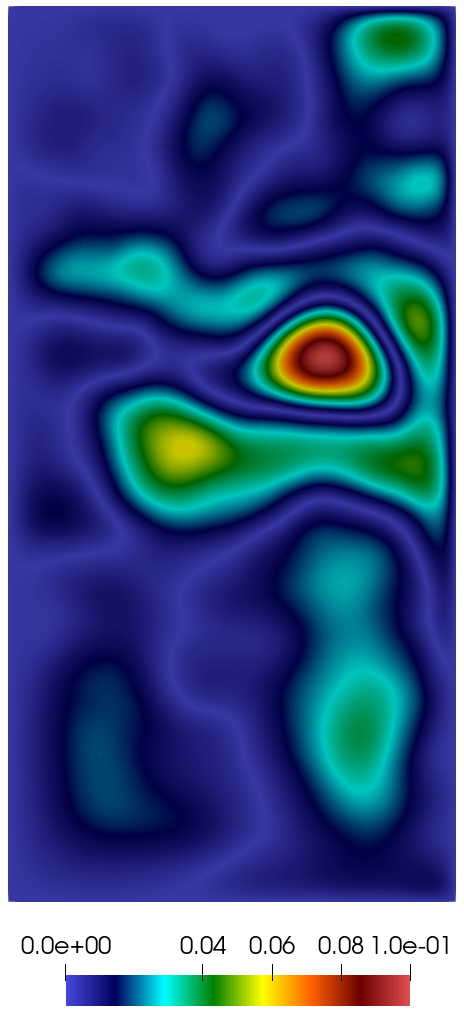}
         \caption{\tiny{$\abs{\tildep_2^\text{FOM} - \tildep_2^\text{ROM}}$}}
     \end{subfigure}
\caption{Parametric study: Absolute difference between the FOM time-averaged solutions and the predicted time-averaged fields computed by the POD-LSTM ROM with $N_{\Phi}^r = 10$ and $\sigma_L = 10$.}
\label{fig:abs_err}
\end{figure}

Next, we show the absolute difference between the true time-averaged fields coming from the FOM and the POD-LSTM ROM with $N_{\Phi}^r = 10$ and $\sigma_L = 10$ in Fig.~\ref{fig:abs_err}. 
We see that the peaks in the absolute difference are localized.
Notice that although it might seem that for $\tildeq_2$, there is a large red region, i.e., large absolute difference, the maximum absolute difference is of order 1E-02.

\begin{rem}
In the case of $\delta = 0.125$, it was obvious to take ${\bm \mu}_c = 0.1$ in \eqref{eq:time_av_mu}. However, if one is interested in 
the ROM approximation for $\delta = 0.2$ with the sampling in Fig.~\ref{fig:train-test-param}, the choice between ${\bm \mu}_c = 0.1$ and ${\bm \mu}_c = 0.3$ is less obvious. We found out that the time-averaged FOM fields for $\delta = 0.2$ are much closer to the time-averaged fields for $\delta = 0.1$ than $\delta = 0.3$ and, hence, 
${\bm \mu}_c = 0.1$ yields a more accurate ROM approximation for $\delta = 0.2$. 
A more robust ROM approximation could be obtained by replacing 
$\tildeq_l^0$ and $\tildep_l^0$ in \eqref{eq:q_approx_mu}-\eqref{eq:psi_approx_mu} with higher
order approximations of $\tildeq_l(\xbf,\bm{\mu})$ and 
$\tildep_l(\xbf, \bm{\mu})$.
\end{rem}

To conclude this section, we comment on the computational time savings. 
All simulations were carried out using the same machine described 
in Sec.~\ref{sec:asses-pod}.
For the offline phase of the parametric study, we need to collect the fluctuation field snapshots over the time interval $[10,50]$ from the FOM for five 
values of $\delta$. Each of such FOM simulations takes around 3.7 days,
so the offline phase is significantly more expensive if not run in parallel. 
The corresponding ROM computations online, with $N_{q_l}^r = N_{\psi_l}^r = 10$ 
modes retained and lookback window of $\sigma_L = 10$,
take  0.51 s for both $q_1$ and $q_2$, 0.52 s for $\psi_1$, and 0.55 s for $\psi_2$.

\section{Conclusions} \label{sec:conclusion}

We proposed a POD-based data-driven ROM to capture the large-scale behavior of ocean dynamics described by 
the two-layer quasi-geostrophic equations. The ROM approximates the 
solution of the 2QGE with a linear combination of POD
modes, with modal coefficients found through a recurrent neural network called long short-term memory. The LSTM architecture
was chosen for its robustness and stability in predicting chaotic dynamical systems. The proposed methodology consists of two phases: offline and online. During the offline phase, snapshots data are collected from a
high-resolution simulation for a given time interval
(and given parameter values, in the case of variable physical parameters), the reduced basis for each field variable is generated via POD, and the LSTM architecture is trained with the time-dependent modal coefficients associated with the snapshots.  A critical model parameter to set in this phase is the lookback window for the LSTM, i.e., the quantity that tells the network how far in the past to look.
During the online phase, the trained model is used to predict the modal coefficients recursively for the remaining of the time interval of interest (and new parameter values, in the case of variable physical parameters). 

We assessed the performance of the POD-LSTM ROM with
an extension of a widely used benchmark called the double-gyre wind forcing experiment. We started by considering time as the only parameter and 
then included the the aspect ratio $\delta$ between the two ocean layers as a variable parameter. 
Our results indicate that, even when retaining only $10-20\%$ of the singular value energy of the system, our POD-LSTM ROM provides an accurate prediction of both time-averaged fields and time-dependent quantities (modal coefficients, enstrophy, and kinetic energy). We observe up to 1E+07 computational speed up in the online phase compared to the FOM.
A drawback of this approach is its sensitivity 
to the lookback window. Our computational study indicates there is an ``optimal value'', past
which the prediction by the POD-LSTM ROM degenerates.
Once this sensitivity is better understood, 
the POD-LSTM ROM will be a methodology, not only accurate and efficient, but also robust for the 
prediction of large-scale flows with choatic spatio-temporal behavior. 


\section*{Acknowledgements}
We acknowledge the support provided by  
PRIN “FaReX - Full and Reduced order modelling of coupled systems: focus on non-matching methods and automatic learning” project, PNRR NGE iNEST “Interconnected Nord-Est Innovation Ecosystem” project, INdAM-GNCS 2019–2020 projects and PON “Research and Innovation on Green related issues” FSE REACT-EU 2021 project. This work was also partially supported by the U.S. National Science Foundation through Grant No. DMS-1953535
(PI A. Quaini).

\bibliography{QGE}
\end{document}